\documentclass{article} 

\usepackage[english]{babel} 
\usepackage{amssymb}
\usepackage{amsmath}
\usepackage{txfonts}
\usepackage{mathdots}
\usepackage[classicReIm]{kpfonts}
\usepackage[dvips]{graphicx} 

\begin{document}


\begin{center}\textbf{Biorthogonal Polynomial System}\end{center}

\begin{center}\textbf{Composed of X-Jacobi Polynomials from Different Sequences}\end{center}

\begin{center} Gregory Natanson\end{center}

\begin{center} ai-solutions Inc.\end{center}

\begin{center} greg\_natanson@yahoo.com\end{center}

\noindent The paper examines rational Darboux transformations (RDTs) of the Jacobi equation written in the canonical form, with emphasis on the Sturm-Liouville problems (SLPs)  solved under the Dirichlet boundary conditions (DBCs) at the ends of the infinite interval [1, \textit{inf}).  To be able to extend the analysis to the Darboux-Crum net of rational SL equations (SLEs) solved under the cited DBCs in terms of multi-indexed \textit{orthogonal} exceptional Romanovski-Jacobi (XR-Jacobi) polynomials, we consider only seed functions which represent \textit{principal }Frobenius solutions (PFSs) near one of the singular endpoints. There are three distinct types of  such solutions with no zeros inside the selected interval: two infinite sequences formed by the PFSs near the lower endpoint and one finite sequence formed by the PFSs near infinity.   It is shown that use of classical Jacobi polynomials as seed functions results in the double-indexed manifold composed of orthogonal X${}_{m}$-Jacobi polynomials in the \textit{reversed} argument. As a result polynomials from this manifold obey the cross-orthogonality relation when integrated from     +1 to \textit{inf}.  As a corollary we assert that each X${}_{m}$-Jacobi polynomial of degree m + n has exactly m \textit{exceptional} zeros between --\textit{inf} and --1 as far as its indexes are restricted by the derived constraints on indexes of XR-Jacobi polynomials.

\noindent \textbf{\eject }

\section*{1. Introduction}

\noindent Though the Schrödinger equation with the hypergeometric Pöschl-Teller (h-PT) potential [1] was a focus of numerous articles for several decades recent study of Odake and Sasaki [2] (see also Appendix D in their joint work with Ho [3]) revealed a remarkable \textit{novel} element completely overlooked in the earlier literature.  Namely, they discovered that this equation, when mapped onto the sub-domain [1, $\mathrm{\infty}$) of the hypergeometric equation by the change of variable z(x) = cosh(2x)  [4]  is quantized by a finite set of orthogonal polynomials which were later identified by us [5] as `Romanovski-Jacobi' polynomials [6-10].  (In following [11] we will refer to the latter simply as `R-Jacobi' polynomials.)  

In our original study [12] on rational potentials exactly solvable by hypergeometric functions the author considered only the change of variable z(x) which converted the Schrödinger equation onto the hypergeometric equation defined on the finite subdomain 0 $\mathrm{<}$ z $\mathrm{<}$ 1, i.e., the change of variable  z(x) = tanh(x) [13] in the h-PT potential.  It was taken for granted that there is no use in discussing alternative representations by mapping the Schrödinger equation onto infinite intervals -$\infty$ $\mathrm{<}$ z $\mathrm{<}$ 0 [14, 15] and 1 $\mathrm{<}$ z $\mathrm{<}$ +$\infty$ since both mappings could be realized via linear fractional transformations of the already constructed rational potential.  The aforementioned restriction of the change of variable z(x) automatically assured that the rational potentials in question were solvable by \textit{classical }Jacobi polynomials (with energy-dependent indexes in general).   However author's original view was changed a few years ago under influence of Odake and Susaki's discovery [2].

Recently that the author [16] developed the general technique for generating rational Darboux-Crum [17, 18] transforms (RDC$\mathscr{T}$s) of the canonical Sturm-Liouville equations (CSLEs) exactly solvable via hypergeometric and confluent hypergeometric functions [14, 15, 12] as well as of the partner Fuschian SLE with two second-order poles at $\pm$\textit{i}  [19, 20].  This allowed us to construct three DC nets of rational CSLEs (RCSLEs) solvable by polynomials which satisfy Heine-type equations [21-23] and for this reason termed by us [5, 24] `Gauss-seed' (GS) Heine polynomials. To be more precise, we refer to these polynomials (with \textit{variable} indexes in general) as Jacobi-seed ($\mathscr{J}$S), Laguerre-seed (LS), and Routh-seed (RS) Heine polynomials, when  Jacobi, generalized Laguerre, and Routh [25] polynomials are used accordingly [26, 27, 20] to construct the appropriate `quasi-rational' [28, 29] seed functions (or, to be more precise, \textit{seed functions with rational log-derivatives}).   If the density function has the simple poles at both singular points of the Fuschian equations with two second-order poles on the either real or imaginary axis then the exponent differences (ExpDiffs) for the common finite persistent singularities of the corresponding RDC net (RDCN) of solvable real RSLEs become energy independent and we come to two RDCNs of Bochner-type [30] eigenequations with polynomial coefficient functions.  Each equation has infinitely many polynomial solutions referred to by us as exceptional Bochner-Jacobi (XB-Jacobi) or exceptional Bochner-Routh (XB-Routh) polynomials accordingly, in following the terminology suggested by Gomez-Ullate, Kamran, and Milson  [31] for orthogonal polynomial systems (OPSs) starting from polynomials of non-zero degree and thereby violating the Bochner theorem [30].

Taking into account that gauge transformation of the RCSLE to an arbitrarily chosen self-adjoint form shifts the linear coefficient of the indicial equation for the second-order pole at the given endpoint, we formulate Sturm-Liouville problem (SLP) of our interest by imposing the Dirichlet boundary conditions (DBCs) on solutions of the so-called `prime' [16] self-adjoint SLE.  This self-adjoint realization of the given RCSLE is chosen in such way that sum of two characteristic exponents (ChExps) is equal to zero at each endpoint pole.  As a result of such a very specific choice of the self-adjoint SLE the DBC unambiguously selects the principal [32-35] Frobenius solution (PFS) for any real (and therefore necessarily positive) value of the corresponding exponent difference (ExpDiff) whether it is the limit-point (LP) or limit-circle (LC) case [36-39]. 

As far as the rational Darboux transforms (RD$\mathscr{T}$s) of the quasi-rational eigenfunctions of the prime SLE satisfy the DBCs at both singular endpoints they must be mutually orthogonal with the weight function unambiguously prescribed by the particular selection of the self-adjoint form from its infinitely many choices.  If the density function has simple poles  at both  persistent singular points in the finite plane then the indexes of GS-Heine polynomials (coincident with the ExpDiffs of the Gauss-reference (GRef) CSLE starting the RDCN in question) become energy independent and as a result the polynomials themselves form orthogonal sets [5, 40, 41], when integrated with the properly chosen weight function. 

For the finite interval [-1, +1] the corresponding prime self-adjoint RSLE was first introduced in [31] to prove the orthogonality of X${}_{1}$-Jacobi polynomials. However the appropriate Jacobi-reference (JRef) CSLE can be alternatively converted to the prime RSLE solved with the DBCs at the ends of the infinite interval [+1, +$\mathrm{\infty}$).   When the new self-adjoint SLE is solved under the DBCs at the singular endpoints of the infinite interval, the resultant eigenfunctions become expressible in terms of R-Jacobi polynomials with \textit{energy-independent} indexes. Since any PFS below the lowest eigenvalue is necessarily nodeless it can be used as a seed functions for the RDCN of SLPs conditionally exactly solvable via multi-index exceptional R-Jacobi (XR-Jacobi) polynomials.

It is shown that there are three distinct types of Jacobi polynomials with no zeros inside the selected interval (R-Jacobi admissible in Duran's terms [42]): two infinite sequences formed by the PFSs near the lower endpoint and one finite sequence formed by the PFSs near infinity.   It is shown that use of classical Jacobi polynomials as seed functions results in the double-indexed manifold composed of orthogonal X${}_{m}$-Jacobi polynomials in the \textit{reversed} argument. As a result polynomials from this manifold obey the cross-orthogonality relation when integrated from +1 to +$\mathrm{\infty}$.   As a corollary we assert that each X${}_{m}$-Jacobi polynomial of degree m + n has exactly m \textit{exceptional} zeros between --$\mathrm{\infty}$ and --1 provided that its indexes are restricted by the derived constraints on indexes of XR-Jacobi polynomials.

\section*{2. Double-Indexed XB-Jacobi eigenpolynomials}

\noindent As outlined in a more general context in [16], we define the monic XB-Jacobi eigenpolynomials 
\[{\rm P}_{{\rm N}_{{\rm \{ }\stackrel{\rightharpoonup}{{\bf \sigma }}{\rm m}\} _{2} } }^{{\rm (}\stackrel{\rightharpoonup}{{\rm \lambda }}_{o} )} [\eta |\stackrel{\rightharpoonup}{{\rm \sigma }}_{1} {\rm ,m}_{1} {\rm ;}\stackrel{\rightharpoonup}{{\rm \sigma }}_{2} {\rm ,m}_{2} ]\equiv \Pi _{{\rm N}_{{\rm \{ }\stackrel{\rightharpoonup}{{\rm \sigma }}{\rm m}\} _{2} } } [\eta ;\bar{\eta }_{{\rm \{ }\stackrel{\rightharpoonup}{{\rm \sigma }}m\} _{2} } (\stackrel{\rightharpoonup}{{\rm \lambda }}_{o} )],   (2.1)\] 
with ${\rm N}_{{\rm \{ }\stackrel{\rightharpoonup}{{\bf \sigma }}{\rm m}\} _{2} } $simple roots $\bar{\eta }_{{\bf \{ }\stackrel{\rightharpoonup}{{\bf \sigma }}m\} _{2} } $and $\stackrel{\rightharpoonup}{{\bf \lambda }}_{o} \equiv \lambda _{o;-} ,\lambda _{o;+} \, \, \, (\lambda _{o;\pm } >0)$, via the monomial decomposition of the so-called `polynomial determinants' (PDs) 
\[{\rm D}_{n+m+1} [\eta ;\stackrel{\rightharpoonup}{{\bf \lambda }}|\stackrel{\rightharpoonup}{{\bf \sigma }}{\bf ,}{\rm m}]\equiv \left|\mathop{}\limits_{S_{m+1}^{(\sigma _{+} \lambda _{+} ,\sigma _{-} \lambda _{-} )} (\eta )\, \, \, \, \, \, S_{n+1}^{(\lambda _{+} ,\lambda _{-} )} (\eta )\, \, }^{P_{m}^{(\sigma _{+} \lambda _{+} ,\sigma _{-} \lambda _{-} )} (\eta )\, \, \, \, \, \, P_{n}^{(\lambda _{+} ,\lambda _{-} )} (\eta )\, \, \, \, } \right|      (2.2)\] 
assuming that the latter have only simple roots.  The `supplementary Jacobi-seed'  polynomials (S$\mathscr{J}$SPs) in the second row of PD (2.2) are defined  as follows
\[\, S_{n+1}^{(\lambda _{+} ,\lambda _{-} )} (\eta )\equiv {\raise0.7ex\hbox{$ 1 $}\!\mathord{\left/ {\vphantom {1 2}} \right. \kern-\nulldelimiterspace}\!\lower0.7ex\hbox{$ 2 $}} [(\lambda _{-} +\lambda _{+} +2)\eta +\lambda _{+} -\lambda _{-} ]\, P_{n}^{(\lambda _{+} ,\lambda _{-} )} (\eta )+(\eta ^{2} -1){\mathop{P}\limits^{\bullet }} \, _{n}^{(\lambda _{+} ,\lambda _{-} )} (\eta )   (2.3)\] 
\[={\raise0.7ex\hbox{$ 1 $}\!\mathord{\left/ {\vphantom {1 2}} \right. \kern-\nulldelimiterspace}\!\lower0.7ex\hbox{$ 2 $}} [(\lambda _{-} +\lambda _{+} +2)\eta +\lambda _{+} -\lambda _{-} ]\, P_{n}^{(\lambda _{+} ,\lambda _{-} )} (\eta )      (2.3*)\] 
\[+{\raise0.7ex\hbox{$ 1 $}\!\mathord{\left/ {\vphantom {1 2}} \right. \kern-\nulldelimiterspace}\!\lower0.7ex\hbox{$ 2 $}} \, (\lambda _{-} +\lambda _{+} +2n+1)(\eta ^{2} -1)\, P_{n-1}^{(\lambda _{+} +1,\lambda _{-} +1)} (\eta ),\] 
with dot standing for the first derivative of the given polynomial with respect to $\eta $.  They satisfy 

\noindent the symmetry relations
\[S_{n}^{(\lambda _{+} ,\lambda _{-} )} (\eta )=(-1)^{n} \, S_{n}^{(\lambda _{-} ,\lambda _{+} )} (-\eta )       (2.3')\] 
which implies that
\[{\rm D}_{n+m+1} [\eta ;\lambda _{+} ,\lambda _{-} |\sigma _{+} \sigma _{-} {\bf ,}{\rm m}]={\bf D}_{n+m+1} [-\eta ;\lambda _{-} ,\lambda _{+} |\sigma _{-} \sigma _{+} {\bf ,}{\rm m}].   (2.2')\] 

By definition
\[{\rm P}_{{\rm N}_{{\rm \{ }\stackrel{\rightharpoonup}{{\bf \sigma }}{\rm m}\} _{2} } }^{{\rm (}\stackrel{\rightharpoonup}{{\rm \lambda }}_{o} )} [\pm |{\rm \{ }\stackrel{\rightharpoonup}{{\rm \sigma }}{\rm m}\} _{2} ]\ne 0        (2.4)\] 
so we set
\[{\rm D}_{m_{1} +m_{2} +1} [\eta ;\stackrel{\rightharpoonup}{{\bf \sigma }}_{2} \times \stackrel{\rightharpoonup}{{\bf \lambda }}_{o} |\stackrel{\rightharpoonup}{{\bf \sigma }}_{1} \times \stackrel{\rightharpoonup}{{\bf \sigma }}_{2} ,{\rm m}_{1} ]\propto \, \prod _{\aleph \, =\, \pm } (1+\aleph \eta )^{\kappa _{\aleph {\rm ;}\stackrel{\rightharpoonup}{{\rm \sigma }}_{1} ,\stackrel{\rightharpoonup}{{\rm \sigma }}_{2} } } {\rm P}_{{\rm N}_{{\rm \{ }\stackrel{\rightharpoonup}{{\rm \sigma }}{\rm m}\} _{2} } }^{{\rm (}\stackrel{\rightharpoonup}{{\rm \lambda }}_{o} )} [\eta |\{ \stackrel{\rightharpoonup}{{\rm \sigma }}{\rm m\} }_{2} ],             (2.5) \] 
where $\kappa _{\aleph {\rm ;}\stackrel{\rightharpoonup}{{\bf \sigma }}_{1} ,\stackrel{\rightharpoonup}{{\bf \sigma }}_{2} } =0\, \, or\, \, 1,$
\[{\rm N}_{{\rm \{ }\stackrel{\rightharpoonup}{{\bf \sigma }}{\rm m}\} _{2} } =m_{1} +m_{2} +1-\kappa _{\stackrel{\rightharpoonup}{{\rm \sigma }}_{1} ,\stackrel{\rightharpoonup}{{\rm \sigma }}_{2} }        (2.6)\] 
with
\[\kappa _{\stackrel{\rightharpoonup}{{\bf \sigma }}_{1} ,\stackrel{\rightharpoonup}{{\bf \sigma }}_{2} } \equiv \sum _{\aleph \, =\, \pm } \kappa _{\aleph {\rm ;}\stackrel{\rightharpoonup}{{\rm \sigma }}_{1} ,\stackrel{\rightharpoonup}{{\rm \sigma }}_{2} } =0,\, 1,\, \, or\, \, 2,      (2.6*)\] 
and 
\[\overrightarrow{{\it a}}\times \overrightarrow{{\it b}}\equiv {\it a}_{-} {\it b}_{-} ,{\it a}_{+} {\it b}_{+}          (2.7)\] 
i.e., ${\bf (}\stackrel{\rightharpoonup}{{\bf \sigma }}_{1} \times \stackrel{\rightharpoonup}{{\bf \sigma }}_{2} )\times \stackrel{\rightharpoonup}{{\bf \sigma }}_{2} =\stackrel{\rightharpoonup}{{\bf \sigma }}_{1} $.  We thus need to consider only four sequences of DPs (2.2), with each of the parameters $\lambda _{-} ,\, \, \lambda _{+} $ changing between ?$\mathrm{\infty}$ and +$\mathrm{\infty}$ or in other words, with the 2D array $\stackrel{\rightharpoonup}{{\bf \lambda }}$ lying in any of four quadrants:

I ($\lambda _{-} >0,\, \, \lambda _{+} >0$), II ($\lambda _{-} <0,\, \, \lambda _{+} >0$), III ($\lambda _{-} <0,\, \, \lambda _{+} <0$), or IV ($\lambda _{-} >0,\, \, \lambda _{+} <0$).

\noindent Let us prove that
\[{\rm P}_{{\rm m}_{1} +{\rm m}_{2} }^{{\rm (|}\lambda _{-} |,|\lambda _{+} |)} [\eta |\stackrel{\rightharpoonup}{{\bf \sigma }}_{1} {\bf ,m}_{1} ;-\stackrel{\rightharpoonup}{{\bf \sigma }}_{1} {\bf ,m}_{2} ]\propto {\bf D}_{m_{1} +m_{2} +1} [\eta ;\stackrel{\rightharpoonup}{{\bf \lambda }}|{\bf -}\, \, {\bf -},{\bf m}_{1} ]   (2.8)\] 
with
\[\stackrel{\rightharpoonup}{{\bf \sigma }}_{1} {\bf =}\left\{\mathop{}\limits_{\mathop{}\limits_{-+\, \, \, \, if\, \, \, \stackrel{\rightharpoonup}{{\bf \lambda }}\in IV;}^{++\, \, \, \, if\, \, \, \stackrel{\rightharpoonup}{{\bf \lambda }}\in III,} }^{\mathop{}\limits_{+-\, \, \, \, if\, \, \, \stackrel{\rightharpoonup}{{\bf \lambda }}\in II,\, \, }^{--\, \, \, if\, \, \, \stackrel{\rightharpoonup}{{\bf \lambda }}\in I,\, \, \, } } \right.         (2.8')\] 

${\rm P}_{{\rm m}_{1} +{\rm m}_{2} }^{{\rm (|}\lambda _{-} |,|\lambda _{+} |)} [\eta |(\stackrel{\rightharpoonup}{{\bf \sigma }}_{1} {\bf ,m}_{1} ;\stackrel{\rightharpoonup}{{\bf \sigma }}_{1}^{*} {\bf ,m}_{2} ]\propto \frac{{\bf D}_{m_{1} +m_{2} +1} [\eta ;\stackrel{\rightharpoonup}{{\bf \lambda }}|{\bf -}\, {\bf +},{\bf m}_{1} ]}{1-\eta } $   (2.8-1i)

\noindent with 

${\bf \sigma }_{1;-}^{*} =-\sigma _{1;-} ,\, \, {\bf \sigma }_{1;+}^{*} =\sigma _{1;+} $        (2.8-1i*)

\noindent and

$\stackrel{\rightharpoonup}{{\bf \sigma }}_{1} {\bf =}\left\{\mathop{}\limits_{\mathop{}\limits_{--\, \, \, \, if\, \, \, \stackrel{\rightharpoonup}{{\bf \lambda }}\in IV;}^{+-\, \, \, \, if\, \, \, \stackrel{\rightharpoonup}{{\bf \lambda }}\in III,} }^{\mathop{}\limits_{++\, \, \, \, if\, \, \, \stackrel{\rightharpoonup}{{\bf \lambda }}\in II,\, \, }^{-+\, \, \, if\, \, \, \stackrel{\rightharpoonup}{{\bf \lambda }}\in I,\, \, \, } } \right. $        (2.8-1i')

${\rm P}_{{\rm m}_{1} +{\rm m}_{2} }^{{\rm (|}\lambda _{-} |,|\lambda _{+} |)} [\eta |(\stackrel{\rightharpoonup}{{\bf \sigma }}_{1} {\bf ,m}_{1} ;\stackrel{\rightharpoonup}{{\bf \sigma }}_{1}^{\dag } {\bf ,m}_{2} ]\propto \frac{{\bf D}_{m_{1} +m_{2} +1} [\eta ;\stackrel{\rightharpoonup}{{\bf \lambda }}|{\bf +}\, {\bf -},{\bf m}_{1} ]}{\eta +1} $   (2.8-1ii)

\noindent with 

${\bf \sigma }_{1;-}^{\dag } =\sigma _{1;-} ,\, \, {\bf \sigma }_{1;+}^{\dag } =-\sigma _{1;+} $        (2.8-1ii*)

\noindent and

$\stackrel{\rightharpoonup}{{\bf \sigma }}_{1} {\bf =}\left\{\mathop{}\limits_{\mathop{}\limits_{++\, \, \, \, if\, \, \, \stackrel{\rightharpoonup}{{\bf \lambda }}\in IV;}^{-+\, \, \, \, if\, \, \, \stackrel{\rightharpoonup}{{\bf \lambda }}\in III,} }^{\mathop{}\limits_{--\, \, \, \, if\, \, \, \stackrel{\rightharpoonup}{{\bf \lambda }}\in II,\, \, }^{+-\, \, \, if\, \, \, \stackrel{\rightharpoonup}{{\bf \lambda }}\in I,\, \, \, } } \right. $        (2.8-1ii')
\[{\rm P}_{{\rm m}_{1} +{\rm m}_{2} -1}^{{\rm (|}\lambda _{-} |,|\lambda _{+} |)} [\eta |\stackrel{\rightharpoonup}{{\bf \sigma }}_{1} {\bf ,m}_{1} ;\stackrel{\rightharpoonup}{{\bf \sigma }}_{1} {\bf ,m}_{2} ]\propto \frac{{\bf D}_{m_{1} +m_{2} +1} [\eta ;\stackrel{\rightharpoonup}{{\bf \lambda }}|{\bf ++},{\bf m}_{1} ]}{1-\eta ^{2} }     (2.8-2)\] 
with
\[\stackrel{\rightharpoonup}{{\bf \sigma }}_{1} {\bf =}\left\{\mathop{}\limits_{\mathop{}\limits_{+-\, \, \, \, if\, \, \, \stackrel{\rightharpoonup}{{\bf \lambda }}\in IV.}^{--\, \, \, \, if\, \, \, \stackrel{\rightharpoonup}{{\bf \lambda }}\in III,} }^{\mathop{}\limits_{-+\, \, \, \, if\, \, \, \stackrel{\rightharpoonup}{{\bf \lambda }}\in II,\, \, }^{++\, \, \, if\, \, \, \stackrel{\rightharpoonup}{{\bf \lambda }}\in I,\, \, \, } } \right.         (2.8-2')\] 

XB-Jacobi eigenpolynomials (2.8) represent the simplest case when the RJ-Bochner eigenpolynomial is nothing but the PD written in the monic form. To our knowledge, they have 

\noindent not been discussed in the literature yet.

XB-Jacobi eigenpolynomials (2.8-1i) are most known since they turn into orthogonal X${}_{m}$-Jacobi OPSs [28, 43] if the parameters$\, \stackrel{\rightharpoonup}{{\bf \lambda }}$ lie in the first quadrant.  To prove this assertion [5] let us first re-write (83) in [28]${}^{\ x)}$ as
\[\lambda _{+} P_{m}^{(-\lambda _{+} ,\lambda _{-} )} (\eta )+(1-\eta ){\mathop{P}\limits^{\bullet }} \, _{m}^{(-\lambda _{+} ,\lambda _{-} )} (\eta )=(\lambda _{+} -m)P_{m}^{(-\lambda _{+} -1,\lambda _{-} +1)} (\eta )  (2.9*)\] 
and then z change $\eta $ for -$\eta $:  

\noindent \_\_\_\_\_\_\_\_\_\_

\noindent ${}^{x)\ }$Note that z-1 should be changed for 1-z in the left-hand side of similar expression (74) in [43].
\[\lambda _{-} P_{m}^{(\lambda _{+} ,-\lambda _{-} )} (\eta )-(\eta +1){\mathop{P}\limits^{\bullet }} \, _{m}^{(\lambda _{+} ,-\lambda _{-} )} (\eta )=(\lambda _{-} -m)P_{m}^{(\lambda _{+} +1,-\lambda _{-} -1)} (\eta ).  (2.9)\] 
As pointed to us by R. Milson, the above expression can be obtained by combining two `obscure' relations between Jacobi polynomials

$P_{n}^{(\alpha ,\beta )} (\eta )=\frac{\alpha +\beta +n+1}{\alpha +n+1} P_{n}^{(\alpha +1,\beta )} (\eta )-\frac{\beta +n}{\alpha +n+1} P_{n}^{(\alpha +1,\beta -1)} (\eta )$   (2.10a)

\noindent and

$P_{n}^{(\alpha +1,\beta )} (\eta )-P_{n}^{(\alpha ,\beta )} (\eta )={\raise0.7ex\hbox{$ 1 $}\!\mathord{\left/ {\vphantom {1 2}} \right. \kern-\nulldelimiterspace}\!\lower0.7ex\hbox{$ 2 $}} (\eta +1)P_{n-1}^{(\alpha +1,\beta +1)} (\eta )$     (2.10b)

\noindent (see 05.06.17.0026.01 and 05.06.17.0008.01 in [44]), coupled with the conventional formula for the first derivative of polynomial (2.10a),
\[{\mathop{P}\limits^{\bullet }} \, _{m}^{(\alpha ,\beta )} (\eta )=\frac{1}{2} (\alpha +\beta +m+1)P_{m-1}^{(\alpha +1,\beta +1)} (\eta ).      (2.11)\] 
After changing $\lambda _{-} $ and n in from the right-hand side of (2.3) for -$\lambda _{-} $ and m accordingly we then use (2.9) to exclude the first derivative of the Jacobi polynomial in such a way:
\[\begin{array}{l} {\, S_{m+1}^{(\lambda _{+} ,-\lambda _{-} )} (\eta )={\raise0.7ex\hbox{$ 1 $}\!\mathord{\left/ {\vphantom {1 2}} \right. \kern-\nulldelimiterspace}\!\lower0.7ex\hbox{$ 2 $}} [(\lambda _{+} +\lambda _{-} +2)\eta +\lambda _{+} -\lambda _{-} ]\, P_{m}^{(\lambda _{+} ,-\lambda _{-} )} (\eta )} \\ {\, \, \, \, \, \, \, \, \, \, \, \, \, \, \, \, \, \, \, \, \, \, \, \, \, \, \, \, \, \, \, \, \, \, \, \, \, \, \, \, \, \, \, \, -(\eta -1)(\lambda _{-} -m)P_{m}^{(\lambda _{+} +1,-\lambda _{-} -1)} (\eta )} \end{array}    (2.12)\] 
that the coefficient of the Jacobi polynomial in the first summand in the sum in the right-hand side of  this relation coincides with the corresponding  coefficient in (2.3).  We thus conclude that PD (2.2) with $\stackrel{\rightharpoonup}{{\bf \sigma }}=-+$can be decomposed according to (2.8-1i) with
\[{\rm P}_{{\rm m}+n}^{{\rm (|}\lambda _{-} |,|\lambda _{+} |)} [\eta |\stackrel{\rightharpoonup}{{\bf \sigma }}_{1} {\bf ,}m;\stackrel{\rightharpoonup}{{\bf \sigma }}_{1}^{*} {\bf ,}n]\propto (\lambda _{-} -m)P_{m}^{(\lambda _{+} +1,-\lambda _{-} -1)} (\eta )P_{n}^{(\lambda _{+} ,\lambda _{-} )} (\eta )\] 
\[+(\eta +1)P_{m}^{(\lambda _{+} ,-\lambda _{-} )} (\eta ){\mathop{P}\limits^{\bullet }} \, _{n}^{(\lambda _{+} ,\lambda _{-} )} (\eta ) (2.13)\] 
assuming that $\stackrel{\rightharpoonup}{{\bf \sigma }}_{1} $ and $\stackrel{\rightharpoonup}{{\bf \sigma }}_{1}^{*} $ are chosen via (2.8-1i') and (2.8-1i*) accordingly, depending on signs of the parameters $\stackrel{\rightharpoonup}{{\bf \lambda }}$.

If the parameters $\stackrel{\rightharpoonup}{{\bf \lambda }}$ lie in the first quadrant then it directly follows from (72) in [43]: 
\[\begin{array}{l} {(\lambda _{-} +n)\stackrel{\frown}{P}_{m,m+n}^{(\lambda _{-} -1,\lambda _{+} +1)} (x)=(\lambda _{-} -m)P_{m}^{(-\lambda _{-} -1,\lambda _{+} +1)} (x)P_{n}^{(\lambda _{-} ,\lambda _{+} )} (x)} \\ {\, \, \, \, \, \, \, \, \, \, \, \, \, \, \, \, \, \, \, \, \, \, \, \, \, \, \, \, \, \, \, \, \, \, \, \, \, \, \, \, \, \, \, \, \, \, \, \, \, \, \, \, \, \, \, \, \, \, \, \, \, \, \, \, \, \, \, \, \, \, \, \, \, \, +(x-1)P_{m}^{(-\lambda _{-} ,\lambda _{+} )} (x){\mathop{P}\limits^{\bullet }} \, _{n}^{(\lambda _{-} ,\lambda _{+} )} (x)} \end{array}  (2.14)\] 
that PD (2.2) is related to the orthogonal X${}_{m}$-Jacobi polynomial in the left side of (2.11) via the

\noindent simple formula  
\[{\rm D}_{m+n+1} [\eta ;\stackrel{\rightharpoonup}{{\bf \lambda }}|-\, +,m{\rm ;++},{\rm n}]=(\lambda _{-} +n)\, (\eta -1)\, \stackrel{\frown}{P}_{m,m+n}^{(\lambda _{-} -1,\lambda _{+} +1)} (-\eta ).   (2.15)\] 
We thus conclude that the XB-Jacobi eigenpolynomial defined via (2.8-1i)-(2.8.1i${}_{\ }$${}_{'}$) is nothing 

\noindent but the monic form of polynomial (2.14) in the reversed argument:  
\[\, (-1)^{m+n} \stackrel{\frown}{k}_{m,m+n}^{(\lambda _{o;-} -1,\lambda _{o;+} +1)} \, {\rm P}_{m+n}^{{\rm (}\stackrel{\rightharpoonup}{{\bf \lambda }}_{o} )} [\eta |-\, +{\bf ,}{\rm m};{\rm ++,}n]=\, \stackrel{\frown}{P}_{m,m+n}^{(\lambda _{o;-} -1,\lambda _{o;+} +1)} (-\eta ) (2.16)\] 
with
\[0\le m<\lambda _{o;-} -1          (2.16o)\] 
(case B in [43]).  The explicit expression for the leading coefficient of the X${}_{m}$-Jacobi polynomial $\stackrel{\frown}{P}_{m,m+n}^{(\alpha ,\beta )} (\eta )$:
\[\stackrel{\frown}{k}_{m,m+n}^{(\alpha ,\beta )} =\frac{\alpha +n+1-m}{\alpha +n+1} k_{m}^{(-\alpha -2,\beta )} k_{n}^{(\alpha +1,\beta -1)} .     (2.16*)\] 
is easily obtained by using conventional formula 
\[k_{n}^{(\alpha ,\beta )} \equiv \frac{{\rm (}\alpha +\beta +n)_{n} \, \, }{n!\, 2^{n} }         (2.17)\] 
(see, i.g., (22.3.1) in [45])  for the leading coefficient of the Jacobi polynomial $P_{n}^{(\alpha ,\beta )} (\eta )$.

Decomposition (2.8-1ii) is obtained in a similar way.  Making use of (2.9*) we again exclude the first derivative of the Jacobi polynomial from the right-hand side of (2.3) with $\lambda _{+} $ and n changed for -$\lambda _{+} $ and m accordingly: 
\[\begin{array}{l} {\, S_{m+1}^{(-\lambda _{+} ,\lambda _{-} )} (\eta )\equiv {\raise0.7ex\hbox{$ 1 $}\!\mathord{\left/ {\vphantom {1 2}} \right. \kern-\nulldelimiterspace}\!\lower0.7ex\hbox{$ 2 $}} [(\lambda _{-} +\lambda _{+} +2)\eta +\lambda _{+} -\lambda _{-} ]\, P_{m}^{(-\lambda _{+} ,\lambda _{-} )} (\eta )} \\ {\, \, \, \, \, \, \, \, \, \, \, \, \, \, \, \, \, \, \, \, \, \, \, \, \, \, \, \, \, \, \, \, \, \, \, \, \, \, -(\eta +1)(\lambda _{+} -m)P_{m}^{(-\lambda _{+} -1,\lambda _{-} +1)} (\eta )} \end{array}    (2.18)\] 
so the coefficient of the Jacobi polynomial in the first summand in the sum in the right-hand side of  this relation coincides with the corresponding coefficient in (2.3).  This brings us to the relation
\[\begin{array}{l} {{\rm P}_{{\rm m}+n}^{{\rm (|}\lambda _{-} |,|\lambda _{+} |)} [\eta |\stackrel{\rightharpoonup}{{\bf \sigma }}_{1} {\bf ,m};\stackrel{\rightharpoonup}{{\bf \sigma }}_{1}^{\dag } {\bf ,}n]\propto P_{m}^{(-\lambda _{+} ,\lambda _{-} )} (\eta )(\eta -1){\mathop{P}\limits^{\bullet }} \, _{n}^{(\lambda _{+} ,\lambda _{-} )} (\eta )} \\ {\, \, \, \, \, \, \, \, \, \, \, \, \, \, \, \, \, \, \, \, \, \, \, \, \, \, \, \, \, \, \, \, \, \, \, \, \, \, \, \, \, \, \, \, \, \, \, \, \, \, \, \, \, \, \, \, \, +(\lambda _{+} -m)P_{m}^{(-\lambda _{+} -1,\lambda _{-} +1)} (\eta )P_{n}^{(\lambda _{+} ,\lambda _{-} )} (\eta )} \end{array}  (2.19)\] 
with $\stackrel{\rightharpoonup}{{\bf \sigma }}_{1}^{\dag } $ defined via (2.8-1ii*).

XB-Jacobi eigenpolynomials (2.8-2)-(2.8-2') with ${\rm m}_{2} ={\rm m}_{1} +1$ and $\stackrel{\rightharpoonup}{{\bf \lambda }}\in I$ were implicitly used in [46, 47] to construct `juxtaposed' [48-50] RD$\mathscr{T}$s of the Darboux/ Pöschl-Teller (D/PT) potential [17, 1] based on the Krein-Adler theorem [51, 52].   The remarkable benchmark of these eigenpolynomials is that they are expressible in terms of Wroskians formed by classical Jacobi polynomials:
\[{\rm D}_{m_{1} +m_{2} +1} [\eta ;\stackrel{\rightharpoonup}{{\bf \lambda }};\stackrel{\rightharpoonup}{{\bf \sigma }}{\bf ,}{\rm m}_{1} ]=(\eta ^{2} -1)\, W\{ P_{m_{1} }^{(\lambda _{+} ,\lambda _{-} )} (\eta ),P_{m_{2} }^{(\lambda _{+} ,\lambda _{-} )} (\eta )\} .    (2.20)\] 
so
\[{\rm P}_{{\rm m}_{1} +{\rm m}_{2} -1}^{{\rm (|}\lambda _{-} |,|\lambda _{+} |)} [\eta |\stackrel{\rightharpoonup}{{\bf \sigma }}_{1} {\bf ,m}_{1} ;\stackrel{\rightharpoonup}{{\bf \sigma }}_{1} {\bf ,m}_{2} ]\propto W\{ P_{m_{1} }^{(\lambda _{+} ,\lambda _{-} )} (\eta ),P_{m_{2} }^{(\lambda _{+} ,\lambda _{-} )} (\eta )\} .   (2.21)\] 
The mentioned benchmark immediately brings to memory the widely cited study of Karlin and 

\noindent Szegö [53] on Wroskians formed by orthogonal polynomials but so far we were unable to 

\noindent obtain any useful results based on this observation.

\noindent \textbf{}

\noindent \textbf{\eject }

\section*{3.  Double-indexed orthogonal polynomial sequences formed by RD$\mathscr{T}$s of R-Jacobi polynomials}

\noindent As originally pointed to by us in [5] (see also [24] for a more general discussion) polynomials (2.1) satisfy the Heine-type equation [21-23]
\[(1-\eta ^{2} )\Pi _{m} [\eta ;\eta _{\stackrel{\rightharpoonup}{{\bf \sigma }}{\rm ,m}} (\stackrel{\rightharpoonup}{{\rm \lambda }}_{o} )]\, {\mathop{{\rm P}}\limits^{\bullet \bullet }} \, _{m+m'+1-\kappa _{\stackrel{\rightharpoonup}{{\rm \sigma }}\stackrel{\rightharpoonup}{{\rm \sigma }}'} }^{{\rm (}\stackrel{\rightharpoonup}{{\rm \lambda }}_{o} )} [\eta |\stackrel{\rightharpoonup}{{\rm \sigma }}{\rm ,m;}\stackrel{\rightharpoonup}{{\rm \sigma }}'{\rm ,}{\rm m}']\] 
\[+\, 2B_{m+1} [\eta ;\stackrel{\rightharpoonup}{{\bf \sigma }}'\times \stackrel{\rightharpoonup}{{\bf \lambda }}_{o} ;\bar{\eta }_{\stackrel{\rightharpoonup}{{\bf \sigma }}{\rm ,m}} (\stackrel{\rightharpoonup}{{\rm \lambda }}_{o} )]\, {\mathop{{\rm P}}\limits^{\bullet }} \, _{m+m'+1-\kappa _{\stackrel{\rightharpoonup}{{\rm \sigma }}\stackrel{\rightharpoonup}{{\rm \sigma }}'} }^{{\rm (}\stackrel{\rightharpoonup}{{\rm \lambda }}_{o} )} [\eta |\stackrel{\rightharpoonup}{{\rm \sigma }}{\rm ,m;}\stackrel{\rightharpoonup}{{\rm \sigma }}'{\rm ,}{\rm m}']  (3.1)\] 
\[+\, \, C_{m} [\eta ;\varepsilon _{\stackrel{\rightharpoonup}{{\bf \sigma }}'m'} ;\stackrel{\rightharpoonup}{{\bf \lambda }}_{o} |\stackrel{\rightharpoonup}{{\bf \sigma }}m;\stackrel{\rightharpoonup}{{\bf \sigma }}']\, {\rm P}_{m+m'+1-\kappa _{\stackrel{\rightharpoonup}{{\rm \sigma }}\stackrel{\rightharpoonup}{{\rm \sigma }}'} }^{{\rm (}\stackrel{\rightharpoonup}{{\rm \lambda }}_{o} )} [\eta |\stackrel{\rightharpoonup}{{\rm \sigma }}{\rm ,m;}\stackrel{\rightharpoonup}{{\rm \sigma }}'{\rm ,}{\rm m}']=0 \] 
with the polynomial coefficients of degree m+2, m+1, and m.   Since the ExpDiffs for the second-order poles of the JRef CSLE (A.1) are energy-independent [5, 40] both polynomial coefficient of the first derivative
\[B_{m+1} [\eta ;\stackrel{\rightharpoonup}{{\bf \lambda }};\bar{\eta }]\equiv \, \Pi _{m} [\eta ;\bar{\eta }]\times \{ \stackrel{\frown}{B}_{1}^{(\stackrel{\rightharpoonup}{{\bf \lambda }})} [\eta ]+(1-\eta ^{2} )\sum _{k=1}^{m}( \eta -\eta _{{\rm k}} )^{-1} \, \, \}   (3.2)\] 
with
\[\stackrel{\frown}{B}_{1}^{(\stackrel{\rightharpoonup}{{\bf \lambda }})} [\eta ]\equiv {\raise0.7ex\hbox{$ 1 $}\!\mathord{\left/ {\vphantom {1 2}} \right. \kern-\nulldelimiterspace}\!\lower0.7ex\hbox{$ 2 $}} [(\lambda _{-} +1)(\eta +1)+(\lambda _{+} +1)(\eta -1)]     (3.21)\] 
and the zero-energy componen 
\[\, C_{m} [\eta ;\stackrel{\rightharpoonup}{{\bf \lambda }}_{o} |\stackrel{\rightharpoonup}{{\bf \sigma }}{\bf ,}m;\stackrel{\rightharpoonup}{{\bf \sigma }}']={\raise0.7ex\hbox{$ 1 $}\!\mathord{\left/ {\vphantom {1 4}} \right. \kern-\nulldelimiterspace}\!\lower0.7ex\hbox{$ 4 $}} \stackrel{\frown}{O}\, _{m}^{\downarrow \, } [\eta ;\stackrel{\rightharpoonup}{{\bf \lambda }}_{o} |\stackrel{\rightharpoonup}{{\bf \sigma }}{\bf ,}m]\, -2\stackrel{\frown}{B}_{1}^{(\stackrel{\rightharpoonup}{{\bf \lambda }}_{o} )} [\eta ]\, \, {\mathop{\Pi }\limits^{\bullet }} _{m} [\eta ;\eta _{\stackrel{\rightharpoonup}{{\bf \sigma }}{\rm ,m}} (\stackrel{\rightharpoonup}{{\rm \lambda }}_{o} )]\] 
\[-{\raise0.7ex\hbox{$ 1 $}\!\mathord{\left/ {\vphantom {1 8}} \right. \kern-\nulldelimiterspace}\!\lower0.7ex\hbox{$ 8 $}} (\lambda _{o;-} +1)(\lambda _{o;+} +1)\, \Pi _{m} [\eta ;\eta _{\stackrel{\rightharpoonup}{{\bf \sigma }}{\rm ,m}} (\stackrel{\rightharpoonup}{{\rm \lambda }}_{o} )]. (3.3)\] 
of the free term
\[C_{m} [\eta ;\varepsilon ;\stackrel{\rightharpoonup}{{\bf \lambda }}_{o} |\stackrel{\rightharpoonup}{{\bf \sigma }}{\bf ,}m;\stackrel{\rightharpoonup}{{\bf \sigma }}']=C_{m} [\eta ;\stackrel{\rightharpoonup}{{\bf \lambda }}_{o} |\stackrel{\rightharpoonup}{{\bf \sigma }}{\bf ,}m;\stackrel{\rightharpoonup}{{\bf \sigma }}']+\varepsilon \, \Pi _{m} [\eta ;\bar{\eta }_{\stackrel{\rightharpoonup}{{\bf \sigma }}m} (\stackrel{\rightharpoonup}{{\bf \lambda }}_{o} )]   (3.4)\] 
are independent of the degree m' of the polynomial forming $\mathscr{J}$S solutions of the given type $\stackrel{\rightharpoonup}{{\bf \sigma }}'$,
\[\phi [\eta ;\stackrel{\rightharpoonup}{{\bf \lambda }}_{o} |\, \stackrel{\rightharpoonup}{{\bf \sigma }}'{\bf ,}m']\propto \, \, \, \prod _{\aleph =-,+} (1+\aleph \eta )^{{\raise0.7ex\hbox{$ 1 $}\!\mathord{\left/ {\vphantom {1 2}} \right. \kern-\nulldelimiterspace}\!\lower0.7ex\hbox{$ 2 $}} (\sigma '_{\aleph } \lambda _{o;\aleph } +1)} \Pi _{m'} [\eta ;\bar{\eta }_{\stackrel{\rightharpoonup}{{\bf \sigma }}'m'} (\stackrel{\rightharpoonup}{{\bf \lambda }}_{o} )].    (3. 5)\] 
Note that the polynomial $\stackrel{\frown}{O}\, _{m}^{\downarrow \, } [\eta ;\stackrel{\rightharpoonup}{{\bf \lambda }}_{o} |\stackrel{\rightharpoonup}{{\bf \sigma }}{\bf ,}m]$in the right-hand side of (3.3) is also independent of the particular choice of second $\mathscr{J}$S solution (3.5) \textit{by definition} because it directly comes from reference polynomial fraction (RefPFr) of transformed CSLE (A.8).

The most important corollary of the stated results is that polynomials (2.1) satisfy the Bochner-type eigenequation with rational coefficients:
\[(1-\eta ^{2} ){\mathop{{\rm P}}\limits^{\bullet \bullet }} \, _{m+m'+1-\kappa _{\stackrel{\rightharpoonup}{{\bf \sigma }}\stackrel{\rightharpoonup}{{\bf \sigma }}'} }^{{\rm (}\stackrel{\rightharpoonup}{{\bf \lambda }}_{o} )} [\eta |\stackrel{\rightharpoonup}{{\bf \sigma }}{\rm ,m;}\stackrel{\rightharpoonup}{{\rm \sigma }}'{\rm ,}{\rm m}']\] 
\[+\, \frac{2B_{m+1} [\eta ;\stackrel{\rightharpoonup}{{\bf \sigma }}'\times \stackrel{\rightharpoonup}{{\bf \lambda }}_{o} ;\bar{\eta }_{\stackrel{\rightharpoonup}{{\bf \sigma }}{\rm ,m}} (\stackrel{\rightharpoonup}{{\rm \lambda }}_{o} )]}{\Pi _{m} [\eta ;\eta _{\stackrel{\rightharpoonup}{{\rm \sigma }}{\rm ,m}} (\stackrel{\rightharpoonup}{{\rm \lambda }}_{o} )]\, } \, {\mathop{{\rm P}}\limits^{\bullet }} \, _{m+m'+1-\kappa _{\stackrel{\rightharpoonup}{{\rm \sigma }}\stackrel{\rightharpoonup}{{\rm \sigma }}'} }^{{\rm (}\stackrel{\rightharpoonup}{{\rm \lambda }}_{o} )} [\eta |\stackrel{\rightharpoonup}{{\rm \sigma }}{\rm ,m;}\stackrel{\rightharpoonup}{{\rm \sigma }}'{\rm ,}{\rm m}']  (3.6)\] 
\[+\left(\frac{C_{m} [\eta ;\varepsilon _{\stackrel{\rightharpoonup}{{\bf \sigma }}'m'} ;\stackrel{\rightharpoonup}{{\bf \lambda }}_{o} |\stackrel{\rightharpoonup}{{\bf \sigma }}m;\stackrel{\rightharpoonup}{{\bf \sigma }}']\, }{\Pi _{m} [\eta ;\eta _{\stackrel{\rightharpoonup}{{\bf \sigma }}{\rm ,m}} (\stackrel{\rightharpoonup}{{\rm \lambda }}_{o} )]} \, -\varepsilon _{\stackrel{\rightharpoonup}{{\rm \sigma }}'{\rm ,}{\rm m}'} \right)\, {\rm P}_{m+m'+1-\kappa _{\stackrel{\rightharpoonup}{{\rm \sigma }}\stackrel{\rightharpoonup}{{\rm \sigma }}'} }^{{\rm (}\stackrel{\rightharpoonup}{{\rm \lambda }}_{o} )} [\eta |\stackrel{\rightharpoonup}{{\rm \sigma }}{\rm ,m;}\stackrel{\rightharpoonup}{{\rm \sigma }}'{\rm ,}{\rm m}']=0.\] 
In next section we formulate the SLP which allows us to select finite subsets of  double-indexed orthogonal polynomials termed `XR-Jacobi polynomials' for briefness.

\noindent 

\section*{4. SUSY pairs of SLPs solvable under the DBCs at the ends of the positive semi-axis}

\noindent The SLP of our current interest is formulated by representing both JRef CSLE (A.1) and its RD$\mathscr{T}$ 

\noindent (A.8) to their `prime' self-adjoint forms
\[\left\{\, \frac{d\, \, }{d\eta } \lower3pt\hbox{\rlap{$\scriptscriptstyle\rightarrow$}}\not {\rm p}_{\diamondsuit } \, [\eta ]\frac{d\, \, }{d\eta } -\lower3pt\hbox{\rlap{$\scriptscriptstyle\rightarrow$}}\not q\, [\eta ;\mathop{\lambda }\limits^{\rightharpoonup}{} _{o} ]\, +\, \lower3pt\hbox{\rlap{$\scriptscriptstyle\rightarrow$}}\varepsilon \, \lower3pt\hbox{\rlap{$\scriptscriptstyle\rightarrow$}}\not {\rm w}_{\diamondsuit } \, [\eta ]\right\}\lower3pt\hbox{\rlap{$\scriptscriptstyle\rightarrow$}}\not \Psi [\eta ;\lower3pt\hbox{\rlap{$\scriptscriptstyle\rightarrow$}}\varepsilon ;\mathop{\lambda }\limits^{\rightharpoonup}{} _{o} ]\, =\, 0\, \, \, \, \, \, (1\le \eta <\infty )    (4.1)\] 
and
\[\left\{\, \frac{d\, \, }{d\eta } \lower3pt\hbox{\rlap{$\scriptscriptstyle\rightarrow$}}\not {\rm p}_{\diamondsuit } \, [\eta ]\frac{d\, \, }{d\eta } -\lower3pt\hbox{\rlap{$\scriptscriptstyle\rightarrow$}}\not q\, [\eta ;\mathop{\lambda }\limits^{\rightharpoonup}{} _{o} |\stackrel{\rightharpoonup}{{\bf \sigma }}{\bf ,}m]\, +\, \, \lower3pt\hbox{\rlap{$\scriptscriptstyle\rightarrow$}}\varepsilon \, \not {\bf w}_{\diamondsuit } \, [\eta ]\right\}\lower3pt\hbox{\rlap{$\scriptscriptstyle\rightarrow$}}\not \Psi [\eta ;\lower3pt\hbox{\rlap{$\scriptscriptstyle\rightarrow$}}\varepsilon ;\mathop{\lambda }\limits^{\rightharpoonup}{} _{o} |\stackrel{\rightharpoonup}{{\bf \sigma }}{\bf ,}m]\, =\, 0\, \, \, \, (1\le \eta <\infty ),  (4.1*)\] 
where the leading coefficient function (LCF) and weight function (WF) are defined as follows
\[\lower3pt\hbox{\rlap{$\scriptscriptstyle\rightarrow$}}\not {\rm p}_{\diamondsuit } \, [\eta ]\equiv \eta -1          (4.2)\] 
and
\[\lower3pt\hbox{\rlap{$\scriptscriptstyle\rightarrow$}}\not {\rm w}_{\diamondsuit } \, [\eta ]\equiv \frac{1}{4(\eta +1)} ,         (4.2')\] 
respectively, and the zero-energy free term is related to the RefPFr in question via the conventional formula
\[\lower3pt\hbox{\rlap{$\scriptscriptstyle\rightarrow$}}\not q\, [\eta ;\mathop{\lambda }\limits^{\rightharpoonup}{} _{o} ]=-\lower3pt\hbox{\rlap{$\scriptscriptstyle\rightarrow$}}\not {\rm p}_{\diamondsuit } {}_{\diamondsuit } \, [\eta ]\, I^{o} [\eta ;\pm 1;\stackrel{\rightharpoonup}{{\bf \lambda }}_{o} ]\, +\frac{1}{4(\eta -1)}       (4.3)\] 
and 
\[\lower3pt\hbox{\rlap{$\scriptscriptstyle\rightarrow$}}\not q\, [\eta ;\mathop{\lambda }\limits^{\rightharpoonup}{} _{o} |\stackrel{\rightharpoonup}{{\bf \sigma }}{\bf ,}m]=-\lower3pt\hbox{\rlap{$\scriptscriptstyle\rightarrow$}}\not {\rm p}_{\diamondsuit } {}_{\diamondsuit } \, [\eta ]\, I^{o} [\eta ;\pm 1;\stackrel{\rightharpoonup}{{\rm \lambda }}_{o} |\stackrel{\rightharpoonup}{{\rm \sigma }}{\rm ,}m]\, +\frac{1}{4(\eta -1)}     (4.3*)\] 
 correspondingly. (Note that we changed $\epsilon$ for $\lower3pt\hbox{\rlap{$\scriptscriptstyle\rightarrow$}}\varepsilon =-\varepsilon $ to make positive the common weight for both RSLEs.)  As discussed by us in [16] in a more general framework  the prime SLE is defined by the requirement that sum of characteristic exponents (ChExps) of the pair of Frobenius solutions at each singular endpoint is equal to zero and as a result the DBC at the given end unambiguously selects the PFS (the Frobenius solution with a larger ChExp)..  

Combining (4.3) and (A.2) one finds
\[\lower3pt\hbox{\rlap{$\scriptscriptstyle\rightarrow$}}\not q\, [\eta ;\mathop{\lambda }\limits^{\rightharpoonup}{} _{o} ]=\, \, \frac{1-\lambda _{o;-}^{2} }{2(\eta +1)^{2} } \, +\frac{\lambda _{o;+}^{2} }{2(\eta ^{2} -1)}        (4.4)\] 
so (as intended) this function does not have a simple pole at infinity.  Making use of (A.13) confirms that the latter assertion also holds for function (4.2*).

Since the RDT \textit{by definition} [16] preserves the LCF and WF of the self-adjoint SLE the assertion that the DBC at the singular endpoint unambiguously selects the corresponding PFS is valid for both SLEs (4.1) and (4.1*).  The SLPs of our interest are thus formulated by imposing the DBCs

\[\lower3pt\hbox{\rlap{$\scriptscriptstyle\rightarrow$}}\not \Psi [1;\lower3pt\hbox{\rlap{$\scriptscriptstyle\rightarrow$}}\varepsilon _{v} ;\mathop{\lambda }\limits^{\rightharpoonup}{} _{o} ]\, =\, 0,\, \, \, {\mathop{\lim }\limits_{\eta \to \infty }} \lower3pt\hbox{\rlap{$\scriptscriptstyle\rightarrow$}}\not \Psi [\eta ;\lower3pt\hbox{\rlap{$\scriptscriptstyle\rightarrow$}}\varepsilon _{v} ;\mathop{\lambda }\limits^{\rightharpoonup}{} _{o} ]\, =\, 0\, \, \, \,       (4.5)\] 
or
\[\lower3pt\hbox{\rlap{$\scriptscriptstyle\rightarrow$}}\not \Psi [1;\lower3pt\hbox{\rlap{$\scriptscriptstyle\rightarrow$}}\varepsilon _{v} ;\mathop{\lambda }\limits^{\rightharpoonup}{} _{o} |\stackrel{\rightharpoonup}{{\bf \sigma }}{\bf ,}m]\, =\, 0,\, \, \, {\mathop{\lim }\limits_{\eta \to \infty }} \lower3pt\hbox{\rlap{$\scriptscriptstyle\rightarrow$}}\not \Psi [\eta ;\lower3pt\hbox{\rlap{$\scriptscriptstyle\rightarrow$}}\varepsilon _{v} ;\mathop{\lambda }\limits^{\rightharpoonup}{} _{o} |\stackrel{\rightharpoonup}{{\bf \sigma }}{\bf ,}m]\, =\, 0     (4.5*)\] 
on solutions of the prime RSLE (4.1) or (4.1*) accordingly.

Let $\lower3pt\hbox{\rlap{$\scriptscriptstyle\rightarrow$}}\not \Psi _{{\it r}} [\eta ;\lower3pt\hbox{\rlap{$\scriptscriptstyle\rightarrow$}}\varepsilon ]$ and $\lower3pt\hbox{\rlap{$\scriptscriptstyle\rightarrow$}}\not \Psi _{{\it r}} [\eta ;\lower3pt\hbox{\rlap{$\scriptscriptstyle\rightarrow$}}\varepsilon ']$ be two `principal' [32-35] solutions of a self-adjoint SLE with a 

\noindent LCF ${}_{1} \not {\rm p}_{\diamondsuit } \, [\eta ]$ at the energies\textit{ $\lower3pt\hbox{\rlap{$\scriptscriptstyle\rightarrow$}}\varepsilon '$and $\lower3pt\hbox{\rlap{$\scriptscriptstyle\rightarrow$}}\varepsilon $, }with \textit{r} specifying the common endpoint where the PFSs in question vanish.  As the direct consequence of such a very specific choice of the self-adjoint 

\noindent modifications of RCSLEs (A.1) and (A.8) the `generalized' [54] Wroskian of two solution

\[\lower3pt\hbox{\rlap{$\scriptscriptstyle\rightarrow$}}\not {\rm W}_{\, \diamondsuit } \{ \lower3pt\hbox{\rlap{$\scriptscriptstyle\rightarrow$}}\not \Psi _{{\it r}} [\eta ;\lower3pt\hbox{\rlap{$\scriptscriptstyle\rightarrow$}}\varepsilon '],\lower3pt\hbox{\rlap{$\scriptscriptstyle\rightarrow$}}\not \Psi _{{\it r}} [\eta ;\lower3pt\hbox{\rlap{$\scriptscriptstyle\rightarrow$}}\varepsilon ]\} \equiv \lower3pt\hbox{\rlap{$\scriptscriptstyle\rightarrow$}}\not {\it p}_{\diamondsuit } \, [\eta ]\, W\{ \lower3pt\hbox{\rlap{$\scriptscriptstyle\rightarrow$}}\not \Psi _{{\it r}} [\eta ;\lower3pt\hbox{\rlap{$\scriptscriptstyle\rightarrow$}}\varepsilon '],\lower3pt\hbox{\rlap{$\scriptscriptstyle\rightarrow$}}\not \Psi _{{\it r}} [\eta ;\lower3pt\hbox{\rlap{$\scriptscriptstyle\rightarrow$}}\varepsilon ]\}     (4.6)\] 
(`Lagrange sesquilinear form or `modified Wronskian determinant' in terms of [55, 56, 39] or [57], respectively) also vanishes at the endpoint \textit{r}:
\[{\mathop{{\it lim}}\limits_{\eta \, \, \to \, {\it r}}} \lower3pt\hbox{\rlap{$\scriptscriptstyle\rightarrow$}}\not {\rm W}_{\, \diamondsuit } \{ \lower3pt\hbox{\rlap{$\scriptscriptstyle\rightarrow$}}\not \Psi _{{\rm r}} [\eta ;\lower3pt\hbox{\rlap{$\scriptscriptstyle\rightarrow$}}\varepsilon '],\lower3pt\hbox{\rlap{$\scriptscriptstyle\rightarrow$}}\not \Psi _{{\rm r}} [\eta ;\lower3pt\hbox{\rlap{$\scriptscriptstyle\rightarrow$}}\varepsilon ]\} =0.\,        (4.7)\] 
The most important consequence of this proposition is that the eigenfunctions ${}_{1} \not \psi _{{\rm v}} [\eta ;\mathop{\lambda }\limits^{\rightharpoonup}{} _{o} ]$and $\lower3pt\hbox{\rlap{$\scriptscriptstyle\rightarrow$}}\not \psi _{{\rm v}} [\eta ;\mathop{\lambda }\limits^{\rightharpoonup}{} _{o} |\stackrel{\rightharpoonup}{{\bf \sigma }}{\bf ,}m]$ of prime SLEs (4.1) and (4.1*) solved under DBCs form orthogonal sets:
\[\int _{1}^{\infty }\lower3pt\hbox{\rlap{$\scriptscriptstyle\rightarrow$}}\not \psi _{{\rm v}'}  [\eta ;\mathop{\lambda }\limits^{\rightharpoonup}{} _{o} ]\, \lower3pt\hbox{\rlap{$\scriptscriptstyle\rightarrow$}}\not \psi _{v} [\eta ;\mathop{\lambda }\limits^{\rightharpoonup}{} _{o} ]\, \lower3pt\hbox{\rlap{$\scriptscriptstyle\rightarrow$}}\not {\rm w}_{\diamondsuit } \, [\eta ]\, d\eta =0\, \, \, \, (v'\ne v)     (4.8)\] 
and
\[\int _{1}^{\infty }\lower3pt\hbox{\rlap{$\scriptscriptstyle\rightarrow$}}\not \psi _{{\rm v}'}  [\eta ;\mathop{\lambda }\limits^{\rightharpoonup}{} _{o} |\stackrel{\rightharpoonup}{{\bf \sigma }}{\bf ,}m]\, \lower3pt\hbox{\rlap{$\scriptscriptstyle\rightarrow$}}\not \psi _{v} [\eta ;\mathop{\lambda }\limits^{\rightharpoonup}{} _{o} |\stackrel{\rightharpoonup}{{\bf \sigma }}{\bf ,}m]\, \lower3pt\hbox{\rlap{$\scriptscriptstyle\rightarrow$}}\not {\rm w}_{\diamondsuit } \, [\eta ]\, d\eta =0\, \, \, \, (v'\ne v)    (4.8*)\] 
for any positive values of the parameters $\lambda _{o;-} $ and  $\lambda _{o;+} $.   After orthogonality conditions (4.8) and (4.8*) has been proven for the eigenfunctions of the prime SLE it can be automatically extended to the eigenfunctions of any other self-adjoint SLE obtained from the given RCSLE via an arbitrarily chosen gauge transformation, including the RCSLE itself.  However, to prove the latter assertion we had to convert the later equation to prime self-adjoint form (4.1*) first.

In contrast with quantization on the finite interval [-1, +1], the SLP in question has a finite number of solutions at the energies
\[\lower3pt\hbox{\rlap{$\scriptscriptstyle\rightarrow$}}\varepsilon _{v} {\rm (}\stackrel{\rightharpoonup}{{\rm \lambda }}_{{\rm o}} {\rm )}=1-(\lambda _{o;-} -\lambda _{o;+} -2v-1)^{2} \, \, \, \, for\, \, \, \, 0\le v\, \le v_{\max } ,    (4.9)\] 
where
\[v_{\max } \equiv \left\lfloor {\raise0.7ex\hbox{$ 1 $}\!\mathord{\left/ {\vphantom {1 2}} \right. \kern-\nulldelimiterspace}\!\lower0.7ex\hbox{$ 2 $}} (\lambda _{o;-} -\lambda _{o;+} -1)\right\rfloor .        (4.9')\] 
The corresponding eigenfunctions expressed in terms of the new variable
\[\lower3pt\hbox{\rlap{$\scriptscriptstyle\rightarrow$}}z={\raise0.7ex\hbox{$ 1 $}\!\mathord{\left/ {\vphantom {1 2}} \right. \kern-\nulldelimiterspace}\!\lower0.7ex\hbox{$ 2 $}} (\eta -1)          (4.10)\] 
are related to R-Jacobi polynomials [8, 58, 59]
\[J_{v}^{(\lambda _{o;+} ,\, -\lambda _{o;-} )} (\lower3pt\hbox{\rlap{$\scriptscriptstyle\rightarrow$}}z)\, \, \equiv \, P_{v}^{(-\lambda _{o;-} ,\, \lambda _{o;+} )} (2\lower3pt\hbox{\rlap{$\scriptscriptstyle\rightarrow$}}z+1)      (4.11)\] 
via the gauge transformation
\[\lower3pt\hbox{\rlap{$\scriptscriptstyle\rightarrow$}}\psi _{v} [2\lower3pt\hbox{\rlap{$\scriptscriptstyle\rightarrow$}}z+1;\mathop{\lambda }\limits^{\rightharpoonup}{} _{o} ]\, \propto \, \, (\lower3pt\hbox{\rlap{$\scriptscriptstyle\rightarrow$}}z+1)^{-{\raise0.7ex\hbox{$ 1 $}\!\mathord{\left/ {\vphantom {1 2}} \right. \kern-\nulldelimiterspace}\!\lower0.7ex\hbox{$ 2 $}} (\lambda _{o;-} -1)} \lower3pt\hbox{\rlap{$\scriptscriptstyle\rightarrow$}}z^{{\raise0.7ex\hbox{$ 1 $}\!\mathord{\left/ {\vphantom {1 2}} \right. \kern-\nulldelimiterspace}\!\lower0.7ex\hbox{$ 2 $}} \lambda _{o;+} } \, J_{v}^{(\lambda _{o;+} ,-\lambda _{o;-} )} (\lower3pt\hbox{\rlap{$\scriptscriptstyle\rightarrow$}}z)   (4.12)\] 
Note that we use Askey's definition of R-Jacobi polynomials according to (1.16) in [8], except that we changed symbol R for J to be able to distinguish orthogonal polynomials (4.11) from their counter-parts: R-Routh polynomials denoted in recent publications [60-63] in exactly the same way.  The described change in the notation seems more consistent with our terminology than touching the notation $R_{v}^{(\alpha ,\beta )} (\eta )$ for R-Routh polynomials.  

Since eigenfunctions (4.12) are orthogonal when integrated with WF (4.2'):
\[\int _{1}^{\infty }\lower3pt\hbox{\rlap{$\scriptscriptstyle\rightarrow$}}\psi _{v'}  [\eta ;\mathop{\lambda }\limits^{\rightharpoonup}{} _{o} ]\, \lower3pt\hbox{\rlap{$\scriptscriptstyle\rightarrow$}}\psi _{v} [\eta ;\mathop{\lambda }\limits^{\rightharpoonup}{} _{o} ]\, \lower3pt\hbox{\rlap{$\scriptscriptstyle\rightarrow$}}\not {\rm w}_{\diamondsuit } \, [\eta ]\, d\eta =0\, \, \, for\, \, 0<v<v'\le v_{\max } ,    (4.13)\] 
polynomials (4.11) form an finite orthogonal basis set with respect to the WF

\noindent \_\_\_\_\_\_\_\_\_\_\_\_\_\_\_\_\_\_\_\_\_\_\_\_\_\_\_

\noindent ${}^{x)}$The argument x of the hypergeometric function in the definition of the R-Jacobi polynomial in [59] must be changed for --x.
\[\lower3pt\hbox{\rlap{$\scriptscriptstyle\rightarrow$}}W(\lower3pt\hbox{\rlap{$\scriptscriptstyle\rightarrow$}}z;\mathop{\lambda }\limits^{\rightharpoonup}{} _{o} )\equiv (\lower3pt\hbox{\rlap{$\scriptscriptstyle\rightarrow$}}z+1)^{-\lambda _{o;-} } \lower3pt\hbox{\rlap{$\scriptscriptstyle\rightarrow$}}z^{\lambda _{o;+} } ,       (4.14)\] 
i.e.,
\[\int _{0}^{\infty } J_{v}^{(\lambda _{o;+} ,-\lambda _{o;-} )} (\lower3pt\hbox{\rlap{$\scriptscriptstyle\rightarrow$}}z)\, J_{v'}^{(\lambda _{o;+} ,-\lambda _{o;-} )} (\lower3pt\hbox{\rlap{$\scriptscriptstyle\rightarrow$}}z)\, \lower3pt\hbox{\rlap{$\scriptscriptstyle\rightarrow$}}W(\lower3pt\hbox{\rlap{$\scriptscriptstyle\rightarrow$}}z;\mathop{\lambda }\limits^{\rightharpoonup}{} _{o} )d\lower3pt\hbox{\rlap{$\scriptscriptstyle\rightarrow$}}z=0\, \, \, for\, \, 0<v<v'\le v_{\max } . (4.15)\] 
It is worth mentioning that upper limit (4.9') for the order of R-Jacobi polynomials is unambiguously determined by the requirement for the integral to converge if v'=v:
\[\int _{0}^{\infty }\,  [J_{v}^{(\lambda _{o;+} ,-\lambda _{o;-} )} (\lower3pt\hbox{\rlap{$\scriptscriptstyle\rightarrow$}}z)]^{2} \lower3pt\hbox{\rlap{$\scriptscriptstyle\rightarrow$}}W(\lower3pt\hbox{\rlap{$\scriptscriptstyle\rightarrow$}}z;\mathop{\lambda }\limits^{\rightharpoonup}{} _{o} )\, d\lower3pt\hbox{\rlap{$\scriptscriptstyle\rightarrow$}}z\, <\infty       (4.16)\] 
as it was originally necessitated by Romanovsky [7, 8].

Our next step is to study the prerequisites for RD$\mathscr{T}$s of R-Jacobi polynomials to form double-indexed finite orthogonal polynomial sequences.  Let ${}_{1} \not \Psi _{{\it r}} [\eta ;\lower3pt\hbox{\rlap{$\scriptscriptstyle\rightarrow$}}\varepsilon ;\mathop{\lambda }\limits^{\rightharpoonup}{} _{o} ]$ be the principal solution of prime RSLE (4.1) near the endpoint \textit{r}, i.e., \textit{by definition}
\[\lower3pt\hbox{\rlap{$\scriptscriptstyle\rightarrow$}}\not \Psi _{+1} [1;\lower3pt\hbox{\rlap{$\scriptscriptstyle\rightarrow$}}\varepsilon ;\mathop{\lambda }\limits^{\rightharpoonup}{} _{o} ]\, =\, 0,\, \, \, \, {\mathop{\lim }\limits_{\eta \to +\infty }} \lower3pt\hbox{\rlap{$\scriptscriptstyle\rightarrow$}}\not \Psi _{+\infty } [\eta ;\lower3pt\hbox{\rlap{$\scriptscriptstyle\rightarrow$}}\varepsilon ;\mathop{\lambda }\limits^{\rightharpoonup}{} _{o} ]\, =\, 0.      (4.17)\] 
Since the general solution $\lower3pt\hbox{\rlap{$\scriptscriptstyle\rightarrow$}}\not \Psi [\eta ;\lower3pt\hbox{\rlap{$\scriptscriptstyle\rightarrow$}}\varepsilon ;\mathop{\lambda }\limits^{\rightharpoonup}{} _{o} ]$of prime RSLE (4.1) can be represented as a superposition of two Frobenius solutions near the given singular endpoint  it necessarily obeys the condition

\noindent 
\[{\mathop{{\it lim}}\limits_{\eta \to \, 1+}} \, \, |(\eta \, -1)\, {\it ld}\, \lower3pt\hbox{\rlap{$\scriptscriptstyle\rightarrow$}}\not \Psi [\eta ;\lower3pt\hbox{\rlap{$\scriptscriptstyle\rightarrow$}}\varepsilon ;\mathop{\lambda }\limits^{\rightharpoonup}{} _{o} ]|\, <+\infty ,  {\mathop{{\it lim}}\limits_{\eta \to \, \infty }} \, \, |\eta \, \, {\it ld}\, \lower3pt\hbox{\rlap{$\scriptscriptstyle\rightarrow$}}\not \Psi [\eta ;\lower3pt\hbox{\rlap{$\scriptscriptstyle\rightarrow$}}\varepsilon ;\mathop{\lambda }\limits^{\rightharpoonup}{} _{o} ]|\, <+\infty .   (4.18)\] 
Keeping in mind that density function (A.2') has the second-order pole at infinity:
\[{\mathop{{\it lim}}\limits_{\eta \to \, \infty }} \, \, \{ \eta ^{2} \, \lower3pt\hbox{\rlap{$\scriptscriptstyle\rightarrow$}}{\it \rho }_{{}_{{\rm \diamondsuit }} } [\eta ]\, \} ={\raise0.7ex\hbox{$ 1 $}\!\mathord{\left/ {\vphantom {1 4}} \right. \kern-\nulldelimiterspace}\!\lower0.7ex\hbox{$ 4 $}} ,        (4.19)\] 
we conclude that the RD$\mathscr{T}$ 
\[\lower3pt\hbox{\rlap{$\scriptscriptstyle\rightarrow$}}\not \Psi _{+\infty } [\eta ;\lower3pt\hbox{\rlap{$\scriptscriptstyle\rightarrow$}}\varepsilon ;\mathop{\lambda }\limits^{\rightharpoonup}{} _{o} |\stackrel{\rightharpoonup}{{\bf \sigma }}{\bf ,}m]=\lower3pt\hbox{\rlap{$\scriptscriptstyle\rightarrow$}}{\it \rho }_{\diamondsuit }^{-{\raise0.7ex\hbox{$ 1 $}\!\mathord{\left/ {\vphantom {1 2}} \right. \kern-\nulldelimiterspace}\!\lower0.7ex\hbox{$ 2 $}} } [\eta ]\, \, \lower3pt\hbox{\rlap{$\scriptscriptstyle\rightarrow$}}\not \Psi _{+\infty } [\eta ;\lower3pt\hbox{\rlap{$\scriptscriptstyle\rightarrow$}}\varepsilon ;\mathop{\lambda }\limits^{\rightharpoonup}{} _{o} ]     (4.20)\] 
\[\times \, \, ({\it ld}\lower3pt\hbox{\rlap{$\scriptscriptstyle\rightarrow$}}\not \Psi _{+\infty } [\eta ;\lower3pt\hbox{\rlap{$\scriptscriptstyle\rightarrow$}}\varepsilon ;\mathop{\lambda }\limits^{\rightharpoonup}{} _{o} ]-{\it ld}\, \lower3pt\hbox{\rlap{$\scriptscriptstyle\rightarrow$}}\not \psi [\eta ;\mathop{\lambda }\limits^{\rightharpoonup}{} _{o} |\stackrel{\rightharpoonup}{{\bf \sigma }}{\bf ,}m])\, \] 
of any principal solution $\lower3pt\hbox{\rlap{$\scriptscriptstyle\rightarrow$}}\not \Psi _{+\infty } [\eta ;\lower3pt\hbox{\rlap{$\scriptscriptstyle\rightarrow$}}\varepsilon ;\mathop{\lambda }\limits^{\rightharpoonup}{} _{o} |\stackrel{\rightharpoonup}{{\bf \sigma }}{\bf ,}m]$ of   PRSLE  (4.1) satisfies the DBC

\noindent 
\[{\mathop{{\it lim}}\limits_{\eta \to \, +\infty }} \, \lower3pt\hbox{\rlap{$\scriptscriptstyle\rightarrow$}}\not \Psi _{+\infty } [\eta ;\lower3pt\hbox{\rlap{$\scriptscriptstyle\rightarrow$}}\varepsilon ;\mathop{\lambda }\limits^{\rightharpoonup}{} _{o} |\stackrel{\rightharpoonup}{{\bf \sigma }}{\bf ,}m]=0,       (4.21)\] 
regardless of the choice of the $\mathscr{J}$S solution
\[\lower3pt\hbox{\rlap{$\scriptscriptstyle\rightarrow$}}\not \psi [\eta ;\mathop{\lambda }\limits^{\rightharpoonup}{} _{o} |\stackrel{\rightharpoonup}{{\bf \sigma }}{\bf ,}m]\propto \, \, (1+\eta )^{{\raise0.7ex\hbox{$ 1 $}\!\mathord{\left/ {\vphantom {1 2}} \right. \kern-\nulldelimiterspace}\!\lower0.7ex\hbox{$ 2 $}} (\sigma _{-} \lambda _{o;-} +1)} (1-\eta )^{{\raise0.7ex\hbox{$ 1 $}\!\mathord{\left/ {\vphantom {1 2}} \right. \kern-\nulldelimiterspace}\!\lower0.7ex\hbox{$ 2 $}} \sigma _{+} \lambda _{o;+} } P_{m}^{(\sigma _{+} \lambda _{o;+} ,\sigma _{-} \lambda _{o;-} )} (\eta ).  (4.22)\] 

Keeping in mind that the lower ChExp for the pole of prime RSLE (4.1) is negative:
\[{}_{1} \not \rho _{{\it r}}^{-} (\lower3pt\hbox{\rlap{$\scriptscriptstyle\rightarrow$}}\varepsilon ;\mathop{\lambda }\limits^{\rightharpoonup}{} _{o} )=-{\raise0.7ex\hbox{$ 1 $}\!\mathord{\left/ {\vphantom {1 2}} \right. \kern-\nulldelimiterspace}\!\lower0.7ex\hbox{$ 2 $}} \lambda _{{\it r}} (-\lower3pt\hbox{\rlap{$\scriptscriptstyle\rightarrow$}}\varepsilon ;\mathop{\lambda }\limits^{\rightharpoonup}{} _{o} )<0,       (4.23)\] 
where $\lambda _{{\it r}} (\varepsilon ;\mathop{\lambda }\limits^{\rightharpoonup}{} _{o} )$ is the ExpDiff for the singular point \textit{r} of the JRef CSLE (\textit{r} = +1 or +$\mathrm{\infty}$), one finds that
\[\, \int _{\eta }^{\infty }\lower3pt\hbox{\rlap{$\scriptscriptstyle\rightarrow$}}\not \Psi ^{2}  [\xi ;\lower3pt\hbox{\rlap{$\scriptscriptstyle\rightarrow$}}\varepsilon ;\mathop{\lambda }\limits^{\rightharpoonup}{} _{o} ]\, \lower3pt\hbox{\rlap{$\scriptscriptstyle\rightarrow$}}\not {\rm w}_{\diamondsuit } \, [\xi ]\, \, d\xi =\int _{\eta }^{\infty } \, \lower3pt\hbox{\rlap{$\scriptscriptstyle\rightarrow$}}\Phi ^{2} [\xi ;\lower3pt\hbox{\rlap{$\scriptscriptstyle\rightarrow$}}\varepsilon ;\mathop{\lambda }\limits^{\rightharpoonup}{} _{o} ]\, \lower3pt\hbox{\rlap{$\scriptscriptstyle\rightarrow$}}{\it \rho }_{{}_{{\rm \diamondsuit }} } [\xi ]d\xi =\infty     (4.24)\] 
where $\lower3pt\hbox{\rlap{$\scriptscriptstyle\rightarrow$}}\Phi [\eta ;\lower3pt\hbox{\rlap{$\scriptscriptstyle\rightarrow$}}\varepsilon ;\mathop{\lambda }\limits^{\rightharpoonup}{} _{o} ]$  and 
\[\lower3pt\hbox{\rlap{$\scriptscriptstyle\rightarrow$}}\not \Psi [\eta ;\lower3pt\hbox{\rlap{$\scriptscriptstyle\rightarrow$}}\varepsilon ;\mathop{\lambda }\limits^{\rightharpoonup}{} _{o} ]=\lower3pt\hbox{\rlap{$\scriptscriptstyle\rightarrow$}}{\it \rho }_{\diamondsuit }^{-{\raise0.7ex\hbox{$ {\rm 1} $}\!\mathord{\left/ {\vphantom {{\rm 1} {\rm 2}}} \right. \kern-\nulldelimiterspace}\!\lower0.7ex\hbox{$ {\rm 2} $}} } [\eta ]\, \, \Phi [\eta ;\lower3pt\hbox{\rlap{$\scriptscriptstyle\rightarrow$}}\varepsilon ;\mathop{\lambda }\limits^{\rightharpoonup}{} _{o} ]       (4.25)\] 
are general solutions of RCSLE (A.1) and prime RSLE (4.1), respectively.  We thus conclude 

\noindent that the pole at infinity is the LP singularity (see [37-39] for the precise definition of the LP/LC cases) regardless of values of $\mathop{\lambda }\limits^{\rightharpoonup}{} _{o} $.  Since RCSLE (A.8) also has a second-order pole at infinity we can apply similar arguments to prime RSLE (4.1*) which gives 
\[\, \int _{\eta }^{\infty }\lower3pt\hbox{\rlap{$\scriptscriptstyle\rightarrow$}}\not \Psi ^{2}  [\xi ;\lower3pt\hbox{\rlap{$\scriptscriptstyle\rightarrow$}}\varepsilon ;\mathop{\lambda }\limits^{\rightharpoonup}{} _{o} |\stackrel{\rightharpoonup}{{\bf \sigma }}{\bf ,}m]\, \lower3pt\hbox{\rlap{$\scriptscriptstyle\rightarrow$}}\not {\rm w}_{\diamondsuit } \, [\xi ]\, \, d\xi =\int _{\eta }^{\infty }\lower3pt\hbox{\rlap{$\scriptscriptstyle\rightarrow$}}\Phi  \, ^{2} [\xi ;\lower3pt\hbox{\rlap{$\scriptscriptstyle\rightarrow$}}\varepsilon ;\mathop{\lambda }\limits^{\rightharpoonup}{} _{o} |\stackrel{\rightharpoonup}{{\rm \sigma }}{\rm ,}m]\, {\it \rho }_{{}_{{\rm \diamondsuit }} } [\xi ]d\xi =\infty ,   (4.24*)\] 
where $\lower3pt\hbox{\rlap{$\scriptscriptstyle\rightarrow$}}\Phi [\eta ;\lower3pt\hbox{\rlap{$\scriptscriptstyle\rightarrow$}}\varepsilon ;\mathop{\lambda }\limits^{\rightharpoonup}{} _{o} |\stackrel{\rightharpoonup}{{\bf \sigma }}{\bf ,}m]$ and 
\[\lower3pt\hbox{\rlap{$\scriptscriptstyle\rightarrow$}}\not \Psi [\eta ;\lower3pt\hbox{\rlap{$\scriptscriptstyle\rightarrow$}}\varepsilon ;\mathop{\lambda }\limits^{\rightharpoonup}{} _{o} |\stackrel{\rightharpoonup}{{\bf \sigma }}{\bf ,}m]=\lower3pt\hbox{\rlap{$\scriptscriptstyle\rightarrow$}}{\it \rho }_{\diamondsuit }^{-{\raise0.7ex\hbox{$ {\rm 1} $}\!\mathord{\left/ {\vphantom {{\rm 1} {\rm 2}}} \right. \kern-\nulldelimiterspace}\!\lower0.7ex\hbox{$ {\rm 2} $}} } [\eta ]\, \, \Phi [\eta ;\lower3pt\hbox{\rlap{$\scriptscriptstyle\rightarrow$}}\varepsilon ;\mathop{\lambda }\limits^{\rightharpoonup}{} _{o} |\stackrel{\rightharpoonup}{{\it \sigma }}{\it ,}m]      (4.25*)\] 
are general solutions of RCSLE (A.8) and prime RSLE (4.1*), respectively.  It is crucial that criterion (4.24) or (4.24*) for the pole of the given SLE at infinity to be a LP singularity is independent of the particular choice of the self-adjoint form.  However we had to use the prime self-adjoint forms (4.1) and (4.1*) of the RCSLEs in question to prove that any RD$\mathscr{T}$ of the PFS $\lower3pt\hbox{\rlap{$\scriptscriptstyle\rightarrow$}}\not \Psi _{+\infty } [\eta ;\lower3pt\hbox{\rlap{$\scriptscriptstyle\rightarrow$}}\varepsilon ;\mathop{\lambda }\limits^{\rightharpoonup}{} _{o} ]$satisfies the DBC at infinity.

The singularities of CSLEs (A.1) and (A.8) at the common finite endpoint of the selected quantization interval represent a more challenging problem because common density function (A.3) has simple poles at $\eta$ = $\pm$1.  As a result the ExpDiffs for the finite singular becomes energy-independent: 
\[\lambda _{\pm 1} (\varepsilon ;\lambda _{o;\pm } |\sigma _{\pm } )=\lambda _{\pm 1} (0|\sigma _{\pm } )\equiv \lambda _{\pm 1;\stackrel{\rightharpoonup}{{\rm \sigma }}_{\pm } }       (4.26)\] 
and therefore change according to (A.10). In particular they change by 1 if the ExpDiff for the singular point $\eta$ = $\pm$1 lies within the LP range [5]:

$\lambda _{o;\pm } >1.$          (4.27P)

\noindent The assertion [63] that inequalities (4.26P) defines the LP range of the ExpDiff for the finite singular points of JRef CSLE (A.1) directly follows from the fact the integral 
\[\, \int _{\pm 1}^{\eta }\lower3pt\hbox{\rlap{$\scriptscriptstyle\rightarrow$}}\not \Psi ^{2}  [\xi ;\vec{\varepsilon };\mathop{\lambda }\limits^{\rightharpoonup}{} _{o} ]\, \lower3pt\hbox{\rlap{$\scriptscriptstyle\rightarrow$}}\not {\rm w}_{\diamondsuit } \, [\xi ]\, \, d\xi =\int _{\pm 1}^{\eta }\lower3pt\hbox{\rlap{$\scriptscriptstyle\rightarrow$}}\Phi  \, ^{2} [\xi ;\lower3pt\hbox{\rlap{$\scriptscriptstyle\rightarrow$}}\varepsilon ;\mathop{\lambda }\limits^{\rightharpoonup}{} _{o} ]\, \lower3pt\hbox{\rlap{$\scriptscriptstyle\rightarrow$}}{\it \rho }_{{}_{{\rm \diamondsuit }} } [\xi ]d\xi      (4.27)\] 
converges if  

\noindent 

$\lambda _{o;\pm } $$\mathrm{<}$ 1.           (4.27C)

Regardless of behavior of the FF near the finite end, the RD$\mathscr{T}$ of the PFS in the neighborhood  of this singular point satisfies the DBC
\[\lower3pt\hbox{\rlap{$\scriptscriptstyle\rightarrow$}}\not \Psi _{+1} [1;\mathop{\lambda }\limits^{\rightharpoonup}{} _{o} |\stackrel{\rightharpoonup}{{\bf \sigma }}{\bf ,}m]=0         (4.28)\] 
if
\[{\mathop{{\it lim}}\limits_{\eta \to 1+}} \, \{ \, (\eta -1)^{-{\raise0.7ex\hbox{$ 1 $}\!\mathord{\left/ {\vphantom {1 2}} \right. \kern-\nulldelimiterspace}\!\lower0.7ex\hbox{$ 2 $}} } \, \lower3pt\hbox{\rlap{$\scriptscriptstyle\rightarrow$}}\not \Psi _{+1} [\eta ;\lower3pt\hbox{\rlap{$\scriptscriptstyle\rightarrow$}}\varepsilon ;\mathop{\lambda }\limits^{\rightharpoonup}{} _{o} |\stackrel{\rightharpoonup}{{\bf \sigma }}{\bf ,}m]\} =0,      (4.28')\] 
i.e., if
\[\lambda _{o;+} >{\raise0.7ex\hbox{$ 1 $}\!\mathord{\left/ {\vphantom {1 2}} \right. \kern-\nulldelimiterspace}\!\lower0.7ex\hbox{$ 2 $}} .          (4.28'')\] 
However this is simply the sufficient condition.  Since ChExps are energy-independent at the given endpoint the right-hand side of the relation
\[\lower3pt\hbox{\rlap{$\scriptscriptstyle\rightarrow$}}\not \Psi _{+1} [\eta ;\lower3pt\hbox{\rlap{$\scriptscriptstyle\rightarrow$}}\varepsilon ;\mathop{\lambda }\limits^{\rightharpoonup}{} _{o} |\stackrel{\rightharpoonup}{{\bf \sigma }}{\bf ,}m]=\lower3pt\hbox{\rlap{$\scriptscriptstyle\rightarrow$}}{\it \rho }_{\, \diamondsuit }^{-{\raise0.7ex\hbox{$ 1 $}\!\mathord{\left/ {\vphantom {1 2}} \right. \kern-\nulldelimiterspace}\!\lower0.7ex\hbox{$ 2 $}} } [\eta ]\, \, \lower3pt\hbox{\rlap{$\scriptscriptstyle\rightarrow$}}\not \Psi _{+1} [\eta ;\lower3pt\hbox{\rlap{$\scriptscriptstyle\rightarrow$}}\varepsilon ;\mathop{\lambda }\limits^{\rightharpoonup}{} _{o} ]\, \, ({\it ld}\lower3pt\hbox{\rlap{$\scriptscriptstyle\rightarrow$}}\not \Psi _{+1} [\eta ;\lower3pt\hbox{\rlap{$\scriptscriptstyle\rightarrow$}}\varepsilon ;\mathop{\lambda }\limits^{\rightharpoonup}{} _{o} ]-{\it ld}\, \lower3pt\hbox{\rlap{$\scriptscriptstyle\rightarrow$}}\not \psi [\eta ;\mathop{\lambda }\limits^{\rightharpoonup}{} _{o} |\stackrel{\rightharpoonup}{{\bf \sigma }}{\bf ,}m])\,              (4.29)\] 
vanishes if the FF $\lower3pt\hbox{\rlap{$\scriptscriptstyle\rightarrow$}}\not \psi [\eta ;\mathop{\lambda }\limits^{\rightharpoonup}{} _{o} |\stackrel{\rightharpoonup}{{\bf \sigma }}{\bf ,}m]$ is the PFS $\lower3pt\hbox{\rlap{$\scriptscriptstyle\rightarrow$}}\not \Psi _{+1} [\eta ;\lower3pt\hbox{\rlap{$\scriptscriptstyle\rightarrow$}}\varepsilon _{\stackrel{\rightharpoonup}{{\bf \sigma }}{\bf ,}m} ;\mathop{\lambda }\limits^{\rightharpoonup}{} _{o} ]$ of PSLE (4.1) in the neighborhood of the finite singular endpoint, i.e., if $\sigma _{+} =+.$ Therefore constraint (4.28'') is applicable only if $\sigma _{+} =-.$  If

${\raise0.7ex\hbox{$ 1 $}\!\mathord{\left/ {\vphantom {1 2}} \right. \kern-\nulldelimiterspace}\!\lower0.7ex\hbox{$ 2 $}} <\lambda _{o;+} <1$ and $\sigma _{+} =-$        (4.30)

\noindent then
\[0<\lambda _{+1;-} <{\raise0.7ex\hbox{$ 1 $}\!\mathord{\left/ {\vphantom {1 2}} \right. \kern-\nulldelimiterspace}\!\lower0.7ex\hbox{$ 2 $}}          (4.30')\] 
which implies that use of any non-principal FS as the  FF for the RDT requires special consideration.

Assuming that $\mathscr{J}$S solution (4.22) is `R-Jacobi admissible' in Duran's terms [42], i.e., that it does not have nodes between +1 and +$\mathrm{\infty}$ in \textit{our} case -- the solutions of prime RSLE (4.1*),
\[\begin{array}{l} {\lower3pt\hbox{\rlap{$\scriptscriptstyle\rightarrow$}}\not \psi _{v} [\eta ;\mathop{\lambda }\limits^{\rightharpoonup}{} _{o} |\stackrel{\rightharpoonup}{{\bf \sigma }}{\bf ,}m]\equiv \lower3pt\hbox{\rlap{$\scriptscriptstyle\rightarrow$}}\not \Psi [\eta ;\lower3pt\hbox{\rlap{$\scriptscriptstyle\rightarrow$}}\varepsilon _{v} ;\mathop{\lambda }\limits^{\rightharpoonup}{} _{o} |\stackrel{\rightharpoonup}{{\bf \sigma }}{\bf ,}m]} \\ {\, \, \, \, \, \, \, \, \, \, \, \, \, \, \, \, \, \, \, \, \propto \, (1+\eta )^{-{\raise0.7ex\hbox{$ 1 $}\!\mathord{\left/ {\vphantom {1 2}} \right. \kern-\nulldelimiterspace}\!\lower0.7ex\hbox{$ 2 $}} \lambda _{o;-} } (1-\eta )^{{\raise0.7ex\hbox{$ 1 $}\!\mathord{\left/ {\vphantom {1 2}} \right. \kern-\nulldelimiterspace}\!\lower0.7ex\hbox{$ 2 $}} (\lambda _{o;+} -1)} \frac{{\rm D}_{m+v+1} [\eta ;-\lambda _{o;-} ,\lambda _{o;+} |\stackrel{\rightharpoonup}{{\rm \sigma }}{\rm ,}m]}{P_{m}^{(\sigma _{+} \lambda _{o;+} ,\sigma _{-} \lambda _{o;-} )} (\eta )} ,} \end{array} (4.31)\] 
which are constructed by applying the RDT to eigenfunctions of prime RSLE (4.1) solved under DBCs (4.5), necessarily satisfy  DBCs (4.5*) and thereby  obey the orthogonality relation
\[\int _{1}^{\infty }\lower3pt\hbox{\rlap{$\scriptscriptstyle\rightarrow$}}\not \psi _{v'}  [\eta ;\mathop{\lambda }\limits^{\rightharpoonup}{} _{o} |\stackrel{\rightharpoonup}{{\bf \sigma }}{\bf ,}m]\, \lower3pt\hbox{\rlap{$\scriptscriptstyle\rightarrow$}}\not \psi _{v} [\eta ;\mathop{\lambda }\limits^{\rightharpoonup}{} _{o} |\stackrel{\rightharpoonup}{{\bf \sigma }}{\bf ,}m]\, \lower3pt\hbox{\rlap{$\scriptscriptstyle\rightarrow$}}\not {\rm w}_{\diamondsuit } \, [\eta ]\, d\eta =0\, \, \, for\, \, 0<v<v'\le v_{\max }   (4. 32)\] 
with the WF $\lower3pt\hbox{\rlap{$\scriptscriptstyle\rightarrow$}}\not {\rm w}_{\diamondsuit } \, [\eta ]\, $.  In next section we present a full list of admissible $\mathscr{J}$S solutions making it possible to solve prime RSLE (4.1*) under DBCs (4.5*) via the double-indexed XR-Jacobi polynomials.  Note that for this purpose we need to consider only PDs (2.2) with $\stackrel{\rightharpoonup}{{\rm \lambda }}$ lying in quadrant II.  

\noindent \textbf{}

\section*{5.  Three sets of nodeless principal $\mathscr{J}$S solutions}

\noindent Our next step is to accumulate all $\mathscr{J}$S solutions (4.22) with positive ChExps near one of the singular endpoints of the infinite interval [+1, +$\mathrm{\infty}$) and then select the subsets of the collected principal solutions below the lowest eigenvalue. To do it we re-label $\mathscr{J}$S solutions (4.22) at the energies 
\[{\rm \varepsilon }_{\lower3pt\hbox{\rlap{$\scriptscriptstyle\rightarrow$}}{\bf t}_{\stackrel{\rightharpoonup}{{\bf \sigma }}} \, {\rm m}} {\bf (}\stackrel{\rightharpoonup}{{\bf \lambda }}_{{\bf o }} {\bf )}=\, -{\bf \varepsilon }_{\stackrel{\rightharpoonup}{{\bf \sigma }}\, {\rm m}} {\bf (}\stackrel{\rightharpoonup}{{\bf \lambda }}_{{\bf o }} {\bf )}\equiv 1-(\sigma _{-} \lambda _{o;-} +\sigma _{+} \lambda _{o;+} +2m+1)^{2}     (5.1)\] 
as specified in Table 5.1.

By analogy with the general case of the second-order pole density function ${}_{1} {\it \rho }[z;{}_{1} {\it T}_{2} ]$ thoroughly examined in [16] there are four primary sequences \textbf{t}m of distinct types \textbf{t }= \textbf{a}, \textbf{b}, \textbf{c}, and \textbf{d }starting from basic solutions \textbf{t}0 so each eigenfunction \textbf{c}v is accompanied by three distinct-type $\mathscr{J}$S solutions formed by Jacobi polynomials of the same order v.

\noindent Table 5.1

\noindent \textbf{Classification of $\mathscr{J}$S solutions based on their asymptotic behavior near the singular endpoints}

\begin{tabular}{|p{0.8in}|p{1.4in}|p{2.0in}|} \hline 
\textbf{$\lower3pt\hbox{\rlap{$\scriptscriptstyle\rightarrow$}}{\bf t}_{\stackrel{\rightharpoonup}{{\bf \sigma }}} $} & \textbf{$\sigma _{-} \, \, \sigma _{+} \, \, \sigma _{\infty } $} & m \\ \hline 
\textbf{a} & +  +   --   & $0\le m<\infty $ \\ \hline 
\textbf{ a}' & --   +   -- & $m>{\raise0.7ex\hbox{$ 1 $}\!\mathord{\left/ {\vphantom {1 2}} \right. \kern-\nulldelimiterspace}\!\lower0.7ex\hbox{$ 2 $}} (\lambda _{o;-} -\lambda _{o;+} -1)$ \\ \hline 
\textbf{b} & -- -- + & $0\le m<{\raise0.7ex\hbox{$ 1 $}\!\mathord{\left/ {\vphantom {1 2}} \right. \kern-\nulldelimiterspace}\!\lower0.7ex\hbox{$ 2 $}} (\lambda _{o;-} +\lambda _{o;+} -1)$ \\ \hline 
\textbf{b}'\textbf{   } & -- +  + & $0\le m<{\raise0.7ex\hbox{$ 1 $}\!\mathord{\left/ {\vphantom {1 2}} \right. \kern-\nulldelimiterspace}\!\lower0.7ex\hbox{$ 2 $}} (\lambda _{o;+} -\lambda _{o;-} -1);$ \\ \hline 
\textbf{c} & -- + + & $0\le m<n_{{\bf c}} \equiv \left\lfloor {\raise0.7ex\hbox{$ 1 $}\!\mathord{\left/ {\vphantom {1 2}} \right. \kern-\nulldelimiterspace}\!\lower0.7ex\hbox{$ 2 $}} (\lambda _{o;-} -\lambda _{o;+} -1)\right\rfloor ;$ \\ \hline 
\textbf{d} & + -- -- & ${\raise0.7ex\hbox{$ 1 $}\!\mathord{\left/ {\vphantom {1 2}} \right. \kern-\nulldelimiterspace}\!\lower0.7ex\hbox{$ 2 $}} (\lambda _{o;+} -\lambda _{o;-} -1)\le m<\infty $ \\ \hline 
\textbf{d}'\textbf{  } & -- -- -- \textbf{   } & $m>{\raise0.7ex\hbox{$ 1 $}\!\mathord{\left/ {\vphantom {1 2}} \right. \kern-\nulldelimiterspace}\!\lower0.7ex\hbox{$ 2 $}} (\lambda _{o;-} +\lambda _{o;+} -1)$ \\ \hline 
\end{tabular}

The infinite sequence \textbf{a}n is formed by classical Jacobi polynomials so all the zeros of these  $\mathscr{J}$S solutions are restricted to the finite interval [?1, +1].  As expected all the $\mathscr{J}$S solutions from this sequence lie below the lowest eigenvalue:
\[\lower3pt\hbox{\rlap{$\scriptscriptstyle\rightarrow$}}{\rm \varepsilon }_{{\bf a}{\rm n}} {\bf (}\stackrel{\rightharpoonup}{{\bf \lambda }}_{{\bf o }} {\bf )-}\lower3pt\hbox{\rlap{$\scriptscriptstyle\rightarrow$}}{\bf \varepsilon }_{{\bf c}{\rm 0}} {\bf (}\stackrel{\rightharpoonup}{{\bf \lambda }}_{{\bf o }} {\bf )}=\, (\lambda _{o;+} -\lambda _{o;-} +1)^{2} -(\lambda _{o;+} +\lambda _{o;-} +2n+1)^{2} <0.  (5.2)\] 
We thus conclude that there is an infinite set of seed function with no nodes for 1 $\leq$ $\eta$ $\mathrm{<}$ $\mathrm{\infty}$ and 

\noindent therefore any of them can be used for constructing \textit{infinite} DCN of finite orthogonal sets 

\noindent composed of XR-Jacobi polynomials.  As stressed at the very end of the previous section we are only interested in PDs (2.2) with $\stackrel{\rightharpoonup}{{\rm \lambda }}\in $ II ($\stackrel{\rightharpoonup}{{\bf \sigma }}_{2} =-+$, \textbf{t${}_{{\mathop{\stackrel{\rightharpoonup}{{\bf \sigma }}_{2} }\limits^{}} }$}= \textbf{c}).  Examination of conditions (2.8'), (2.8-1i${}_{'}$), (2.8-1ii${}_{'}$), and (2.8-2') reveals that we deal with XB-Jacobi polynomials (2.8-1i). 

\noindent  Defining XR-Jacobi polynomials of this type as
\[{\rm J}_{n+v}^{{\rm (}\stackrel{\rightharpoonup}{{\bf \lambda }}_{o} )} [\eta |{\bf a}{\rm n}]\equiv {\rm P}_{n+v}^{{\rm (}\stackrel{\rightharpoonup}{{\rm \lambda }}_{o} )} [\eta |{\rm ++,}n;-\, +{\rm ,}{\rm v}].      (5.3)\] 
with
\[0\le v<{\raise0.7ex\hbox{$ 1 $}\!\mathord{\left/ {\vphantom {1 2}} \right. \kern-\nulldelimiterspace}\!\lower0.7ex\hbox{$ 2 $}} (\lambda _{o;-} -\lambda _{o;+} -1)<\lambda _{o;-} -1      (5.3o)\] 
and substituting (2.8-1), coupled with (5.3), into the right hand side of (4.31) brings us to the finite set of quasi-rational solutions of prime SLE (4.1*),
\[\begin{array}{l} {\lower3pt\hbox{\rlap{$\scriptscriptstyle\rightarrow$}}\not \psi _{v} [\eta ;\mathop{\lambda }\limits^{\rightharpoonup}{} _{o} |{\bf a}n]\equiv \lower3pt\hbox{\rlap{$\scriptscriptstyle\rightarrow$}}\not \psi _{v} [\eta ;\mathop{\lambda }\limits^{\rightharpoonup}{} _{o} |{\bf +}\, {\bf +,}n]} \\ {\, \, \, \, \, \, \, \, \, \, \, \, \, \, \, \, \, \, \, \, \, \, \, \, \, \, \, \, \, \, \, \, \, \, \, \propto \, (1+\eta )^{-{\raise0.7ex\hbox{$ 1 $}\!\mathord{\left/ {\vphantom {1 2}} \right. \kern-\nulldelimiterspace}\!\lower0.7ex\hbox{$ 2 $}} \lambda _{o;-} } (\eta -1)^{{\raise0.7ex\hbox{$ 1 $}\!\mathord{\left/ {\vphantom {1 2}} \right. \kern-\nulldelimiterspace}\!\lower0.7ex\hbox{$ 2 $}} (\lambda _{o;+} +1)} \frac{{\rm J}_{n+v}^{{\rm (}\stackrel{\rightharpoonup}{{\bf \lambda }}_{o} )} [\eta |{\bf a}{\rm n}]}{P_{n}^{(\lambda _{o;+} ,\lambda _{o;-} )} (\eta )} ,} \end{array}  (5.4)\] 
orthogonal with the WF $\lower3pt\hbox{\rlap{$\scriptscriptstyle\rightarrow$}}\not {\rm w}_{\diamondsuit } \, [\eta ]$.  Representing $-{\raise0.7ex\hbox{$ 1 $}\!\mathord{\left/ {\vphantom {1 2}} \right. \kern-\nulldelimiterspace}\!\lower0.7ex\hbox{$ 2 $}} \lambda _{o;-} $ as ${\raise0.7ex\hbox{$ 1 $}\!\mathord{\left/ {\vphantom {1 2}} \right. \kern-\nulldelimiterspace}\!\lower0.7ex\hbox{$ 2 $}} -{\raise0.7ex\hbox{$ 1 $}\!\mathord{\left/ {\vphantom {1 2}} \right. \kern-\nulldelimiterspace}\!\lower0.7ex\hbox{$ 2 $}} (\lambda _{o;-} +1)$ we confirm that the RDT in question increases the ExpDiffs for the poles at $\eta$=$\pm$1 by 1.   It is remarkable that R-Jacobi polynomials (5.3) coincide with polynomials (2.14) from different X${}_{m}$-OPSs for
\[0\le m<{\raise0.7ex\hbox{$ 1 $}\!\mathord{\left/ {\vphantom {1 2}} \right. \kern-\nulldelimiterspace}\!\lower0.7ex\hbox{$ 2 $}} (\lambda _{o;-} -\lambda _{o;+} -1)<\lambda _{o;-} -1,      (5.5o)\] 
After the latter polynomials are written in its monic form
\[\, {\rm J}_{n+m}^{{\rm (}\stackrel{\rightharpoonup}{{\bf \lambda }}_{o} )} [\eta |{\bf a}{\rm n}]=\, \, (-1)^{m+n} \, \stackrel{\frown}{P}_{m,m+n}^{(\lambda _{o;-} -1,\lambda _{o;+} +1)} (-\eta )/\stackrel{\frown}{k}_{m,m+n}^{(\lambda _{o;-} -1,\lambda _{o;+} +1)}   (5.5)\] 
with $\stackrel{\frown}{k}_{m,m+n}^{(\lambda _{o;-} -1,\lambda _{o;+} +1)} $ standing for the leading coefficient of the X${}_{m}$-Jacobi polynomial in $\eta$ [see (2.16*) in Section 2].   Note that this assertion is not applicable to the X${}_{m}$-Jacobi OPSs with 

\noindent m lying within the range
\[0<{\raise0.7ex\hbox{$ 1 $}\!\mathord{\left/ {\vphantom {1 2}} \right. \kern-\nulldelimiterspace}\!\lower0.7ex\hbox{$ 2 $}} (\lambda _{o;-} -\lambda _{o;+} -1)<m<\lambda _{o;-} -1.      (5.6)\] 

Substituting (5.4) into (4.32) we conclude that the double-indexed XR-Jacobi polynomials (5.3) form a finite orthogonal polynomial set:
\[\int _{1}^{\infty }J \, _{n+v}^{{\rm (}\stackrel{\rightharpoonup}{{\bf \lambda }}_{o} )} [\eta |{\bf a}{\rm n}]\, {\rm J}_{n+\tilde{v}}^{{\rm (}\stackrel{\rightharpoonup}{{\bf \lambda }}_{o} )} [\eta |{\bf a}{\rm n}]\, W[\eta ;\stackrel{\rightharpoonup}{{\bf \lambda }}_{o} |{\bf a}n]d\eta =0     (5.7)\] 

       for $0\le v<\tilde{v}<{\raise0.7ex\hbox{$ 1 $}\!\mathord{\left/ {\vphantom {1 2}} \right. \kern-\nulldelimiterspace}\!\lower0.7ex\hbox{$ 2 $}} (\lambda _{o;-} -\lambda _{o;+} -1)<\lambda _{o;-} -1,$  

\noindent with the weight function
\[W[\eta ;\stackrel{\rightharpoonup}{{\bf \lambda }}_{o} |{\bf a}n]\equiv \frac{(\eta +1)^{-\lambda _{o;-} -1} (\eta -1)^{\lambda _{o;+} +1} }{[P_{n}^{(\lambda _{o;-} ,\lambda _{o;+} )} (\eta )]^{2} }  .     (5.8)\] 
As expected selected range (5.3${}^{o}$) for v is equivalent to the condition
\[2v+\lambda _{o;+} -\lambda _{o;-} <-1         (5.8o)\] 
for the square integral of polynomials (5.3) to converge:
\[\int _{1}^{\infty }( {\rm J}_{n+v}^{{\rm (}\stackrel{\rightharpoonup}{{\bf \lambda }}_{o} )} [\eta |{\bf a}{\rm n}])\, ^{2} \, W[\eta ;\stackrel{\rightharpoonup}{{\bf \lambda }}_{o} |{\bf a}n]\, d\eta <\infty .      (5.9)\] 

Substituting (2.9) into (2.13) with $\mathop{\lambda }\limits^{\rightharpoonup}{} $= ?$\lambda _{o;-} ,\lambda _{o;+} $and m and n changed for n and v 

\noindent accordingly, one finds
\[{\rm P}_{v}^{{\rm (}\lambda _{o;-} ,\lambda _{o;+} )} [\eta |++{\bf ,}n;-+{\bf ,}v]\propto -(\lambda _{o;-} +n)P_{n}^{(\lambda _{o;+} -1,\lambda _{o;-} +1)} (\eta )P_{v}^{(\lambda _{o;+} ,-\lambda _{o;-} )} (\eta )\] 
\[+(\eta +1)P_{n}^{(\lambda _{o;+} ,\lambda _{o;-} )} (\eta ){\mathop{P}\limits^{\bullet }} \, _{v}^{(\lambda _{o;+} ,-\lambda _{o;-} )} (\eta ) (5.10)\] 
so the sequence of XR-Jacobi polynomials (5.3) starts from the monic Jacobi polynomial 
\[{\rm P}_{n}^{{\rm (}\lambda _{o;-} ,\lambda _{o;+} )} [\eta |++{\bf ,}n;-+{\bf ,}0]=P_{n}^{(\lambda _{o;+} -1,\lambda _{o;-} +1)} (\eta )/k_{n}^{(\lambda _{o;+} -1,\lambda _{o;-} +1)} .  (5.11)\] 

On other hand, substituting (2.9) into (5.10) and setting m=0 gives
\[{\rm P}_{n}^{{\rm (}\lambda _{o;-} ,\lambda _{o;+} )} [\eta |++{\bf ,}0;-+{\bf ,}n]=P_{n}^{(\lambda _{o;+} +1,-\lambda _{o;-} -1)} (\eta )/k_{n}^{(\lambda _{o;+} +1,-\lambda _{o;-} -1)} . (5.12)\] 
The RDT with the basic FF ${\bf a}0$ thus simply increases by 1 the ExpDiffs for the singular points $\eta =\pm 1$,  without creating new poles -- the direct consequence of the `shape-invariance' [65] of the appropriate Liouville potential.  Making use of (4.14) one finds that
\[W[\eta ;\stackrel{\rightharpoonup}{{\bf \lambda }}_{o} |{\bf a}0]\equiv 2^{\lambda _{o;+} -\lambda _{o;-} } \lower3pt\hbox{\rlap{$\scriptscriptstyle\rightarrow$}}W({\raise0.7ex\hbox{$ 1 $}\!\mathord{\left/ {\vphantom {1 2}} \right. \kern-\nulldelimiterspace}\!\lower0.7ex\hbox{$ 2 $}} \eta -{\raise0.7ex\hbox{$ 1 $}\!\mathord{\left/ {\vphantom {1 2}} \right. \kern-\nulldelimiterspace}\!\lower0.7ex\hbox{$ 2 $}} ;\lambda _{o;-} +1,\lambda _{o;+} +1),   (5.12*)\] 
in agreement with (5.12).

In next section we will come back to discussion of this remarkable polynomial subset formed 

\noindent by orthogonal X-Jacobi polynomials from different OPSs.  As for the rest of this section we 

\noindent solely focus to two other sets of nodeless principal $\mathscr{J}$S solutions, \textbf{b}m\textbf{${}_{b}$} and \textbf{a}'m\textbf{${}_{a}$}${}_{'}$,  which can be 

\noindent added to the aforementioned infinite sequence \textbf{a}m to further expand the RDCN of the SLPs solvable via finite sets of multi-index XR-Jacobi polynomials.  (The $\mathscr{J}$S solutions \textbf{b}'m\textbf{${}_{\ }$}does not\textbf{${}_{\ }$}co-exist with the discrete spectrum and for this reason are excluded from the further analysis.)

The $\mathscr{J}$S solutions \textbf{a}'m\textbf{${}_{a}$}${}_{'}$${}_{\ }$from the secondary sequence \textbf{a}'m are nodeless iff
\[\lower3pt\hbox{\rlap{$\scriptscriptstyle\rightarrow$}}{\rm \varepsilon }_{{\bf a}'{\rm m}_{{\bf a}'} } {\bf (}\stackrel{\rightharpoonup}{{\bf \lambda }}_{{\bf o }} {\bf )-}\lower3pt\hbox{\rlap{$\scriptscriptstyle\rightarrow$}}{\bf \varepsilon }_{{\bf c}{\rm 0}} {\bf (}\stackrel{\rightharpoonup}{{\bf \lambda }}_{{\bf o }} {\bf )}=\, (\lambda _{o;+} -\lambda _{o;-} +1)^{2} -(\lambda _{o;+} -\lambda _{o;-} +2{\rm m}_{{\bf a}'} +1)^{2} \] 
\[=\, -4{\rm m}_{{\bf a}'} ({\rm m}_{{\bf a}'} +1-\lambda _{o;+} +\lambda _{o;-} )<0    (5.13)\] 
so
\[{\rm m}_{{\bf a}'} >\lambda _{o;-} -\lambda _{o;+} -1>0.      (5.13o)\] 
Note that the lower bound for ${\rm m}_{{\bf a}'} $ is either twice larger than the number of eigenfunctions, ${\rm n}_{{\bf c}} $ 

\noindent (if $\left\lfloor \lambda _{o;-} -\lambda _{o;+} -1\right\rfloor $ is even) or is equal to ${\rm 2n}_{{\bf c}} +1$(if $\left\lfloor \lambda _{o;-} -\lambda _{o;+} -1\right\rfloor $ is odd).  It directly follows from (5.13${}^{o}$) that prime RSLE (4.1) has infinitely many R-Jacobi admissible $\mathscr{J}$S solutions \textbf{a}'m\textbf{${}_{a}$}${}_{'}$${}_{.\ \ }$${}_{\ }$Since, according to Table 5.1,  $\stackrel{\rightharpoonup}{{\bf \sigma }}_{1} =-+$ in this case and $\stackrel{\rightharpoonup}{{\rm \lambda }}\in {\rm {\rm I} {\rm I} }$ we come to double 

\noindent indexed XB-Jacobi (2.8-2):
\[{\rm J}_{m_{{\bf a}'} +v-1}^{{\rm (}\stackrel{\rightharpoonup}{{\bf \lambda }}_{o} )} [\eta |{\bf a}'m_{{\bf a}'} ]\equiv {\rm P}_{m_{{\rm a}'} +v-1}^{{\rm (}\stackrel{\rightharpoonup}{{\rm \lambda }}_{o} )} [\eta |-\, +{\rm ,}m_{{\rm a}'} ;{\rm -}\, {\rm +,}v]    (5.14)\] 
with v and m\textbf{${}_{a}$}${}_{'}$${}_{\ }$restricted by constraints (5.3${}^{o}$) and (5.13${}^{o}$), respectively.  Substituting (2.8-2), 

\noindent coupled with (5.14), into the right hand side of (4.31) we come to another set of quasi-rational solutions of prime RSLE (4.1*),

\noindent 
\[\begin{array}{l} {\lower3pt\hbox{\rlap{$\scriptscriptstyle\rightarrow$}}\not \psi _{v} {}_{v} [\eta ;\mathop{\lambda }\limits^{\rightharpoonup}{} _{o} |{\bf a}'m_{{\bf a}'} ]\equiv \lower3pt\hbox{\rlap{$\scriptscriptstyle\rightarrow$}}\not \psi _{v} [\eta ;\mathop{\lambda }\limits^{\rightharpoonup}{} _{o} |{\bf -}\, {\bf +,}m_{{\bf a}'} ]} \\ {\, \, \, \, \, \, \, \, \, \, \, \, \, \, \, \, \, \, \, \, \propto \, (1+\eta )^{-{\raise0.7ex\hbox{$ 1 $}\!\mathord{\left/ {\vphantom {1 2}} \right. \kern-\nulldelimiterspace}\!\lower0.7ex\hbox{$ 2 $}} \lambda _{o;-} +1} (\eta -1)^{{\raise0.7ex\hbox{$ 1 $}\!\mathord{\left/ {\vphantom {1 2}} \right. \kern-\nulldelimiterspace}\!\lower0.7ex\hbox{$ 2 $}} (\lambda _{o;+} +1)} \frac{{\rm J}_{m_{{\bf a}'} +v-1}^{{\rm (}\stackrel{\rightharpoonup}{{\bf \lambda }}_{o} )} [\eta |{\bf a}'m_{{\bf a}'} ]}{P_{n}^{(\lambda _{o;+} ,-\lambda _{o;-} )} (\eta )} ,} \end{array}  (5.15)\] 
orthogonal with the WF $\lower3pt\hbox{\rlap{$\scriptscriptstyle\rightarrow$}}\not {\rm w}_{\diamondsuit } \, [\eta ]\, $.  As expected the RDT in question decreases (increases) by 1 the ExpDiff for the pole at $\eta$ = ?1 ($\eta$ = +1).   

Substituting (5.15) into (4.32) we conclude that the double-indexed XR-Jacobi polynomials (5.14) form a finite orthogonal polynomial set:
\[\int _{1}^{\infty }J \, _{m_{{\bf a}'} +v-1}^{{\rm (}\stackrel{\rightharpoonup}{{\bf \lambda }}_{o} )} [\eta |{\bf a}'m_{{\bf a}'} ]\, \, {\rm J}_{m_{{\bf a}'} +\tilde{v}-1}^{{\rm (}\stackrel{\rightharpoonup}{{\bf \lambda }}_{o} )} [\eta |{\bf a}'m_{{\bf a}'} ]\, W[\eta ;\stackrel{\rightharpoonup}{{\bf \lambda }}_{o} |{\bf a}'m_{{\bf a}'} ]d\eta =0   (5.16)\] 

       for $0\le v<\tilde{v}<{\raise0.7ex\hbox{$ 1 $}\!\mathord{\left/ {\vphantom {1 2}} \right. \kern-\nulldelimiterspace}\!\lower0.7ex\hbox{$ 2 $}} (\lambda _{o;-} -\lambda _{o;+} -1)<\lambda _{o;-} -1\, $  

\noindent with the weight function
\[W[\eta ;\stackrel{\rightharpoonup}{{\bf \lambda }}_{o} |{\bf a}'m_{{\bf a}'} ]\equiv \frac{(\eta +1)^{1-\lambda _{o;-} } (\eta -1)^{\lambda _{o;+} +1} }{[P_{m_{{\bf a}'} }^{(-\lambda _{o;-} ,\lambda _{o;+} )} (\eta )]^{2} }      (5.17)\] 
Again, selected range (5.3${}^{o}$) for v is equivalent to condition (5.8${}^{o}$) for the square integral of 

\noindent polynomials (5.14) to converge:
\[\int _{1}^{\infty }( {\rm J}_{m_{{\bf a}'} +v-1}^{{\rm (}\stackrel{\rightharpoonup}{{\bf \lambda }}_{o} )} [\eta |{\bf a}'m_{{\bf a}'} ])\, ^{2} \, W[\eta ;\stackrel{\rightharpoonup}{{\bf \lambda }}_{o} |{\bf a}'m_{{\bf a}'} ]\, d\eta <\infty .     (5.18)\] 
It directly follows from (2.21), coupled with the conventional formula for the first derivatives of Jacobi polynomials, that the sequence of XR-Jacobi polynomials (5.14) starts from the monic Jacobi polynomial
\[{\rm P}_{m_{{\bf a}'} -1}^{{\rm (}\stackrel{\rightharpoonup}{{\bf \lambda }}_{o} )} [\eta |-\, +{\bf ,}m_{{\bf a}'} ;{\bf -}\, {\bf +,}0]=P\, _{m_{{\bf a}'} }^{(1-\lambda _{o;-} ,\lambda _{o;+} +1)} (\eta )/k\, _{m_{{\bf a}'} }^{(1-\lambda _{o;-} ,\lambda _{o;+} +1)}   (5.19)\] 

All of the mentioned $\mathscr{J}$S solutions \textbf{a}'m\textbf{${}_{a}$}${}_{'}$${}_{\ }$can be added${}_{\ }$to the set${}_{\ }$$\mathrm{\{}$\textbf{a}m\textbf{${}_{a}$}${}_{\}}$${}_{\ }$${}_{{}_{p_{{\bf a}} } }$which results in an extended RDCN $\mathrm{\{}$\textbf{a}m\textbf{${}_{a}$}${}_{\}}$${}_{{}_{p_{{\bf a}} } }$$\mathrm{\{}$\textbf{a}'m\textbf{${}_{a}$}${}_{'}$${}_{\}}$${}_{{}_{p_{{\bf a}'} } }$of  R-Jacobi admissible $\mathscr{J}$S solutions${}_{\ }$vanishing at the 

\noindent singular endpoint $\eta$ = +1.  While each RDT in the net increases by 1 the ExpDiff for the pole at $\eta$ = +1 the RDTs using seed solutions \textbf{a}'m\textbf{${}_{a}$}${}_{'}$${}_{\ }$decreases the ExpDiff for the pole at $\eta$ = ?1 until it 

\noindent becomes smaller than 1.  After that this ExpDiff oscillates between the values 

\noindent $\lambda _{o;-} -\, \left\lfloor \lambda _{o;-} \right\rfloor $ and $\left\lceil \lambda _{o;-} \right\rceil -\lambda _{o;-} \, $.

The finite subset \textbf{b}m\textbf{${}_{b}$} of nodeless $\mathscr{J}$S solutions of type  \textbf{b} co-existent with the discrete spectrum is defined by the condition
\[\lower3pt\hbox{\rlap{$\scriptscriptstyle\rightarrow$}}{\rm \varepsilon }_{{\bf b}{\rm m}_{{\bf b}} } {\bf (}\stackrel{\rightharpoonup}{{\bf \lambda }}_{{\bf o }} {\bf )-}\lower3pt\hbox{\rlap{$\scriptscriptstyle\rightarrow$}}{\bf \varepsilon }_{{\bf c}{\rm 0}} {\bf (}\stackrel{\rightharpoonup}{{\bf \lambda }}_{{\bf o }} {\bf )}=\, (\lambda _{o;+} -\lambda _{o;-} +1)^{2} -(-\lambda _{o;+} -\lambda _{o;-} +2{\rm m}_{{\bf b}} +1)^{2} \] 
\[=\, -4(\lambda _{o;-} -1-{\rm m}_{{\bf b}} )(\lambda _{o;+} -{\rm m}_{{\bf b}} )<0    (5.20)\] 
so

\noindent $0\le {\rm m}_{{\bf b}} <\lambda _{o;+} <\lambda _{o;-} -1$ $\mathrm{<}$${\raise0.7ex\hbox{$ 1 $}\!\mathord{\left/ {\vphantom {1 2}} \right. \kern-\nulldelimiterspace}\!\lower0.7ex\hbox{$ 2 $}} (\lambda _{o;-} +\lambda _{o;+} -1)$.     (5.20${}^{o}$)

\noindent Based on the analysis presented in previous section one can use FFs from this subset only if 

\noindent $\lambda _{o;+} >{\raise0.7ex\hbox{$ 1 $}\!\mathord{\left/ {\vphantom {1 2}} \right. \kern-\nulldelimiterspace}\!\lower0.7ex\hbox{$ 2 $}} .$   Since, according to Table 5.1, $\stackrel{\rightharpoonup}{{\bf \sigma }}_{1} =--$ for any primary $\mathscr{J}$S solution of this type we 

\noindent come to double indexed XB-Jacobi polynomials (2.8-1ii) with $\stackrel{\rightharpoonup}{{\rm \lambda }}\in \, {\rm {\rm I} {\rm I} }$:
\[{\rm J}_{{\rm m}_{{\bf b}} +v}^{{\rm (}\stackrel{\rightharpoonup}{{\bf \lambda }}_{o} )} [\eta |{\bf b}{\rm m}_{{\bf b}} {\rm ;}v]\equiv {\rm P}_{{\rm m}_{{\rm b}} +v}^{{\rm (}\stackrel{\rightharpoonup}{{\rm \lambda }}_{o} )} [\eta |-\, -{\rm ,m}_{{\rm b}} ;{\rm -+,}v]     (5.21)\] 
with v and m\textbf{${}_{b}$}${}_{\ }$restricted by constraints (5.3${}^{o}$) and (5.20${}^{o}$), respectively.  Substituting (2.8-1ii)  coupled with (5.21), into the right hand side of (4.31) we come to the third finite set of quasi-

\noindent rational solutions of prime SLE (4.1*),

\noindent 
\[\begin{array}{l} {\lower3pt\hbox{\rlap{$\scriptscriptstyle\rightarrow$}}\not \psi _{v} {}_{v} [\eta ;\mathop{\lambda }\limits^{\rightharpoonup}{} _{o} |{\bf b}m_{{\bf b}} ]\equiv \lower3pt\hbox{\rlap{$\scriptscriptstyle\rightarrow$}}\not \psi _{v} {}_{v} [\eta ;\mathop{\lambda }\limits^{\rightharpoonup}{} _{o} |{\bf -}\, {\bf +,}m_{{\bf b}} ]} \\ {\, \, \, \, \, \, \, \, \, \, \, \, \, \, \, \, \, \, \, \, \propto \, (1+\eta )^{-{\raise0.7ex\hbox{$ 1 $}\!\mathord{\left/ {\vphantom {1 2}} \right. \kern-\nulldelimiterspace}\!\lower0.7ex\hbox{$ 2 $}} \lambda _{o;-} +1} (1-\eta )^{{\raise0.7ex\hbox{$ 1 $}\!\mathord{\left/ {\vphantom {1 2}} \right. \kern-\nulldelimiterspace}\!\lower0.7ex\hbox{$ 2 $}} (\lambda _{o;+} -1)} \frac{{\rm J}_{m_{{\bf b}} +v}^{{\rm (}\stackrel{\rightharpoonup}{{\bf \lambda }}_{o} )} [\eta |{\bf b}m_{{\bf b}} ]}{P_{n}^{(-\lambda _{o;+} ,-\lambda _{o;-} )} (\eta )} ,} \end{array}   (5.22)\] 
orthogonal with the WF ${}_{1} \not {\rm w}_{\diamondsuit } \, [\eta ]\, $.  Again, as expected the RDT in question decreases by 1 the ExpDiffs for the poles at $\eta$=$\pm$1.   

Substituting (5.22) into (4.32) we conclude that the double-indexed XR-Jacobi polynomials (5.21) form a finite orthogonal polynomial set:
\[\int _{1}^{\infty }J \, _{m_{{\bf b}} +v}^{{\rm (}\stackrel{\rightharpoonup}{{\bf \lambda }}_{o} )} [\eta |{\bf b}m_{{\bf b}} ]\, {\rm J}_{m_{{\bf b}} +\tilde{v}}^{{\rm (}\stackrel{\rightharpoonup}{{\bf \lambda }}_{o} )} [\eta |{\bf b}m_{{\bf b}} ]\, W[\eta ;\stackrel{\rightharpoonup}{{\bf \lambda }}_{o} |{\bf b}m_{{\bf b}} ]d\eta =0    (5.23)\] 

       for $0\le v<\tilde{v}<{\raise0.7ex\hbox{$ 1 $}\!\mathord{\left/ {\vphantom {1 2}} \right. \kern-\nulldelimiterspace}\!\lower0.7ex\hbox{$ 2 $}} (\lambda _{o;-} -\lambda _{o;+} -1)<\lambda _{o;-} -1\, $  (5.23${}^{o}$)

\noindent with the weight function
\[W[\eta ;\stackrel{\rightharpoonup}{{\bf \lambda }}_{o} |{\bf b}m_{{\bf b}} ]\equiv \frac{(\eta +1)^{-\lambda _{o;-} +1} (\eta -1)^{\lambda _{o;+} -1} }{[P_{m_{{\bf b}} }^{(-\lambda _{o;-} ,-\lambda _{o;+} )} (\eta )]^{2} }  .    (5.24)\] 
Similarly to WFs (5.8) and (5.17) selected range (5.3${}^{o}$) for v is equivalent to condition (5.8${}^{o}$) for the square integral of polynomials (5.21) to converge:
\[\int _{1}^{\infty }( {\rm J}_{m_{{\bf b}} +v}^{{\rm (}\stackrel{\rightharpoonup}{{\bf \lambda }}_{o} )} [\eta |{\bf b}m_{{\bf b}} ])\, ^{2} \, W[\eta ;\stackrel{\rightharpoonup}{{\bf \lambda }}_{o} |{\bf b}m_{{\bf b}} ]\, d\eta <\infty .     (5.25)\] 

Setting $\stackrel{\rightharpoonup}{{\bf \sigma }}_{1} =-\, -,\, \, \stackrel{\rightharpoonup}{{\bf \sigma }}_{1}^{\dag } =-+$, $m=m_{{\bf b}} $, and n = v in (2.19) gives
\[\begin{array}{l} {{\rm P}_{m_{{\bf b}} +v}^{{\rm (}\lambda _{o;-} ,\lambda _{o;+} )} [\eta |+\, -{\bf ,}m_{{\bf b}} ;-+{\bf ,}v]\propto (\eta -1)P_{m_{{\bf b}} }^{(-\lambda _{o;+} ,-\lambda _{o;-} )} (\eta ){\mathop{P}\limits^{\bullet }} \, _{v}^{(\lambda _{o;+} ,-\lambda _{o;-} )} (\eta )} \\ {\, \, \, \, \, \, \, \, \, \, \, \, \, \, \, \, \, \, \, \, \, \, \, \, \, \, \, \, \, \, \, \, \, \, \, \, \, \, \, \, \, \, \, \, \, \, \, \, \, \, \, \, \, \, \, \, \, +(\lambda _{o;+} -m_{{\bf b}} )P_{m_{{\bf b}} }^{(-\lambda _{o;+} -1,1-\lambda _{o;-} )} (\eta )P_{v}^{(\lambda _{o;+} ,-\lambda _{o;-} )} (\eta )} \end{array}(5.26)\] 
so the sequence of XR-Jacobi polynomials (5.21) starts from the monic Jacobi polynomial 
\[{\rm P}_{m_{{\bf b}} }^{{\rm (}\lambda _{o;-} ,\lambda _{o;+} )} [\eta |+\, -{\bf ,}m_{{\bf b}} ;-+{\bf ,}0]\propto (\lambda _{o;+} -m_{{\bf b}} )P_{m_{{\bf b}} }^{(-\lambda _{o;+} -1,1-\lambda _{o;-} )} (\eta ).  (5.27)\] 

On other hand, substituting (2.9*), with $\lambda _{\pm } =-\lambda _{o;\pm } ,\, $ into the right-hand side of (5.10) and setting $m_{{\bf b}} $to 0 gives
\[{\rm P}_{v}^{{\rm (}\lambda _{o;-} ,\lambda _{o;+} )} [\eta |+\, -{\bf ,}0;-+{\bf ,}v]=P_{v}^{(\lambda _{o;+} -1,1-\lambda _{o;-} )} (\eta )/k_{v}^{(\lambda _{o;+} -1,1-\lambda _{o;-} )}   (5.28)\] 
As expected, the RDT with the basic FF ${\bf b}0$ simply decreases by 1 the ExpDiffs for the singular point $\eta =\pm 1$, without creating new poles.   Making use of (4.14) one finds that
\[W[\eta ;\stackrel{\rightharpoonup}{{\bf \lambda }}_{o} |{\bf b}{\rm 0}]\equiv 2^{\lambda _{o;+} -\lambda _{o;-} } \lower3pt\hbox{\rlap{$\scriptscriptstyle\rightarrow$}}W({\raise0.7ex\hbox{$ 1 $}\!\mathord{\left/ {\vphantom {1 2}} \right. \kern-\nulldelimiterspace}\!\lower0.7ex\hbox{$ 2 $}} \eta -{\raise0.7ex\hbox{$ 1 $}\!\mathord{\left/ {\vphantom {1 2}} \right. \kern-\nulldelimiterspace}\!\lower0.7ex\hbox{$ 2 $}} ;\lambda _{o;-} -1,\lambda _{o;+} -1)   (5.29)\] 
in agreement with (5.28).

One can use all three sets of quasi-rational PFSs of JRef  CSLE (A.1) with no zeros between 1 and $\mathrm{\infty}$ to construct the DCN of orthogonal polynomial sequences

\noindent ${\rm J}_{\tilde{n}_{\{ {\bf t}m\} _{p} } }^{{\rm (}\stackrel{\rightharpoonup}{{\bf \lambda }}_{o} )} [\eta |$$\{ {\bf t}m\} _{p} ;{\bf c}v]$   v= 0, 1, {\dots} , $\left\lfloor {\raise0.7ex\hbox{$ 1 $}\!\mathord{\left/ {\vphantom {1 2}} \right. \kern-\nulldelimiterspace}\!\lower0.7ex\hbox{$ 2 $}} (\lambda _{o;-} -\lambda _{o;+} -1)\right\rfloor $,   (5.30${}^{\dagger }$) 

\noindent where
\[\{ {\bf t}m\} _{p} \equiv{\rm \{ }{\bf a}n\} _{p_{{\bf a}} } ;{\rm \{ }{\bf a}'m_{{\bf a}'} \} _{p_{{\bf a}'} } {\rm \{ }{\bf b}m_{{\bf b}} \} _{p_{{\bf b}} }       (5.30)\] 
with 

\noindent p = $p_{{\bf a}} $+$p_{{\bf a}'} $+$p_{{\bf b}} $+1.         (5.30')

\noindent The only constraint on the choice of the quasi-rational PFSs in question is that the ExpDiffs for the persistent singularities $\eta$ = $\pm$1 of the CSLE generated using $p_{{\bf b}} $-1 $\mathscr{J}$S solutions ${\bf b}m_{{\bf b}} $ must be 

\noindent larger than ${\raise0.7ex\hbox{$ 1 $}\!\mathord{\left/ {\vphantom {1 2}} \right. \kern-\nulldelimiterspace}\!\lower0.7ex\hbox{$ 2 $}} $:

\noindent 

\noindent $\lambda _{o;-} $+$p_{{\bf a}} $?$p_{{\bf a}'} $?$p_{{\bf b}} $$\mathrm{>}$${\raise0.7ex\hbox{$ 3 $}\!\mathord{\left/ {\vphantom {3 2}} \right. \kern-\nulldelimiterspace}\!\lower0.7ex\hbox{$ 2 $}} $, $\lambda _{o;+} $+$p_{{\bf a}} $?$p_{{\bf a}'} $?$p_{{\bf b}} $$\mathrm{>}$${\raise0.7ex\hbox{$ 3 $}\!\mathord{\left/ {\vphantom {3 2}} \right. \kern-\nulldelimiterspace}\!\lower0.7ex\hbox{$ 2 $}} $.     (5.30b)

\noindent If the PFSs at infinity are not used ($p_{{\bf b}} $= 0) then the number $p_{{\bf a}'} $ of $\mathscr{J}$S solutions ${\bf a}'m_{{\bf a}'} $is constrained by the single requirement that the ExpDiff for the persistent singularity $\eta$ = -1 of the CSLE generated using $p_{{\bf a}'} $-1 $\mathscr{J}$S solutions ${\bf a}'m_{{\bf a}'} $ must be larger than${\raise0.7ex\hbox{$ 1 $}\!\mathord{\left/ {\vphantom {1 2}} \right. \kern-\nulldelimiterspace}\!\lower0.7ex\hbox{$ 2 $}} $:

\noindent $\lambda _{o;-} $+$p_{{\bf a}} $?$p_{{\bf a}'} $$\mathrm{>}$${\raise0.7ex\hbox{$ 3 $}\!\mathord{\left/ {\vphantom {3 2}} \right. \kern-\nulldelimiterspace}\!\lower0.7ex\hbox{$ 2 $}} $.        (5.30a')

If there is $\mathscr{J}$S solutions ${\bf d}m$ or ${\bf d}'m$ with no zeros between 1 and $\mathrm{\infty}$ (case III in Quesne's terms [66]) then it can be used to extend the constructed DCN of orthogonal polynomial sequences (5.30${}^{\dagger }$)  . An analysis of Klein's [67] formulas (6.72.4)-(6.72.8) in [22] shows that the m-degree Jacobi polynomial with indexes $\lambda _{+} $ and $\lambda _{-} $do not have zero between +1 to +$\mathrm{\infty}$ iff [5]   
\[(\lambda _{+} +\lambda _{-} +m+1)_{m} (\lambda _{+} +1)_{m} >0    (5.31)\] 
and
\[u_{3} (\lambda _{+} ,\lambda _{-} ;m)\equiv {\raise0.7ex\hbox{$ 1 $}\!\mathord{\left/ {\vphantom {1 2}} \right. \kern-\nulldelimiterspace}\!\lower0.7ex\hbox{$ 2 $}} \left(-|2m+\lambda _{+} +\lambda _{-} +1|+\lambda _{o;+} -\lambda _{o;-} +1\right)\le 1,   (5.31')\] 
where $\lambda _{o;\pm } =\, |\lambda _{\pm } |$.  Since the $\mathscr{J}$S solution in question ir required to co-exist with the discrete spectrum condition (5.31') is automatically fulfilled.  For any $\mathscr{J}$S solution \textbf{d}m\textbf{${}_{\ }$}formed by the               m-degree Jacobi polynomial with indexes $\lambda _{+} =-\lambda _{o;+} $and $\lambda _{-} =\lambda _{o;-} $ the rising factorial $(\lambda _{o;-} -\lambda _{o;+} +m+1)_{m} $ is the product of positive number as far as the SLP in question has discrete spectrum.  The second rising factorial also becomes positive for any even m if we 

\noindent choose the ExpDiff $\lambda _{o;+} $ larger than m:
\[(1-\lambda _{o;+} )_{m} =(-1)^{m} (\lambda _{o;+} -m)_{m} >0\, \, \, \, for\, \, m\, =2j_{{\bf d}} <\lambda _{o;+} .   (5.32)\] 
The inequality
\[\begin{array}{l} {\lower3pt\hbox{\rlap{$\scriptscriptstyle\rightarrow$}}\varepsilon \, _{{\bf c}0} \, -\, \lower3pt\hbox{\rlap{$\scriptscriptstyle\rightarrow$}}\varepsilon \, _{{\bf d}{\rm m}} =(\lambda _{o;-} -\lambda _{o;+} +2m+1)^{2} -(\lambda _{o;-} -\lambda _{o;+} -1)^{2} } \\ {\, \, \, \, \, \, \, \, \, \, \, \, \, \, \, \, \, \, \, \, \, \, \, \, =4(m+1)(\lambda _{o;-} -\lambda _{o;+} +m)>0} \end{array}   (5.33)\] 
shows that all the $\mathscr{J}$S\textbf{ }solutions \textbf{d}m ( including the R-Jacobi admissible solutions \textbf{d}m\textbf{${}_{d}$} with        m\textbf{${}_{d}$${}_{\ }$}$\mathrm{<}$$\lambda _{o;+} $) lie below the lowest eigenvalue $\lower3pt\hbox{\rlap{$\scriptscriptstyle\rightarrow$}}\varepsilon \, _{{\bf c}0} $, in agreement with the predictions of the 

\noindent conventional SUSY theory.

\section*{6. Cross-orthogonality relations for X${}_{m}$-Jacobi polynomials from different OPSs using             R-Jacobi polynomials as seed solutions}

\noindent We finally came to the problem specified by the title of the paper: the infinite biorthogonal 

\noindent polynomial system composed of X-Jacobi polynomials.  As pointed to in previous section the 

\noindent XR-Jacobi polynomial ${\rm J}_{n+m}^{{\rm (}\stackrel{\rightharpoonup}{{\bf \lambda }}_{o} )} [\eta |{\rm an}]$ for
\[0\le m<\lambda _{o;-} -1         (6.1o)\] 
coincides with the n${}^{th}$ X${}_{m}$ --Jacobi polynomial of degree n after the latter is written in its monic form and therefore must be orthogonal to any other polynomial from this sequence 
\[\int _{-1}^{1} \stackrel{\frown}{W}_{\lambda _{o;-} \, -1,\lambda _{o;+} +1,m} (-\eta )\, \, {\rm J}_{n+m}^{{\rm (}\stackrel{\rightharpoonup}{{\bf \lambda }}_{o} )} [\eta |{\rm an}]\, {\rm J}_{\tilde{n}+m}^{{\rm (}\stackrel{\rightharpoonup}{{\rm \lambda }}_{o} )} [\eta |{\rm an}]\, d\eta =0 (\tilde{{\rm n}}>n\ge 0) (6.1)\] 
  with the WF given by (88) in [43]
\[\stackrel{\frown}{W}_{\alpha ,\beta ,m} (x)\equiv \frac{(1-x)^{\alpha } (1+x)^{\beta } }{[P_{m}^{(-\alpha -1,\, \beta -1)} (x)]^{2} } .       (6.2)\] 
Here we are only interested in case (B) in [43]:
\[\alpha =\lambda _{o;-} \, -1>m-1,\, \, \, \beta =\lambda _{o;+} +1>0.       (6.3)\] 
Therefore, as far as m lies within range (6.1${}^{o}$) the XR-Jacobi polynomial in question has exactly m exceptional zeros between 1 and +$\mathrm{\infty}$, in addition to n conventional zeros between ?1 and +1.  

Alternatively we can reformulate this statement by saying that X --Jacobi polynomials from 

\noindent the OPSs using R-Jacobi polynomials as seed functions satisfy the cross-orthogonality relation:
\[\int _{-\infty }^{-1}\stackrel{\frown}{P} \, _{m,m+n}^{(\alpha ,\beta )} (x)\stackrel{\frown}{P}_{\tilde{m},\tilde{m}+n}^{(\alpha ,\beta )} (x)\, {\mathop{W}\limits^{\leftrightarrow }} _{\alpha ,\beta ;n} [\eta ]\, \, dx=0.     (6.4) \] 
for $0\le m<\tilde{m}<{\raise0.7ex\hbox{$ 1 $}\!\mathord{\left/ {\vphantom {1 2}} \right. \kern-\nulldelimiterspace}\!\lower0.7ex\hbox{$ 2 $}} (\alpha -\beta +1)$    (6.4${}^{o}$)

\noindent with the WF defined via (5.8):
\[{\mathop{W}\limits^{\leftrightarrow }} _{\alpha ,\beta ;n} [\eta ]\equiv W[-\eta ;\, \alpha +1,\beta -1|{\rm a}n]\equiv \frac{(1-\eta )^{-\alpha -2} (\eta +1)^{\beta } }{[P_{n}^{(\beta -1,\alpha +1)} (\eta )]^{2} } .    (6.5)\] 
Any X${}_{m}$-Jacobi polynomial of degree m + n with m smaller than ${\raise0.7ex\hbox{$ 1 $}\!\mathord{\left/ {\vphantom {1 2}} \right. \kern-\nulldelimiterspace}\!\lower0.7ex\hbox{$ 2 $}} (\alpha -\beta +1)$ must thereby have exactly m negative exceptional zeros with absolute values larger than 1, in addition to n 

\noindent conventional zeros between ?1 and +1.  Indeed, examination of the first two X${}_{1}$-Jacobi polynomials of degree 1 and 2 in [31] shows that they have only negative zeros, in agreement with the above assertion.

It has been proven in [43] that all m exceptional zeros of the X${}_{m}$-Jacobi polynomial 

\noindent $\, \stackrel{\frown}{P}_{m,m+n}^{(\alpha ,\beta )} (x)$ outside the interval [?1, +1] tend to zeros of the seed Jacobi polynomial of degree m which do lie inside the infinite negative interval (-$\mathrm{\infty}$, -1) for m $\mathrm{<}$${\raise0.7ex\hbox{$ 1 $}\!\mathord{\left/ {\vphantom {1 2}} \right. \kern-\nulldelimiterspace}\!\lower0.7ex\hbox{$ 2 $}} (\alpha -\beta +1)$.  A more detailed comparison between our results for exceptional zeros of X${}_{m}$-Jacobi polynomials and 

\noindent predictions based on the so-called `electrostatic' properties of their zeros [68-70] will be done in a separate study.

To complete this analysis it seems informative to outline the route leading from real 

\noindent XB-Jacobi polynomials to X-Jacobi OPSs in framework of a more general formalism developed 

\noindent in [16].  Based on our technique we come to X${}_{m}$-Jacobi OPSs by converting RCSLEs (A.1) and 

\noindent (A.1*) to the prime self-adjoint SLEs 
\[\left\{\, \frac{d\, \, }{d\eta } (1-\eta ^{2} )\frac{d\, \, }{d\eta } -\not q\, [\eta ;\pm 1;\mathop{\lambda }\limits^{\rightharpoonup}{} _{o} ]\, +{\raise0.7ex\hbox{$ 1 $}\!\mathord{\left/ {\vphantom {1 4}} \right. \kern-\nulldelimiterspace}\!\lower0.7ex\hbox{$ 4 $}} \, \varepsilon \right\}\not \Psi [\eta ;\varepsilon ;\mathop{\lambda }\limits^{\rightharpoonup}{} _{o} ]\, =\, 0    (6.6)\] 
and 
\[\left\{\, \frac{d\, \, }{d\eta } (1-\eta ^{2} )\frac{d\, \, }{d\eta } -\not q\, [\eta ;\pm 1;\mathop{\lambda }\limits^{\rightharpoonup}{} _{o} |+\, -;m]\, +{\raise0.7ex\hbox{$ 1 $}\!\mathord{\left/ {\vphantom {1 4}} \right. \kern-\nulldelimiterspace}\!\lower0.7ex\hbox{$ 4 $}} \, \varepsilon \right\}\not \Psi [\eta ;\varepsilon ;\mathop{\lambda }\limits^{\rightharpoonup}{} _{o} |+\, -;m]\, =\, 0  (6.6*)\] 
solved under the DBCs
\[\not \Psi [\pm 1;\varepsilon _{v} ;\mathop{\lambda }\limits^{\rightharpoonup}{} _{o} ]\, =\, 0         (6.7)\] 
and
\[\not \Psi [\pm 1;\varepsilon _{v} ;\mathop{\lambda }\limits^{\rightharpoonup}{} _{o} |-\, +;m]\, =\, 0.        (6.7*)\] 
To our knowledge the SLPs of this type we first introduced in by Gomez-Ullate, Kamran, and Milson in their revolutionary study [31] on orthogonal X${}_{1}$-Jacobi polynomials.  Two SLEs (which happened to be eigenequations in this particular case) share the same discrete spectrum 
\[\, {\rm \varepsilon }_{{\rm ++,}\, {\rm v}} {\rm (}\stackrel{\rightharpoonup}{{\rm \lambda }}_{{\rm o}} {\rm )}=(\lambda _{o;-} +\lambda _{o;+} +2v+1)^{2}        (6.8)\] 
independent of the choice of the $\mathscr{J}$S solution $-\, +;m$. 

\noindent Since the RDTs in question reflects the ExpDiff for the singular endpoint $\eta =+1$:
\[\lambda _{+1;+} =1-\lambda _{o;+} >0         (6.9)\] 
if $\lambda _{o;+} <1$,  the RD$\mathscr{T}$ of the eigenfunction $+\, +;v$of eigenequation (6.6) \textit{necessarily} satisfies DBCs (6.7*) only for $\lambda _{o;+} >{\raise0.7ex\hbox{$ 1 $}\!\mathord{\left/ {\vphantom {1 2}} \right. \kern-\nulldelimiterspace}\!\lower0.7ex\hbox{$ 2 $}} $. (There can be cases when this is also true for $0<\lambda _{o;+} <{\raise0.7ex\hbox{$ 1 $}\!\mathord{\left/ {\vphantom {1 2}} \right. \kern-\nulldelimiterspace}\!\lower0.7ex\hbox{$ 2 $}} $ but we currently skip the examination of this issue since these solutions cannot be automatically 

\noindent used as seed functions for extending the RDCN of solvable SLPs of our interest.)   Note that conditions (6.3) for the solution +?,m to have zeros between -1 and +1 is applicable for any positive value of $\beta$ which implies that our technique does not cover the range of the index $\beta$ between 0 and ${\raise0.7ex\hbox{$ 3 $}\!\mathord{\left/ {\vphantom {3 2}} \right. \kern-\nulldelimiterspace}\!\lower0.7ex\hbox{$ 2 $}} $.

\noindent 

\section*{7.  Quantization of the rationally deformed h-PT potentials by double-indexed XR-Jacobi polynomials}

\noindent One comes to the Schrödinger equation with the h-PT potential
\[\, V_{h-PT} (r;\stackrel{\rightharpoonup}{{\bf \lambda }}_{o} )\equiv \lower3pt\hbox{\rlap{$\scriptscriptstyle\rightarrow$}}V^{{\rm L}} [\eta (r);\stackrel{\rightharpoonup}{{\rm \lambda }}_{o} ]=\, \frac{{\raise0.7ex\hbox{$ 1 $}\!\mathord{\left/ {\vphantom {1 2}} \right. \kern-\nulldelimiterspace}\!\lower0.7ex\hbox{$ 2 $}} -2\lambda _{o;-}^{{\rm 2}} }{2{\it ch}^{2} r} +\frac{2\lambda _{{\rm o ;+}}^{{\rm 2}} -{\raise0.7ex\hbox{$ 1 $}\!\mathord{\left/ {\vphantom {1 2}} \right. \kern-\nulldelimiterspace}\!\lower0.7ex\hbox{$ 2 $}} }{2{\it sh}^{2} r}     (7. 1)\] 
by applying the Liouville transformation [71, 72, 37, 19] to JRef CSLE (A.1) on the infinite interval $+1<\eta <+\infty $, i.e., by representing the former equation in the Schrödinger form 
\[\left\{\, \frac{d\, \, }{d\eta } \lower3pt\hbox{\rlap{$\scriptscriptstyle\rightarrow$}}{\it \rho }_{\, \diamondsuit }^{-{\raise0.7ex\hbox{$ 1 $}\!\mathord{\left/ {\vphantom {1 2}} \right. \kern-\nulldelimiterspace}\!\lower0.7ex\hbox{$ 2 $}} } [\eta ]\frac{d\, \, }{d\eta } -\lower3pt\hbox{\rlap{$\scriptscriptstyle\rightarrow$}}{\it \rho }_{\, \diamondsuit }^{{\raise0.7ex\hbox{$ 1 $}\!\mathord{\left/ {\vphantom {1 2}} \right. \kern-\nulldelimiterspace}\!\lower0.7ex\hbox{$ 2 $}} } [\eta ]\left(\lower3pt\hbox{\rlap{$\scriptscriptstyle\rightarrow$}}V^{{\rm L}} [\eta ;\stackrel{\rightharpoonup}{{\bf \lambda }}_{o} ]-\, \, \lower3pt\hbox{\rlap{$\scriptscriptstyle\rightarrow$}}\varepsilon \, \right)\, \, \right\}\lower3pt\hbox{\rlap{$\scriptscriptstyle\rightarrow$}}\Psi ^{{\bf S}} [\eta ;\lower3pt\hbox{\rlap{$\scriptscriptstyle\rightarrow$}}\varepsilon ;\stackrel{\rightharpoonup}{{\bf \lambda }}_{o} ]=\, 0,   (7.2)\] 
with the algebraic Liouville potential [16]
\[\lower3pt\hbox{\rlap{$\scriptscriptstyle\rightarrow$}}V^{{\rm L}} [\eta ;\stackrel{\rightharpoonup}{{\bf \lambda }}_{o} ]=-\lower3pt\hbox{\rlap{$\scriptscriptstyle\rightarrow$}}{\it \rho }_{\, \diamondsuit }^{-1} [\eta ]\, I^{o} [\eta ;\mathop{\lambda }\limits^{\rightharpoonup}{} _{o} ]+\lower3pt\hbox{\rlap{$\scriptscriptstyle\rightarrow$}}{\it \rho }_{\, \diamondsuit }^{-{\raise0.7ex\hbox{$ 1 $}\!\mathord{\left/ {\vphantom {1 2}} \right. \kern-\nulldelimiterspace}\!\lower0.7ex\hbox{$ 2 $}} } [\eta ]{\it I}\{ \lower3pt\hbox{\rlap{$\scriptscriptstyle\rightarrow$}}{\it \rho }_{\, \diamondsuit }^{-{\raise0.7ex\hbox{$ 1 $}\!\mathord{\left/ {\vphantom {1 2}} \right. \kern-\nulldelimiterspace}\!\lower0.7ex\hbox{$ 2 $}} } [\eta ]\}     (7.3)\] 
and then converting (7.2) to the conventional Schrödinger equation via the change of variable
\[\lower3pt\hbox{\rlap{$\scriptscriptstyle\rightarrow$}}\eta (r)={\it ch}{\rm (}2r)           (7.4)\] 
determined by the first-order differential equation [12]
\[\frac{d\lower3pt\hbox{\rlap{$\scriptscriptstyle\rightarrow$}}\eta }{d\, r} =\lower3pt\hbox{\rlap{$\scriptscriptstyle\rightarrow$}}{\rm p}_{{\rm S}} [\lower3pt\hbox{\rlap{$\scriptscriptstyle\rightarrow$}}\eta (r)].         (7.4')\] 
with
\[\lower3pt\hbox{\rlap{$\scriptscriptstyle\rightarrow$}}{\rm p}_{{\rm S}} [\lower3pt\hbox{\rlap{$\scriptscriptstyle\rightarrow$}}\eta ]=\sqrt{\eta ^{2} -1} ={\raise0.7ex\hbox{$ 1 $}\!\mathord{\left/ {\vphantom {1 2}} \right. \kern-\nulldelimiterspace}\!\lower0.7ex\hbox{$ 2 $}} \lower3pt\hbox{\rlap{$\scriptscriptstyle\rightarrow$}}{\it \rho }_{\, \diamondsuit }^{-{\raise0.7ex\hbox{$ 1 $}\!\mathord{\left/ {\vphantom {1 2}} \right. \kern-\nulldelimiterspace}\!\lower0.7ex\hbox{$ 2 $}} } [\eta ]       (7.5)\] 
standing for the LCF of self-adjoint SLE (7.2).   The second summand in the right-hand side of (7.3),
\[{\rm p}_{{\rm S}} [\eta ]\, \tilde{{\rm I}}\{ {\rm p}_{{\rm S}} [\eta ]\} \equiv {\raise0.7ex\hbox{$ 1 $}\!\mathord{\left/ {\vphantom {1 4}} \right. \kern-\nulldelimiterspace}\!\lower0.7ex\hbox{$ 4 $}} {\mathop{{\rm p}}\limits^{\bullet }} \, _{{\rm S}}^{2} \, [\eta ]-{\raise0.7ex\hbox{$ 1 $}\!\mathord{\left/ {\vphantom {1 2}} \right. \kern-\nulldelimiterspace}\!\lower0.7ex\hbox{$ 2 $}} {\rm p}_{{\rm S}} [\eta ]{\mathop{{\rm p}}\limits^{\bullet \bullet }} \, _{{\rm S}} [\eta ],     (7.6)\] 
turns into the conventional Schwarzian derivative $\mathrm{\{}$$\eta$,r$\mathrm{\}}$ when converted to r.  Substituting (7.5) into (7.6), one finds
\[\lower3pt\hbox{\rlap{$\scriptscriptstyle\rightarrow$}}{\rm p}_{{\rm S}} [\lower3pt\hbox{\rlap{$\scriptscriptstyle\rightarrow$}}\eta ]\tilde{{\rm I}}\{ \lower3pt\hbox{\rlap{$\scriptscriptstyle\rightarrow$}}{\rm p}_{{\rm S}} [\eta ]\} =1+\frac{3}{\eta ^{2} -1}         (7.7)\] 
so
\[\lower3pt\hbox{\rlap{$\scriptscriptstyle\rightarrow$}}V^{{\rm L}} [\eta ;\stackrel{\rightharpoonup}{{\bf \lambda }}_{o} ]=\, \frac{{\raise0.7ex\hbox{$ 1 $}\!\mathord{\left/ {\vphantom {1 2}} \right. \kern-\nulldelimiterspace}\!\lower0.7ex\hbox{$ 2 $}} -2\lambda _{o;-}^{{\rm 2}} }{\eta +1} +\frac{2\lambda _{{\rm o ;+}}^{{\rm 2}} -{\raise0.7ex\hbox{$ 1 $}\!\mathord{\left/ {\vphantom {1 2}} \right. \kern-\nulldelimiterspace}\!\lower0.7ex\hbox{$ 2 $}} }{\eta -1} .      (7.8)\] 
which immediately brings us to potential function (8.1).  In terms of our earlier works [5, 63]${}^{\ }$

$\lambda _{o;+} =\lambda _{o} ,\, \, \lambda _{o;-} =\mu _{o} $,        (7.8p)

\noindent whereas the parameters used by Odake and Sasaki [2] are related to ExpDiffs (7.9) as follows

${\rm g}=\lambda _{o;+} +{\raise0.7ex\hbox{$ 1 $}\!\mathord{\left/ {\vphantom {1 2}} \right. \kern-\nulldelimiterspace}\!\lower0.7ex\hbox{$ 2 $}} $, ${\rm h}=\lambda _{o;-} -{\raise0.7ex\hbox{$ 1 $}\!\mathord{\left/ {\vphantom {1 2}} \right. \kern-\nulldelimiterspace}\!\lower0.7ex\hbox{$ 2 $}} .$      (7.8p')

\noindent (There is a misprint in (3.11b) in [5].)  Note that, compared with (38) in [2], we shifted the energy-reference point (ERefP): 
\[V_{h-PT} (r;{\rm h},{\rm g})\equiv -\frac{{\rm h}({\rm h}+1)}{{\it cosh}^{2} r} +\frac{{\it g}({\it g}-1)}{{\it sinh}^{2} r} ,       (7.8')\] 
to make the potential vanish at infinity, in following Pöschl and Teller's  original definition [1] of this potential.  Combining (4.12) with [16]
\[\lower3pt\hbox{\rlap{$\scriptscriptstyle\rightarrow$}}\psi ^{{\rm S}} [\eta ;\stackrel{\rightharpoonup}{{\bf \lambda }}_{o} |\, \stackrel{\rightharpoonup}{{\bf \sigma }}{\bf ,}m]\propto \, \, \sqrt[{4}]{\lower3pt\hbox{\rlap{$\scriptscriptstyle\rightarrow$}}\not {\bf p}[\eta ]\lower3pt\hbox{\rlap{$\scriptscriptstyle\rightarrow$}}\not {\bf w}[\eta ]} \, \lower3pt\hbox{\rlap{$\scriptscriptstyle\rightarrow$}}\not \psi [\eta ;\stackrel{\rightharpoonup}{{\bf \lambda }}_{o} |\, \stackrel{\rightharpoonup}{{\bf \sigma }}{\bf ,}m]     (7.9)\] 
\[\propto \, \, \sqrt[{4}]{\frac{\eta -1}{\eta +1} } \not \psi [\eta ;\stackrel{\rightharpoonup}{{\bf \lambda }}_{o} |\, \stackrel{\rightharpoonup}{{\bf \sigma }}{\bf ,}m]\, \,       (7.9')\] 
gives
\[\psi _{v}^{{\rm S}} [{\it cosh}(\, 2r);{\rm h}+{\raise0.7ex\hbox{$ 1 $}\!\mathord{\left/ {\vphantom {1 2}} \right. \kern-\nulldelimiterspace}\!\lower0.7ex\hbox{$ 2 $}} ,{\rm g}-{\raise0.7ex\hbox{$ 1 $}\!\mathord{\left/ {\vphantom {1 2}} \right. \kern-\nulldelimiterspace}\!\lower0.7ex\hbox{$ 2 $}} ]\propto \, \frac{{\rm sinh}^{{\rm g}} r}{{\rm cosh}^{{\rm h}} r} P_{v}^{({\rm g}-{\raise0.7ex\hbox{$ 1 $}\!\mathord{\left/ {\vphantom {1 2}} \right. \kern-\nulldelimiterspace}\!\lower0.7ex\hbox{$ 2 $}} ,-{\rm h}-{\raise0.7ex\hbox{$ 1 $}\!\mathord{\left/ {\vphantom {1 2}} \right. \kern-\nulldelimiterspace}\!\lower0.7ex\hbox{$ 2 $}} )} \left({\rm cosh}(\, 2r)\right).   (7.10)\] 
in agreement with (39) and (41) in [2]. 

Making use of nodeless $\mathscr{J}$S solutions of three types introduced in Section 5 we come to three families of rationally deformed h-PT potentials
\[\lower3pt\hbox{\rlap{$\scriptscriptstyle\rightarrow$}}V^{{\rm L}} [\eta ;\stackrel{\rightharpoonup}{{\bf \lambda }}_{o} |{\bf t}m]=-\lower3pt\hbox{\rlap{$\scriptscriptstyle\rightarrow$}}{\it \rho }_{\, \diamondsuit }^{-1} [\eta ]\, \, I^{o} [\eta ;\mathop{\lambda }\limits^{\rightharpoonup}{} _{o} |{\it t}m]+\lower3pt\hbox{\rlap{$\scriptscriptstyle\rightarrow$}}{\it \rho }_{\, \diamondsuit }^{-{\raise0.7ex\hbox{$ 1 $}\!\mathord{\left/ {\vphantom {1 2}} \right. \kern-\nulldelimiterspace}\!\lower0.7ex\hbox{$ 2 $}} } [\eta ]{\it I}\{ \lower3pt\hbox{\rlap{$\scriptscriptstyle\rightarrow$}}{\it \rho }_{\, \diamondsuit }^{-{\raise0.7ex\hbox{$ 1 $}\!\mathord{\left/ {\vphantom {1 2}} \right. \kern-\nulldelimiterspace}\!\lower0.7ex\hbox{$ 2 $}} } [\eta ]\} ,   (7.11)\] 
or alternatively
\[\lower3pt\hbox{\rlap{$\scriptscriptstyle\rightarrow$}}V^{{\rm L}} [\eta ;\stackrel{\rightharpoonup}{{\bf \lambda }}_{o} |{\bf t}m]=\lower3pt\hbox{\rlap{$\scriptscriptstyle\rightarrow$}}V^{{\bf L}} [\eta ;\stackrel{\rightharpoonup}{{\bf \lambda }}_{o} ]-\lower3pt\hbox{\rlap{$\scriptscriptstyle\rightarrow$}}{\it \rho }_{\, \diamondsuit }^{-1} [\eta ]\, \, \left(I^{o} [\eta ;\mathop{\lambda }\limits^{\rightharpoonup}{} _{o} |{\it t}m]-\, I^{o} [\eta ;\mathop{\lambda }\limits^{\rightharpoonup}{} _{o} ]\right),   (7.11*)\] 
where

${\rm m}\, \equiv \, n\ge 1$$\mathrm{>}$ 0       for ${\bf t}$=${\bf a}$;     (7.11a)

${\rm m}\, \equiv \, {\rm m}_{{\bf a}'} >\lambda _{o;-} -1-\lambda _{o;+} >0$  for ${\bf t}$=${\bf a}'$;     (7.11a')

$1\le m\equiv {\rm m}_{{\bf b}} <\lambda _{o;+} <\lambda _{o;-} -1$  for ${\bf t}$=${\bf b}$.     (7.11b).

\noindent Substituting (5.4), (5.15), and (5.22) into (7.9') we come to the following expressions for the eigenfunctions of the Schrödinger equation with the potential $\lower3pt\hbox{\rlap{$\scriptscriptstyle\rightarrow$}}V^{{\rm L}} [\eta ;\stackrel{\rightharpoonup}{{\bf \lambda }}_{o} |{\bf t}m]$:

$\psi _{v}^{{\rm S}} [{\it cosh}\, (2r);{\rm h}+{\raise0.7ex\hbox{$ 1 $}\!\mathord{\left/ {\vphantom {1 2}} \right. \kern-\nulldelimiterspace}\!\lower0.7ex\hbox{$ 2 $}} ,{\rm g}-{\raise0.7ex\hbox{$ 1 $}\!\mathord{\left/ {\vphantom {1 2}} \right. \kern-\nulldelimiterspace}\!\lower0.7ex\hbox{$ 2 $}} |\, {\bf a}n]\, \propto \, \frac{{\bf sinh}^{{\bf g+1}} r}{{\bf cosh}^{{\bf h+1}} r} \frac{\stackrel{\frown}{{\rm P}}_{v,n+v}^{{\rm (h}+{\raise0.7ex\hbox{$ 1 $}\!\mathord{\left/ {\vphantom {1 2}} \right. \kern-\nulldelimiterspace}\!\lower0.7ex\hbox{$ 2 $}} ,{\bf g}-{\raise0.7ex\hbox{$ 1 $}\!\mathord{\left/ {\vphantom {1 2}} \right. \kern-\nulldelimiterspace}\!\lower0.7ex\hbox{$ 2 $}} )} \left({\rm -cosh}(\, 2r)\right)}{P_{n}^{({\rm g}-{\raise0.7ex\hbox{$ 1 $}\!\mathord{\left/ {\vphantom {1 2}} \right. \kern-\nulldelimiterspace}\!\lower0.7ex\hbox{$ 2 $}} ,{\rm h}+{\raise0.7ex\hbox{$ 1 $}\!\mathord{\left/ {\vphantom {1 2}} \right. \kern-\nulldelimiterspace}\!\lower0.7ex\hbox{$ 2 $}} )} \left({\rm cosh}(\, 2r)\right)} ,$  (7.12a)

$\psi _{v}^{{\rm S}} [{\it cosh}\, (2r);{\rm h}+{\raise0.7ex\hbox{$ 1 $}\!\mathord{\left/ {\vphantom {1 2}} \right. \kern-\nulldelimiterspace}\!\lower0.7ex\hbox{$ 2 $}} ,{\rm g}-{\raise0.7ex\hbox{$ 1 $}\!\mathord{\left/ {\vphantom {1 2}} \right. \kern-\nulldelimiterspace}\!\lower0.7ex\hbox{$ 2 $}} |\, {\bf a}'m_{{\bf a}'} ]\, \propto \, \frac{{\bf sinh}^{{\bf g+1}} r}{{\bf cosh}^{{\bf h-1}} r} \frac{{\rm J}_{m_{{\bf a}'} +v-1}^{({\bf g}-{\raise0.7ex\hbox{$ 1 $}\!\mathord{\left/ {\vphantom {1 2}} \right. \kern-\nulldelimiterspace}\!\lower0.7ex\hbox{$ 2 $}} ,{\bf h}+{\raise0.7ex\hbox{$ 1 $}\!\mathord{\left/ {\vphantom {1 2}} \right. \kern-\nulldelimiterspace}\!\lower0.7ex\hbox{$ 2 $}} )} {\rm [}{\bf cosh}(2r)|{\bf a}'m_{{\bf a}'} ]}{P_{n}^{({\bf g}-{\raise0.7ex\hbox{$ 1 $}\!\mathord{\left/ {\vphantom {1 2}} \right. \kern-\nulldelimiterspace}\!\lower0.7ex\hbox{$ 2 $}} ,-{\bf h}-{\raise0.7ex\hbox{$ 1 $}\!\mathord{\left/ {\vphantom {1 2}} \right. \kern-\nulldelimiterspace}\!\lower0.7ex\hbox{$ 2 $}} )} \left({\bf cosh}(2r)\right)} ,$(7.12a')

\noindent and

$\psi _{v}^{{\rm S}} [{\it cosh}\, (2r);{\rm h}+{\raise0.7ex\hbox{$ 1 $}\!\mathord{\left/ {\vphantom {1 2}} \right. \kern-\nulldelimiterspace}\!\lower0.7ex\hbox{$ 2 $}} ,{\rm g}-{\raise0.7ex\hbox{$ 1 $}\!\mathord{\left/ {\vphantom {1 2}} \right. \kern-\nulldelimiterspace}\!\lower0.7ex\hbox{$ 2 $}} |\, {\bf b}m_{{\bf b}} ]\, \propto \, \frac{{\bf sinh}^{{\bf g-1}} r}{{\bf cosh}^{{\bf h-1}} r} \frac{{\rm J}_{m_{{\bf b}} +v}^{{\rm (g}-{\raise0.7ex\hbox{$ 1 $}\!\mathord{\left/ {\vphantom {1 2}} \right. \kern-\nulldelimiterspace}\!\lower0.7ex\hbox{$ 2 $}} ,{\bf h}+{\raise0.7ex\hbox{$ 1 $}\!\mathord{\left/ {\vphantom {1 2}} \right. \kern-\nulldelimiterspace}\!\lower0.7ex\hbox{$ 2 $}} )} [\eta |{\bf b}m_{{\bf b}} ]}{P_{n}^{({\raise0.7ex\hbox{$ 1 $}\!\mathord{\left/ {\vphantom {1 2}} \right. \kern-\nulldelimiterspace}\!\lower0.7ex\hbox{$ 2 $}} -{\bf g},{\raise0.7ex\hbox{$ 1 $}\!\mathord{\left/ {\vphantom {1 2}} \right. \kern-\nulldelimiterspace}\!\lower0.7ex\hbox{$ 2 $}} -{\bf h})} (\eta )} $  (7.12b)

\noindent for ${\bf t}$= ${\bf a}{\rm ,}\, \, {\bf a}'{\rm ,}$ and ${\bf b}$ accordingly.

Odake and Sasaki [2] discussed only the branch $\lower3pt\hbox{\rlap{$\scriptscriptstyle\rightarrow$}}V^{{\rm L}} [{\it cosh}\, (2r);\stackrel{\rightharpoonup}{{\bf \lambda }}_{o} |{\bf b}m_{{\bf b}} ]$.  It should be stressed in this connection that the XR-Jacobi polynomials appearing in the right-hand side of (7.12b) essentially differ from X-Jacobi polynomials forming OPSs.   To term them X${}_{m}$-Jacobi polynomials as done in [73] seems confusing even if one uses this term in a broader sense.  The reader should remember that we deal with finite subsets of double-indexed XB-Jacobi polynomials specified by two indexes v and  $m_{{\bf b}} $ in our notation.

 Even in case of the Bagchi-Quesne-Roychoudhury (BQR) potential [74] generated by the first-degree polynomial
\[\Pi _{1} [\eta ;\eta _{{\bf t}1} (\stackrel{\rightharpoonup}{{\bf \lambda }}_{o} )]\equiv \eta -\eta _{{\bf t}1} (\stackrel{\rightharpoonup}{{\bf \lambda }}_{o} ),       (7.13)\] 
where
\[\eta _{{\bf t}1} (\lambda _{{\rm 0;t1}} ,\lambda _{{\rm 1;t1}} )=\frac{\lambda _{{\rm 0;t1}} -\lambda _{{\rm 1;t1}} }{\lambda _{{\rm 0;t1}} +\lambda _{{\rm 1;t1}} +2} <1,      (7.13')\] 
the XR-Jacobi polynomials appearing in the right-hand side of (7.12a) coincide with orthogonal X${}_{m}$-Jacobi polynomials with indexes m = n generally larger than 1 so the statement in [75] that eigenfunctions of the BQR potential are expressible in terms of  X${}_{1}$-Jacobi polynomials is incorrect even if one considers an extension of three solvable families of rationally deformed          h-PT potential using rational PFr (7.11*) parametrized  for m = 1 by three independent parameters $\lambda _{{\rm 0;}{\bf t}{\rm 1}} \, (\, |\lambda _{{\rm 0;t1}} |\, =\lambda _{o;-} ),\, \, \lambda _{{\rm 1;t1}} \, (|\lambda _{{\rm 1;t1}} |=\lambda _{o;+} ),$ and $\eta _{{\bf t}1} $.  The latter is conditionally exactly solvable via XR-Jacobi polynomials if the third parameter is constrained by condition (7.13') which gives

\noindent $\lambda _{{\rm 1;}{\bf a}{\rm 1}} =\, \lambda _{o;+} ,\, \, \, \lambda _{{\rm 0;a1}} =\lambda _{o;-} ;\, $         (7.14a)

\noindent $\lambda _{{\rm 1;}{\bf b}{\rm 1}} =\, -\lambda _{o;+} <-1,\, \, \, \lambda _{{\rm 0;b1}} =-\lambda _{o;-} ,\, $$\lambda _{o;-} +\lambda _{o;+} >3$;    (7.14b)

\noindent $\lambda _{{\rm 1;}{\bf d}'{\rm 1}} =\, -\lambda _{o;+} <-1,\, \, \lambda _{{\rm 0;}{\bf d}'{\rm 1}} =-\lambda _{o;-} ,\, \, \, and\, \, \, 2<\lambda _{o;-} +\lambda _{o;+} <3;$   (7.14d')

 $\lambda _{{\rm 1;}{\bf d}\, {\rm 1}} =\, -\lambda _{o;+} ,\, \, \lambda _{{\rm 0;d}\, {\rm 1}} =\lambda _{o;-} ,\, \, and\, \, \, {\raise0.7ex\hbox{$ 1 $}\!\mathord{\left/ {\vphantom {1 2}} \right. \kern-\nulldelimiterspace}\!\lower0.7ex\hbox{$ 2 $}} <\lambda _{o;+} <1;$     (7.14d)

\noindent and

$\lambda _{{\rm 1;}{\bf d}\, {\rm 1*}} =\, -\lambda _{o;+} ,\, \, \lambda _{{\rm 0;d}\, {\rm 1*}} =\lambda _{o;-} ,\, \, and\, \, \, 0<\lambda _{o;-} <\lambda _{o;+} -2.$    (7.14d*)

\noindent For completeness we also included the nodeless $\mathscr{J}$S solutions ${\bf d}'{\rm 1}$ and ${\bf d}\, {\rm 1*}$ which themselves do not co-exist with the discrete spectrum of the h-PT potential but lead to the single-energy-level BQR potentials when used as the FFs for the RDT.

\noindent 

\noindent \textbf{\eject }

\section*{8. Conclusions and further developments}

\noindent As illuminated by the author separately [40, 41], reduction of the \textit{complex} JRef eigenequation of Bochner type [30] to the real field results in two \textit{distinct }second-order differential equations solvable by polynomials.   Namely, in addition to its conventional reduction to the real Jacobi equation, the mentioned reduction leads to the second less-known second-order differential equation solvable by Routh polynomials [25].

Both real Jacobi polynomials with non-zero indexes and (real-by-definition) Routh polynomials can be used as seed functions to construct two RDCNs of SLEs conditionally exactly solvable by multi-index XB-Jacobi and XB-Routh polynomials.  The indexes of all $\mathscr{J}$S polynomials used to generate the corresponding RDCN may differ only by sign so each sequence of XB-Jacobi polynomials is specified by four positive integers which determine the numbers of positive/ negative first and second Jacobi indexes in the selected set of $\mathscr{J}$S solutions. The Routh-seed (RS) polynomials can be conveniently specified by a single complex index [20].   With this convention, the sets of RS solutions are formed by Routh polynomials with complex index equal to $\alpha$ or -$\alpha$ so each sequence of XB-Routh polynomials is specified by two positive integers which represent the numbers of complex indexes with positive and negative real parts in the selected set of RS solutions.

The next step is to formulate the SLP which allows one to select infinite or finite mutually orthogonal subsets from the constructed sequences of XB-Jacobi and XB-Routh polynomials.  Making use the general technique developed by us in [16] we impose the DBCs on solutions of the \textit{prime} SLE, namely, such that sum of two ChExps is equal to zero at each endpoint pole.  As a result of this very specific choice of the self-adjoint SLE the DBC unambiguously selects the PFS for any real (positive by definition) value of the corresponding ExpDiff whether it is the LP or LC case. 

As far as the RCSLEs of our interest and therefore the corresponding prime SLEs (which turned out be rational for X-Jacobi OPSs and for XB-Jacobi polynomials but algebraic for XB-Routh polynomials) have energy-independent ExpDiffs the SLPs in question are solvable via orthogonal (finite or infinite) polynomial sequences.  If the SLP for RDC$\mathscr{T}$s of the real Jacobi equation is formulated on the finite interval [?1, +1] we come to the broadly discussed (infinite) X-Jacobi OPSs (see, i.g., [29] and references therein).

An extension of Gomez-Ullate, Kamran, and Milson's ideas [31] to the SLP on the infinite 

\noindent interval [+1, +$\mathrm{\infty}$) leads to the RDCN of prime SLEs generated using three distinct sets of 

\noindent Jacobi polynomials with no zeros between 1 and $\mathrm{\infty}$ (`R-Jacobi admissible' in following Duran's 

\noindent terminology [42]).  Two sequences of the admissible $\mathscr{J}$S solutions constructed in such a way, \textbf{a}n and \textbf{a}'m\textbf{${}_{a}$}${}_{'}$, are composed of infinitely many Frobenius solutions of the prime SLE which all 

\noindent satisfy the DBC at the singular endpoint +1.  The $\mathscr{J}$S solutions from the sequence \textbf{a}n are formed 

\noindent by classical Jacobi polynomials starting from the polynomial of degree 1.  On other hand, the lowest degree of Jacobi polynomials from the second infinite sequence \textbf{a}'m\textbf{${}_{a}$}${}_{'}$${}_{\ }$must be at least twice larger than the number of eigenfunctions, n\textbf{${}_{c}$}. The third sequence of the R-Jacobi 

\noindent admissible $\mathscr{J}$S solutions is composed of a finite number of Frobenius solutions of the prime SLE 

\noindent vanishing at infinity.  This is the sequence originally discovered by Odake and Sasaki [2].  To our knowledge, both infinite sequences of the R-Jacobi admissible $\mathscr{J}$S solutions were first brought to light in [5]. 

In this paper we restricted the discussion solely to the RDTS using the aforementioned three sequences of the admissible $\mathscr{J}$S solutions as FFs, while postponing examination the whole RDCN of the SLPs solvable via muli-index XR-Jacobi polynomials for a separate publication.  The   main reason for treating these first-generation RD$\mathscr{T}$s of R-Jacobi polynomials as a special issue is our observation that infinitely many finite sequences (5.3) of double-indexed XR-Jacobi polynomials using classical Jacobi polynomials as seed functions are composed by X-Jacobi polynomials from different OPSs.  As a result polynomials from the manifold combining the         X-Jacobi OPSs in question obey the cross-orthogonality relation when integrated from +1 to +$\mathrm{\infty}$.   As a corollary we assert that each X${}_{m}$-Jacobi polynomial of degree m + n has exactly m \textit{exceptional} zeros between --$\mathrm{\infty}$ and --1 as far as its indexes are restricted by the constraints defining the given sequence of XR-Jacobi polynomials.

In following [46, 47], one can also extend the RDCN $\mathrm{\{}$\textbf{a}m\textbf{${}_{a}$}${}_{\}}$${}_{{}_{p_{{\bf a}} } }$$\mathrm{\{}$\textbf{a}'m\textbf{${}_{a}$}${}_{'}$${}_{\}}$${}_{{}_{p_{{\bf a}'} } }$$\mathrm{\{}$\textbf{b}m\textbf{${}_{b}$}${}_{\}}$${}_{{}_{p_{{\bf b}} } }$of 

\noindent solvable SLPs by using as seed solutions  pairs of `juxtaposed' [48-50] eigenfunctions \textbf{c},v and 

\noindent \textbf{c},v+1 (v$\mathrm{>}$0), as it has been sketched by us in [40, 41].

Odake and Sasaki [46] (see also [47] for a similar conjecture) speculated that nodeless quasi-rational solutions infinite at both endpoints exist for all the shape-invariant potentials, including the h-PT potential.  The only restriction cited by these authors is that the polynomial forming the given solution has \textit{even} degree.  The mentioned solutions (when exist) can be used for extending the DC net of rationally deformed h-PT potentials quantized by multi-index XR-Jacobi polynomials.  

In particular we can then use any ${\bf d},2j_{{\bf d}} $ ($2j_{{\bf d}} <\lambda _{o;+} $) introduced at the end of Section 5 to start a new RDCN of SLPs solvable via multi-index XR-Jacobi polynomials.  We refer to the latter as a new net (rather than an extension of the one using only PFSs) because we have to reduce subsets of R-Jacobi admissible seed polynomials.  Indeed since the RDT in question inserts the new eigen-solution below the energy $\lower3pt\hbox{\rlap{$\scriptscriptstyle\rightarrow$}}\varepsilon \, _{{\bf c}0} $ we should require that
\[\, \lower3pt\hbox{\rlap{$\scriptscriptstyle\rightarrow$}}\varepsilon \, _{{\bf d}{\rm ,2j}_{{\bf d}} } -\lower3pt\hbox{\rlap{$\scriptscriptstyle\rightarrow$}}\varepsilon \, _{{\bf a}\tilde{n}_{{\bf a}} } \, =(2\tilde{n}_{{\bf a}} +1+\lambda _{o;-} +\lambda _{o;+} )^{2} -(\lambda _{o;-} -\lambda _{o;+} +4{\rm j}_{{\bf d}} +1)^{2} \] 
$=4(\lambda _{o;-} +\tilde{n}_{{\bf a}} +2{\rm j}_{{\bf d}} +1)(\lambda _{o;+} +\tilde{n}_{{\bf a}} -2{\rm j}_{{\bf d}} )>0$,    (8.1a)
\[\, \lower3pt\hbox{\rlap{$\scriptscriptstyle\rightarrow$}}\varepsilon \, _{{\bf d}{\rm ,2j}_{{\bf d}} } -\lower3pt\hbox{\rlap{$\scriptscriptstyle\rightarrow$}}\varepsilon \, _{{\bf a}'\tilde{m}_{{\bf a}'} } \, =(2\tilde{m}_{{\bf a}'} +1-\lambda _{o;-} +\lambda _{o;+} )^{2} -(\lambda _{o;-} -\lambda _{o;+} +4{\rm j}_{{\bf d}} +1)^{2} \] 
$=4(\tilde{m}_{{\bf a}'} +2{\rm j}_{{\bf d}} +1)(\tilde{m}_{{\bf a}'} +\lambda _{o;+} -2{\rm j}_{{\bf d}} -\lambda _{o;-} )>0$;   (8.1a')
\[\, \lower3pt\hbox{\rlap{$\scriptscriptstyle\rightarrow$}}\varepsilon \, _{{\bf d}{\rm ,2j}_{{\bf d}} } -\lower3pt\hbox{\rlap{$\scriptscriptstyle\rightarrow$}}\varepsilon \, _{{\bf b}\tilde{m}_{{\bf b}} } \, =(2\tilde{m}_{{\bf b}} +1-\lambda _{o;-} -\lambda _{o;+} )^{2} -(\lambda _{o;-} -\lambda _{o;+} +4{\rm j}_{{\bf d}} +1)^{2} \] 
$=4(\tilde{m}_{{\bf b}} +2{\rm j}_{{\bf d}} +1-\lambda _{o;+} )(\tilde{m}_{{\bf b}} -2{\rm j}_{{\bf d}} -\lambda _{o;-} )>0$.                         (8.1b)

\noindent An analysis of inequality (8.1a) reveals that inserting the new eigenfunction does not affect       

\noindent R-Jacobi admissibility of $\mathscr{J}$S solutions \textbf{a}n, i.e.,  any classical Jacobi polynomials of non-zero 

\noindent degree can be still used to form the nodeless FF for the RDT.  On other hand, degree of any Jacobi polynomials forming the sequence \textbf{a}'$\tilde{m}_{{\bf a}'} $ necessarily exceeds 2(n\textbf{${}_{c}$}$+{\rm j}_{{\bf d}} )$, or to be more precise

\noindent 

$\tilde{m}_{{\bf a}'} >\lambda _{o;-} -\lambda _{o;+} +2{\rm j}_{{\bf d}} $.        (8.2a')

\noindent Finally, since $\tilde{m}_{{\bf b}} <\lambda _{o;+} $ the ExpDiff $\lambda _{o;-} $ must be larger than $\tilde{m}_{{\bf b}} -2{\rm j}_{{\bf d}} $ so (8.1b) holds only if

$\lambda _{o;+} >\tilde{m}_{{\bf b}} +2{\rm j}_{{\bf d}} +1,\, \, \, \lambda _{o;-} >\lambda _{o;+} +1$.      (8.2b)

 One cannot however guarantee that any pair of nodeless solutions ${\bf d},2j_{{\bf d}} $ and ${\bf d},2j'_{{\bf d}} $ can be used for constructing an admissible FF for the double-step RDCT.  From author's point of view, the suggestion [46, 47] to use `case  III' solutions as seed functions for RDC$\mathscr{T}$s of shape-invariant potentials was put forward mainly based on the tremendous progress in the theory of X-Hermite polynomials (see [76] and references therein).  However the very specific feature of this RDCN is that the harmonic oscillator is a symmetric potential and therefore this is also true for any of its RDC$\mathscr{T}$s if one uses even seed solutions.  Excluding the harmonic oscillator itself, the only symmetric potential on Odake and Sasaki' list of `shape-invariant' potentials [46] is the hyperbolic-secant-squared (HSS) potential termed by them `exciton potential'.  (Ironically the sample of illustrative examples examined in [47] does not include this potential.)

One of remarkable features of the HSS potential is that it can be treated as the symmetric limiting case of both Rosen-Morse [77] and Gendenshtein [65] (`Scarf  II' in terms of  [78, 79]) potentials and as a result can be alternatively quantized either via classical Gegenbauer polynomials [22] \textit{with energy-dependent indexes} (as originally pointed to by Ginocchio [80]) or via even/odd R-Routh polynomials expressible in terms of Gegenbauer polynomials in imaginary argument [81]
\[I_{v}^{(\lambda _{o} +{\raise0.7ex\hbox{$ 1 $}\!\mathord{\left/ {\vphantom {1 2}} \right. \kern-\nulldelimiterspace}\!\lower0.7ex\hbox{$ 2 $}} )} (\eta )=({\it i}\, )^{v} v\, C_{v}^{({\raise0.7ex\hbox{$ 1 $}\!\mathord{\left/ {\vphantom {1 2}} \right. \kern-\nulldelimiterspace}\!\lower0.7ex\hbox{$ 2 $}} -\lambda _{o} )} {\rm [i}\, \eta {\rm ]}\, .       (8.3)\] 
We [20] refer to (8.3) as the `Majedjamei polynomials' to get credit to Majedjamei [82] who initially introduced this very specific case of R-Routh polynomials.  By using even Routh polynomials of degree $2j_{{\bf d}} $$\mathrm{>}$ 0 as seed solutions one comes to finite sequences of multi-indexed exceptional Majedjamei polynomials (`X-M  polynomials' for briefness)
\[{\rm M}_{n_{\{ 2j_{{\bf d}} \} } }^{{\rm (}\lambda _{o} )} [\eta |\{ {\bf d},2j_{{\bf d}} \} _{p_{{\bf d}} } ;{\bf c}v]={\rm P}_{n_{\{ 2j_{{\bf d}} \} } }^{{\rm (}\lambda _{o} ,\lambda _{o} )} [{\it i}\eta |\{ {\rm ++,}2j_{{\it d}} \} _{p_{{\it d}} } ;{\rm -}\, {\rm -,}v]   (8.4)\] 
with a positive index $\lambda _{o} $ and a nonnegative integer v labeling bound energy states in the HSS 

\noindent potential.  It should be stressed that only eigenfunctions with more than $p_{{\bf d}} $-1 zeros are expressible in terms of orthogonal polynomials (8.4). The first $p_{{\bf d}} $ eigenfunctions (with the eigenvalues coincident with the energies of seed solutions $\{ {\bf d},2j_{{\bf d}} \} _{p_{{\bf d}} } $) form a separate subset of mutually orthogonal quasi-rational functions.  One can also use pairs of juxtaposed eigenfunctions $\{ \tilde{v},\tilde{v}+1\} _{2j_{{\bf c}} } $ to extend the net of solvable symmetric (non-shape-invariant) rationally deformed HSS potentials which brings us to the generalized sequences of multi-index X-M polynomials:
\[{\rm M}_{n_{\{ 2j_{{\bf d}} \} ;\{ \tilde{v},\tilde{v}+1\} ;v} }^{{\rm (}\lambda _{o} )} [\eta |\{ {\bf d},2j_{{\bf d}} \} _{j} ;\{ \tilde{v},\tilde{v}+1\} _{p_{{\bf c}} } ;{\bf c}v]     (8.4*)\] 
\[={\rm P}_{n_{\{ 2j_{{\bf d}} \} ;\{ \tilde{v},\tilde{v}+1\} _{p_{{\bf c}} } } }^{{\rm (}\lambda _{o} ,\lambda _{o} )} [{\it i}\eta |\{ {\rm ++,}2j_{{\it d}} \} _{p_{{\it d}} } ;-\, -,\{ \tilde{v}_{k} ,\tilde{v}_{k} +1\} _{k=1,...,2j_{{\it c}} } ;{\rm -}\, {\rm -,}v]\] 

 for $\tilde{v}_{k+1} >\tilde{v}_{k} +1\ge 2,\, \, \, v=0,...,\tilde{v}_{k} -1,\tilde{v}_{k} +2,...,\tilde{v}_{j_{{\bf c}} } -1,\tilde{v}_{j_{{\bf c}} } +2,...,v_{\max } .$

\noindent This generalized RDCN of analytically solvable potentials presents the excellent example of the conjecture in [46, 47] that such a net can be generated for any shape-invariant GRef potential.  However existence of similar RDCNs for other shape-invariant GRef potentials has still to be independently confirmed.

\section*{Acknowledgements}

\noindent The author thanks Robert Milson for providing the links to some obscure properties of Jacobi polynomials used by him to prove the \textit{fundamentally important }expression relating the first derivative of a Jacobi polynomial to a superposition of two Jacobi polynomials of exactly the same degree.

\noindent 

\noindent \textbf{\eject }

\section*{Appendix}

\noindent \textbf{Analysis of ExpDiffs for persistent singularities of the transformed CSLE }

\noindent Let us consider the JRef CSLE 

$\, \left\{\frac{d^{2} \, \, }{d\eta ^{2} } +I^{o} [\eta ;\stackrel{\rightharpoonup}{{\bf \lambda }}_{o} ]+\lower3pt\hbox{\rlap{$\scriptscriptstyle\rightarrow$}}\varepsilon \, \lower3pt\hbox{\rlap{$\scriptscriptstyle\rightarrow$}}{\it \rho }_{{}_{{\rm \diamondsuit }} } [\eta ]\, \right\}\, \, \lower3pt\hbox{\rlap{$\scriptscriptstyle\rightarrow$}}\Phi [\eta ;\lower3pt\hbox{\rlap{$\scriptscriptstyle\rightarrow$}}\varepsilon ;\stackrel{\rightharpoonup}{{\rm \lambda }}_{o} \, ]=\, 0$      (A.1)

\noindent with the JRef polynomial fraction (PFr)

$I^{o} [\eta ;\mathop{\lambda }\limits^{\rightharpoonup}{} _{o} ]\, \, \equiv \frac{1}{2(1-\eta ^{2} )} \, \, \, \sum _{\aleph \, =-,+} \frac{1-\lambda _{o;\aleph }^{2} }{1+\aleph \eta } -\frac{1}{4(1-\eta ^{2} )} $     (A.2)

\noindent and the simple-pole density function 

$\lower3pt\hbox{\rlap{$\scriptscriptstyle\rightarrow$}}{\it \rho }_{{}_{{\rm \diamondsuit }} } [\eta ]=\, \, \frac{1}{4(\eta ^{2} -1)} >0$     (1 $\mathrm{<}$ $\eta$ $\mathrm{<}$ $\mathrm{\infty}$).      (A.2')

\noindent The ERefP for the JRef CSLE was chosen in such a way that the ExpDiff for the second-order pole at infinity vanish for $\epsilon$ = 0 or, to be more precise, we required that

\noindent ${\mathop{{\it lim}}\limits_{\eta \to \infty }} \{ \eta ^{2} \, I^{o} [\eta ;\mathop{\lambda }\limits^{\rightharpoonup}{} _{o} ]\} ={\raise0.7ex\hbox{$ 1 $}\!\mathord{\left/ {\vphantom {1 4}} \right. \kern-\nulldelimiterspace}\!\lower0.7ex\hbox{$ 4 $}} $.        (A.3)

\noindent This brought us to JRef PFr (A.2) which depends only on the ExpDiffs $\lambda _{o;\pm } $ for second-order poles at $\eta =\pm 1$.  As the direct corollary of such a choice of the ERefP we find that the energy of any $\mathscr{J}$S solution

\noindent $\phi [\eta ;\stackrel{\rightharpoonup}{{\bf \lambda }}_{o} |\, \stackrel{\rightharpoonup}{{\bf \sigma }}{\bf ,}m]\propto \, \, \, \prod _{\aleph =-,+}( 1+\aleph \eta )^{{\raise0.7ex\hbox{$ 1 $}\!\mathord{\left/ {\vphantom {1 2}} \right. \kern-\nulldelimiterspace}\!\lower0.7ex\hbox{$ 2 $}} (\sigma _{\aleph } \lambda _{o;\aleph } +1)} \Pi _{m} [\eta ;\eta _{\stackrel{\rightharpoonup}{{\bf \sigma }}{\rm ,m}} (\stackrel{\rightharpoonup}{{\rm \lambda }}_{o} )],$    (A.4)

\noindent where 

\noindent $\Pi _{m} [\eta ;\eta _{\stackrel{\rightharpoonup}{{\bf \sigma }}{\rm ,m}} (\stackrel{\rightharpoonup}{{\rm \lambda }}_{o} )]\equiv P_{m}^{(\sigma _{+} \lambda _{o;+} ,\sigma _{-} \lambda _{o;-} )} (\eta )/k_{m}^{(\sigma _{+} \lambda _{o;+} ,\sigma _{-} \lambda _{o;-} )} ,$   (A.5)

\noindent is given by the quadratic formula

$\lower3pt\hbox{\rlap{$\scriptscriptstyle\rightarrow$}}\varepsilon \, _{\stackrel{\rightharpoonup}{{\bf \sigma }}{\bf ,}m} {\rm (}\stackrel{\rightharpoonup}{{\rm \lambda }}_{{\rm o}} {\rm )}=1-(\sigma _{-} \lambda _{o;-} +\sigma _{+} \lambda _{o;+} +2m+1)^{2} $.     (A.6)

Since density function (A.2') has the second-order pole at infinity RDTs of our interest it must preserve the ExpDiff at infinity [24]

${\mathop{{\it lim}}\limits_{\eta \to \infty }} \{ \eta ^{2} \, I^{o} [\eta ;\stackrel{\rightharpoonup}{{\bf \lambda }}_{o} |\stackrel{\rightharpoonup}{{\bf \sigma }}{\bf ,}m]\} ={\mathop{{\bf lim}}\limits_{\eta \to \infty }} \{ \eta ^{2} \, I^{o} [\eta ;\stackrel{\rightharpoonup}{{\bf \lambda }}_{o} \} $,     (A.7)

\noindent where symbol $\stackrel{\rightharpoonup}{{\bf \sigma }}{\bf ,}m$ identifies Jacobi-seed solution (A.3) used as the factorization function (FF) to construct the given RD$\mathscr{T}$ (A.5) of JRef CSLE (A.1),  

\noindent 

$\, \left\{\frac{d^{2} \, \, }{d\eta ^{2} } +I^{o} [\eta ;\stackrel{\rightharpoonup}{{\bf \lambda }}_{o} |\stackrel{\rightharpoonup}{{\bf \sigma }}{\bf ,}m]+\lower3pt\hbox{\rlap{$\scriptscriptstyle\rightarrow$}}\varepsilon \, \lower3pt\hbox{\rlap{$\scriptscriptstyle\rightarrow$}}{\it \rho }_{{}_{{\rm \diamondsuit }} } [\eta ]\, \right\}\, \, \lower3pt\hbox{\rlap{$\scriptscriptstyle\rightarrow$}}\Phi [\eta ;\lower3pt\hbox{\rlap{$\scriptscriptstyle\rightarrow$}}\varepsilon ;\stackrel{\rightharpoonup}{{\rm \lambda }}_{o} |\stackrel{\rightharpoonup}{{\rm \sigma }}{\rm ,}m\, ]\, =\, 0$.    (A.8)

\noindent As a result the zero-energy ExpDiff for the second-order pole of the RCSLE (A.8) at infinity must be also equal to 0.  

On other hand, since density function (A.2') has simple poles at both persistent finite singular points, the RDT changes ExpDiffs for both poles [5, 64]. One can directly verify the latter assertion by making use of Rudyak and Zakhariev's reciprocal formula [82]${}^{\ }$

${\rm *}\phi [\eta ;\stackrel{\rightharpoonup}{{\bf \lambda }}_{o} |\, \stackrel{\rightharpoonup}{{\bf \sigma }}{\bf ,}m]\, =\, \, {\it \rho }\, _{{\rm \diamondsuit }}^{-{\raise0.7ex\hbox{$ 1 $}\!\mathord{\left/ {\vphantom {1 2}} \right. \kern-\nulldelimiterspace}\!\lower0.7ex\hbox{$ 2 $}} } [\eta ]/\phi [\eta ;\stackrel{\rightharpoonup}{{\rm \lambda }}_{o} |\, \stackrel{\rightharpoonup}{{\rm \sigma }}{\rm ,}m]$      (A.9)

$\propto \, \, \, \prod _{\aleph =-,+} (1+\aleph \eta )^{-{\raise0.7ex\hbox{$ 1 $}\!\mathord{\left/ {\vphantom {1 2}} \right. \kern-\nulldelimiterspace}\!\lower0.7ex\hbox{$ 2 $}} \sigma _{\aleph } \lambda _{o;\aleph } } /\Pi _{m} [\eta ;\eta _{\stackrel{\rightharpoonup}{{\bf \sigma }}{\rm ,m}} (\stackrel{\rightharpoonup}{{\rm \lambda }}_{o} )]$   (A.9*)

\noindent for the FF of the inverse RDT from  CSLE (A.8) back to the original CSLE (A.1). We can then 

\noindent directly express the ExpDiff  $\lambda _{\pm } (\lambda _{o;\pm } |\sigma _{\pm } )$for the second order pole at the persistent singular point  $\eta$ = $\pm$1 in terms of the corresponding ExpDiff  $\lambda _{o;\pm } $[5]:

\noindent $\lambda _{\pm } (\lambda _{o;\pm } |\sigma _{\pm } )=\, |\sigma _{\pm } \lambda _{o;\pm } +1|$.  (A.10)

In the current context it is convenient to represent the RefPFr for RCSLE (A.8) as [24]

\[I^{o} [\eta ;\stackrel{\rightharpoonup}{{\bf \lambda }}_{o} |\stackrel{\rightharpoonup}{{\bf \sigma }}{\bf ,}m]=\frac{1-\lambda _{-}^{2} (\lambda _{o;-} |\sigma _{-} )}{4(\eta +1)^{2} } \, \, +\frac{1-\lambda _{+}^{2} (\lambda _{o;+} |\sigma _{+} )}{4(1-\eta )^{2} } \, \, -\frac{\, \stackrel{\frown}{O}\, _{m}^{\downarrow \, } [\eta ;\stackrel{\rightharpoonup}{{\bf \lambda }}_{o} |\stackrel{\rightharpoonup}{{\bf \sigma }}{\bf ,}m]\, }{4(1-\eta ^{2} )\Pi _{m} [\eta ;\eta _{\stackrel{\rightharpoonup}{{\bf \sigma }}{\rm ,m}} (\stackrel{\rightharpoonup}{{\rm \lambda }}_{o} )]} \]

\noindent $+{\mathop{\Pi }\limits^{\bullet \bullet }} [\eta ;\bar{\eta }_{\stackrel{\rightharpoonup}{{\bf \sigma }}{\rm ,m}} ]/\Pi [\eta ;\bar{\eta }_{\stackrel{\rightharpoonup}{{\rm \sigma }}{\rm ,m}} ]-2{\mathop{\Pi }\limits^{\bullet }} \, ^{2} [\eta ;\bar{\eta }_{\stackrel{\rightharpoonup}{{\rm \sigma }}{\rm ,m}} ]/\Pi ^{2} [\eta ;\bar{\eta }_{\stackrel{\rightharpoonup}{{\rm \sigma }}{\rm ,m}} ],$  (A.11)

\noindent where

\noindent $\stackrel{\frown}{O}\, _{m}^{\downarrow \, } [\eta ;\stackrel{\rightharpoonup}{{\bf \lambda }}_{o} |\stackrel{\rightharpoonup}{{\bf \sigma }}{\bf ,}m]\, =-(\, \lower3pt\hbox{\rlap{$\scriptscriptstyle\rightarrow$}}\varepsilon _{\, \stackrel{\rightharpoonup}{{\bf \sigma }}{\bf ,}m} +2\sigma _{-} \sigma _{+} \lambda _{o;-} \lambda _{o;+} )\, \Pi _{m} [\eta ;\eta _{\stackrel{\rightharpoonup}{{\bf \sigma }}{\rm ,m}} (\stackrel{\rightharpoonup}{{\rm \lambda }}_{o} )]\; $   (A.12)

\noindent 
\[-4\, [(\eta +1)\sigma _{+} \lambda _{o;+} +(\eta -1)\sigma _{-} \lambda _{o;-} ]\, \, {\mathop{\Pi }\limits^{\bullet }} _{m} [\eta ;\eta _{\stackrel{\rightharpoonup}{{\bf \sigma }}{\rm ,m}} (\stackrel{\rightharpoonup}{{\rm \lambda }}_{o} )].\]

\noindent In particular, substituting (A.6) into the leading coefficient of polynomial (A.11), 

\noindent $\stackrel{\frown}{O}\, _{m;m}^{\downarrow \, } (\stackrel{\rightharpoonup}{{\bf \lambda }}_{o} |\stackrel{\rightharpoonup}{{\bf \sigma }}{\bf ,}m)=\, -\lower3pt\hbox{\rlap{$\scriptscriptstyle\rightarrow$}}\varepsilon _{\stackrel{\rightharpoonup}{{\bf \sigma }}{\bf ,}m} -2\sigma _{-} \sigma _{+} \lambda _{o;-} \lambda _{o;+} +4m(\sigma _{-} \lambda _{o;-} +\sigma _{+} \lambda _{o;+} +2)\; $, (A.12')

\noindent and making use of  (A.10) one can directly verify that  

\noindent ${\mathop{{\it lim}}\limits_{\eta \to \infty }} \{ \eta ^{2} \, I^{o} [\eta ;\stackrel{\rightharpoonup}{{\bf \lambda }}_{o} |\stackrel{\rightharpoonup}{{\bf \sigma }}{\bf ,}m]\} ={\raise0.7ex\hbox{$ 1 $}\!\mathord{\left/ {\vphantom {1 4}} \right. \kern-\nulldelimiterspace}\!\lower0.7ex\hbox{$ 4 $}} $,       (A.13)

\noindent in agreement with (A.3) and (A.7).

The quasi-rational seed solutions of the algebraic Schrödinger equation  [16]${}_{ }$

${}_{ }$$\left\{\, \sqrt{\eta ^{2} -1} \, \, \frac{d\, \, }{d\eta } \sqrt{\eta ^{2} -1} \frac{d\, \, }{d\eta } -\left(\lower3pt\hbox{\rlap{$\scriptscriptstyle\rightarrow$}}V^{{\rm L}} [\eta ;\stackrel{\rightharpoonup}{{\bf \lambda }}_{o} ]-\, \, \lower3pt\hbox{\rlap{$\scriptscriptstyle\rightarrow$}}\varepsilon _{\stackrel{\rightharpoonup}{{\bf \sigma }}{\bf ,}m} \, \right)\, \, \right\}\lower3pt\hbox{\rlap{$\scriptscriptstyle\rightarrow$}}\psi ^{{\bf S}} [\eta ;\stackrel{\rightharpoonup}{{\bf \lambda }}_{o} |\, \stackrel{\rightharpoonup}{{\bf \sigma }}{\bf ,}m]=0$   (A.14)

\noindent are related to $\mathscr{J}$S solutions (A.4) of JRef CSLE (A.1) via the generic formula

\noindent $\lower3pt\hbox{\rlap{$\scriptscriptstyle\rightarrow$}}\psi ^{{\rm S}} [\eta ;\stackrel{\rightharpoonup}{{\bf \lambda }}_{o} |\, \stackrel{\rightharpoonup}{{\bf \sigma }}{\bf ,}m]=\, \lower3pt\hbox{\rlap{$\scriptscriptstyle\rightarrow$}}{\it \rho }_{\diamondsuit }^{{\raise0.7ex\hbox{$ {\rm 1} $}\!\mathord{\left/ {\vphantom {{\rm 1} {\rm 4}}} \right. \kern-\nulldelimiterspace}\!\lower0.7ex\hbox{$ {\rm 4} $}} } [\eta ]\, \phi [\eta ;\stackrel{\rightharpoonup}{{\it \lambda }}_{o} |\, \stackrel{\rightharpoonup}{{\it \sigma }}{\it ,}m]$,      (A.15)

\noindent which can be also used to obtain  eigenfunctions of the rationally deformed h-PT potentials discussed in Section  7.

\section*{References}

\noindent  [1] G. Pöschl and F. Teller, ``Bemerkungen zur Quantenmechanik des anharmonischen Oszillators,'' \textit{Zs. Phys. }\textbf{83 }(1933), 143 -151

\noindent  [2] S. Odake and R. Sasaki, ``Infinitely many shape invariant potentials and new orthogonal polynomials,'' \textit{Phys. Lett. }\textbf{B679 }(2009), 414-417, \underbar{arXiv: 0906.0142}

\noindent  [3] C.-L. Ho, S. Odake and R. Sasaki, ``Properties of the Exceptional ((X\textit{${}_{l}$}) Laguerre and Jacobi Polynomials,''\textit{SIGMA} \textbf{7 }(2011), 107, 24 pages, arXiv: 0912.5447 

\noindent  [4] J. W. Dabrowska, A. Khare and U. P. Sukhatme, ``Explicit Wavefunctions for Shape-Invariant Potentials by Operator Techniques,'' \textit{J. Phys. }A\textit{ }\textbf{21 }(1988),\textbf{ }L195-L200 

\noindent  [5] G. Natanson, \textbf{``}Gauss-seed nets of Sturm-Liouville problems with energy-independent characteristic exponents and related sequences of exceptional orthogonal polynomials  I. Canonical Darboux transformations using almost-everywhere holomorphic factorization functions'' (2013) arXiv:1305.7453v1\underbar{}

\noindent  [6] V. Romanovsky, ``Sur la généralisation des courbes de Pearson,'' \textit{Atti del Congresso Intern. dei Matem.} (\textit{Bologna)}\textbf{ 6 }(1928), 107, 4 pages, \newline\underbar{mathunion.org/ICM/ICM1928.6/Main/icm1928.6.0107.0110.ocr.pdf}

\noindent  [7] V. I. Romanovski, ``Sur quelques classes nouvelles de polynomes orthogonaux,`` \textit{C. R.  Acad.  Sci. (Paris)} \textbf{188 }(1929), 1023-1025 

\noindent  [8] R. Askey, ``An integral of Ramanujan and orthogonal polynomials,'' \textit{J. Indian Math. Soc. } \textbf{51 }(1987), 27-36 

\noindent  [9] W. Lesky, ``Endliche und unendliche Systeme von kontinuierlichen klassichen Othogonalpolynomen,'' \textit{Z.  Angew.  Math.  Mech.} \textbf{76 }(1996), 181-184 

\noindent [10] W. Koepf and M. Masjed-Jamei, ``A generic polynomial solution for the differential equation of hypergeometric type and six sequences of orthogonal polynomials related to it,'' \textit{Integral Transforms and Special Functions} \textbf{17 }(2006), 559-576 researchgate.net/publication/36409317\textit{\underbar{}}

\noindent [11] S. Das and A. Swaminathan ``Higher order derivatives of R-Jacobi polynomials,'' \textit{AIP Conference Proceed.} \textbf{1739} (2016), 020058, 8 pages

\noindent [12] G. A. Natanzon, ``Study of the one-dimensional Schrödinger equation generated from the hypergeometric equation,'' \textit{Vestn.  Leningr. Univ.} No 10 (1971), 22-28 [see arxiv.org/PS\_cache/physics/pdf/9907/9907032v1.pdf\underbar{ for }English translation]

\noindent [13] N. I. Zhirnov and O. P. Shadrin, ``Calculations of Franck-Condon factors with Poschl-Teller 

\noindent wavefunctions. I. The probabilities of some vibrational transitions in the $\, D{}^{1} \Sigma _{u}^{+} -B{}^{1} \Pi _{g} $ band system of the He${}_{2}$ molecule,'' \textit{Opt. Spectrosc.}\textbf{ 24} (1968), 478-481 

\noindent [14] A. K. Bose, ``Solvable potentials,'' \textit{Phys. Lett.} \textbf{7 }(1963), 245-246 \textbf{}

\noindent [15] A. K. Bose, ``A class of solvable potentials,'' \textit{Nuov. Cim.  }\textbf{32 }(1964), 679-688

\underbar{link.springer.com/article/10.1007/BF02735890}

\noindent [16] G. Natanson, ``Darboux-Crum Nets of Sturm-Liouville Problems Solvable by Quasi-Rational Functions I. General Theory'' (2018) \newline researchgate.net/publication/XXXXXXXX

\noindent [17] G. Darboux, \textit{Leçons sur la théorie générale des surfaces et les applications géométriques du calcul infinitésimal}, Vol. 2 (Paris, Gauthier-Villars, 1915) pp. 210-215

\noindent [18] M. M. Crum, ``Associated Sturm-Liouville systems,'' \textit{Quart. J. Math. Oxford} (2) \textbf{6} (1955), 121-127 arXiv:physics/9908019v1\underbar{}

\noindent [19] R. Milson, ``Liouville transformation and exactly solvable Schrödinger equations,'' \textit{Int. J. Theor. Phys. }\textbf{37 }(1998) 1735-1752, arXiv:solv-int/9706007v2\underbar{}

\noindent [20] G. Natanson, ``Exact quantization of the Milson potential via Romanovski-Routh polynomials,'' (2015) arXiv:1310.0796v3\underbar{}

\noindent [21] E. Heine, \textit{Handbuch der Kugelfunctionen, Theorie und Anwendungen, }vol. 1 (Berlin, G. Reimer, 1878)  p. 172 

\noindent [22] G. Szego, \textit{Orthogonal Polynomials} (New York: American Mathematical  Society, 1959 ) p.150

\noindent [23] B. Shapiro, ``Algebro-geometric aspects of Heine-Stieltjes theory,'' \textit{J. London Math. Soc.} \textbf{83} (2011), 36--56, researchgate.net/publication/23678760

\noindent [24] G. Natanson, ``Single-Source Nets of Algebraically-Quantized Reflective Liouville Potentials on the Line I. Almost-Everywhere Holomorphic Solutions of Rational Canonical Sturm-Liouville Equations with Second-Order Poles,'' (2015)  arXiv:1503.04798v2\underbar{}

\noindent [25] E. J. Routh, ``On some properties of certain solutions of a differential equation of second order,'' \textit{Proc.  London Math. Soc. }\textbf{16 }(1884), 245-261 

\noindent [26] G. Natanson, ``Survey of Nodeless Regular Almost-Everywhere Holomorphic Solutions  for Exactly Solvable Gauss-Reference Liouville Potentials on the Line I. Subsets of Nodeless Jacobi-Seed Solutions Co-Existent with Discrete  Energy Spectrum,'' (2016) arXiv:1606.08758\underbar{}

\noindent [27] G. Natanson, ``Survey of Nodeless Regular Almost-Everywhere Holomorphic Solutions for Exactly Solvable Gauss-Reference Liouville Potentials on the Line II. Selection of Nodeless Laguerre-Seed Solutions for Two Branches of the Confluent Gauss-Reference Potential,'' (2016) researchgate.net/publication/313842091\underbar{}

\noindent [28] D. Gomez-Ullate, N. Kamran, and R. Milson, ``On orthogonal polynomials spanning a non-standard flag,'' Contemp Math. \textbf{563} (2012) 51-71 arXiv:1101.5584v5\underbar{}

\noindent [29] M. Garcia-Ferrero, D. Gomez-Ullate, and R. Milson, ``A Bochner type classification theorem for exceptional orthogonal polynomials,`` (2017)\underbar{  }arXiv:1603.04358v2

\noindent [30] S. Bochner, ``Über Sturm-Liouvillesche Polynomsysteme,'' Math. Z.  \textbf{29 }(1929), 730-736, eudml.org/doc/168099\underbar{}

\noindent [31] D. Gomez-Ullate, N. Kamran, and  R. Milson, ``An extension of Bochner's problem: exceptional invariant subspaces,'' \textit{J. Approx. Theory }\textbf{162} (2010), 987-1006, arXiv:0805.3376v3\underbar{}

\noindent [32] P. Hartman, \textit{Ordinary Differential Equations} (Bikh?user, Boston, 1982)

\noindent [33] W. T. Reid, \textit{Sturmian Theory for Ordinary Differential Equations} (Springer, N-Y 1981)

\noindent [34] W. T. Reid, \textit{Ordinary Differential Equations} (John Wiley \& Sons, New York, 1971).

\noindent [35] J. Weidmann, \textit{Spectral Theory of Ordinary Differential Operators,} Lecture Notes in Math. \textbf{1258} (Springer-Verlag, Berlin, 1987)

\noindent [36] H. Weyl, ``Ramifications, old and new, of the eigenvalue problem,`` Bull. Amer. Math. Soc.\textit{}

\noindent \textbf{56} (1950), 115-139 \newline\underbar{pdfs.semanticscholar.org/93bb/1a7ee8d3b4c869638f9cb7528802d70eb245.pdf}

\noindent [37] W. N. Everitt,  \textit{``A catalogue of Sturm-Liouville differential equations'',} in \textit{Sturm-Liouville Theory, Past and Present}, p. 271-331, Birkhäuser Verlag, Basel 2005 (edited by W.O. Amrein, A.M. Hinz and D.B. Pearson.) 

\noindent [38] P.B. Bailey, W.N. Everitt, and A. Zettl, ``\textbf{Algorithm 810: The SLEIGN2 Sturm-Liouville Code,'' }\textit{ACM Transactions on Mathematical Software (TOMS)} \textbf{27}, 143-192 (2001)\newline\underbar{researchgate.net/publication/265681928\_Algorithm\_810\_The\_SLEIGN2\_Sturm-Liouville\_code}

\noindent [39] A. Zettl, \textit{Sturm-Liouville Theory}, Mathematical Surveys and Monographs Vol. 121,

\noindent [40] G. Natanson, ``Single-source nets of Fuschian rational canonical Sturm-Liouville equations with common simple-poles density functions and related sequences of multi-indexed orthogonal Heine eigenpolynomials,'' Presentation at the 14th International Symposium on Orthogonal Polynomials, Special Functions and Applications / 3-7 July 2017 \newline researchgate.net/publication/317643100\underbar{}

\noindent [41] G. Natanson, ``Single-Source Nets of Fuschian Rational Canonical Sturm-Liouville Equations with Common Simple-Poles Density Functions and Related Sequences of Multi-indexed Orthogonal Eigenpolynomials'', \newline researchgate.net/publication/317643178\underbar{}

\noindent [42] A. J. Duran, ``Exceptional Hahn and Jacobi orthogonal polynomials,`` J. Approx. Theory

214 (2017), 9-48

\noindent [43] D. Gomez-Ullate, F. Marcellan, and R. Milson, ``Asymptotic and interlacing properties of zeros of exceptional\underbar{ }Jacobi and Laguerre polynomials,'' \textit{J. Math. Anal. Appl.} \textbf{399} (2013), 480-495, arXiv:1204.2282v1\underbar{}

\noindent [44] S. Wolfram, ``Jacobi polynomials,'' 2018 \newline functions.wolfram.com/PDF/JacobiP.pdf\underbar{}

\noindent [45] M. Abramowitz and I. A. Stegun, \textit{Handbook of Mathematical Functions} (Washington DC:\textbf{}

\noindent NBS, Applied Mathematics Series -- 55, 1972) \newline\underbar{files.eric.ed.gov/fulltext/ED250164.pdf}

\noindent [46] S. Odake and R. Sasaki, ``Krein--Adler transformations for shape-invariant potentials and  pseudo virtual states,'' \textit{J. Phys. }\textbf{A46 }(2013), 245201, 24 pages, arXiv:1212.6595v2\underbar{}

\noindent [47] D. Gomez-Ullate, Y. Grandati, and R. Milson, ``Extended Krein-Adler theorem for the translationally shape invariant potentials,'' \textit{J. Math. Phys. }\textbf{55 }(2014), 043510, 30 pages arXiv:1309.3756\underbar{}

\noindent [48] B. F. Samsonov, ``On the equivalence of the integral and the differential exact solution generation methods for the one-dimensional Schrodinger equation,'' \textit{J. Phys.} \textbf{A28} (1995), 6989-6998

\noindent [49] B. F. Samsonov, ``New features in supersymmetry breakdown in quantum mechanics,'' \textit{Mod. Phys. Lett.} \textbf{A11} (1996), 1563-1567

\noindent [50] V. G. Bagrov and B. F. Samsonov, ``Darboux transformation and elementary exact solutions of the Schrödinger equation,'' \textit{Pramana J. Phys. }\textbf{49} (1997), 563-580

\noindent  [51] M. Krein, ``On a continuous analogue of the Christoffel formula from the theory of orthogonal polynomials,'' \textit{Dokl. Akad. Nauk SSSR }\textbf{113 }(1957), 970-973 

\noindent [52] V. E. Adler, ``A modification of Crum's method,'' \textit{Theor. Math. Phys.} \textbf{101}, 1381-1386 (1994)

\noindent [53] S. Karlin and G. Szego, ``On certain determinants whose elements are orthogonal polynomials\textit{,'' J. Analyse Math.} \textbf{8} (1960/1961), 1--157

\noindent [54] J. Weidmann, ``Spectral Theory of Sturm-Liouville Operators Approximation by Regular Problems'', in \textit{Sturm-Liouville Theory, Past and Present, }ed. W.O. Amrein, A. M. Hinz and D.B. Pearson (Birkhäuser Verlag, Basel),~pp 75-98 (2005) researchgate.net/publication/226812971\underbar{}

\noindent [55] P. B. Bailey, W. N. Everitt, J. Weidmann, and A. Zettl, ''Regular approximations of singular Sturm-Liouville problems'', Res. Math. 23,~3 (1993) researchgate.net/publication/2610149

\noindent [56] A. Zettl, ``Sturm-Liouville problems'', in \textit{Spectral Theory and Computational Methods of Sturm-Liouville problems}, ed. D. Hinton and P. W.Schaefer (Marcel Dekker, Inc., New York,1997), pp. 1-104

\noindent [57] J. Weidmann, \textit{Linear Operators in Hilbert Spaces}  (Springer, New York, 1980) archive.org/stream/springer\_10.1007-978-1-4612-6027-1/10.1007-978-1-4612-6027-1\#page/n9/mode/2up

\noindent [58] M. P. Chen and H. M. Srivastava, Orthogonality relations and generating functions for Jacobi polynomials and related hypergeometric functions, Appl. Math. Comput. 68 (1995), 153--188. \newline sciencedirect.com/science/article/pii/009630039400092I\underbar{}

\noindent [59] G. Hetyei, ``Shifted Jacobi polynomials and Delannoy numbers,'' (2009) arXiv:0909.5512v2

\noindent [60] A. P. Raposo, H. J. Weber, D. E. Alvarez-Castillo, and M. Kirchbach. ``Romanovski polynomials in selected physics problems,'' \textit{Centr. Eur. J. Phys.} \textbf{5} (2007)\textit{,} 253-284  arXiv:0706.3897v1 

\noindent [61] D. E. Avarez-Castillo and M. Kirchbach, ``Exact spectrum and wave functions of the hyperbolic Scarf potential in terms of finite Romanovski polynomials,''  \textit{Rev. Mex. Fis.} \textbf{E53 }(2007), 143-154 \newline\underbar{researchgate.net/publication/242322566}

\noindent [62]\textbf{ }H. J. Weber, ``Connections between Romanovski and other polynomials,'' C. Eur. J.  Math. \textbf{5} (2007),\textbf{ }581-595 \textbf{arXiv:0706.3153v1}\underbar{}

\noindent [63] C. Quesne\textit{, }``Extending Romanovski polynomials in quantum mechanics,'' \textit{J. Math. Phys. }\textbf{54 }(2013)\textit{,} 122103, 15 pages arXiv:1308.2114v2

\noindent [64] G. Natanson, ``Breakup of SUSY quantum mechanics in the limit-circle region of the reflective Kratzer oscillator,'' (2014) arXiv:1405.2059\underbar{v1}

\noindent [65] L. E. Gendenshtein, ``Derivation of exact spectra of the Schrödinger equation by means of supersymmetry,'' \textit{JETP Lett.} \textbf{38} (1983),\textbf{ }356-359

\noindent [66] C. Quesne, ``Solvable Rational Potentials and Exceptional Orthogonal Polynomials in Supersymmetric Quantum Mechanics,'' \textit{SIGMA }\textbf{5 }(2009), 084, 24 pages   arXiv:0906.2331v3\underbar{}

\noindent [67] F. C. Klein, ``Über die Nullstellen der hypergeometrischenReihe,'' \textit{Math. Ann}. \textbf{37 }(1890), 573-590 

\noindent [68] D. K. Dimitrov and Y. C. Lun, ``Monotonicity, interlacing and electrostatic interpretation of zeros of exceptional Jacobi polynomials,'' \textit{J. Approx. Theory} \textbf{181} (2014) 18-29

\noindent [69] A. P. Horvath, ``The electrostatic properties of zeros of exceptional Laguerre and Jacobi polynomials and stable interpolation, \textit{J. Approx. Theory} \textbf{194} (2015), 87--107.

\noindent [70] Y Luo, ``Zeros of exceptional orthogonal polynomials and the maximum of the modulus of an energy function'' (2017) \underbar{arXiv:1707.05202}

\noindent [71] J\textit{. Liouville, ``}Sur le développment des fonctionsou parties de fonctionsen series,'' \textit{ J. Math. Pure Appl. }\textbf{\textit{2}}\textit{ (1837), 16-35 }

\noindent [72] W. N. Everitt, ``On the transformation theory of ordinary second-order linear symmetric differential expressions,'' \textit{Czech. Math. J}. \textbf{32}\textit{,} 275 (1982) \underbar{eudml.org/doc/13312}

\noindent [73] R. K. Yadav, A. Khare, and B. P. Mandal, ``The scattering amplitude for one parameter family of shape invariant potentials related to X${}_{m}$ Jacobi polynomials,'' Phys. Lett. B723 (2013), 433-435 arXiv:1303.3669v1\underbar{}

\noindent [74] B. Bagchi, C. Quesne, and R. Roychoudhury  \textit{Pramana J. Phys.} \textbf{73 (}2009) 337-347

\noindent [75] R. K. Yadav, A. Khare, and B. P. Mandal, ``The scattering amplitude for newly found exactly solvable potential,'' \textit{Ann. Phys.} \textbf{331} (2013),\textbf{ }313-316 arXiv:1212.4251v1

\noindent [76] D. Gomez-Ullate, Y. Grandati, and R. Milson, ``Rational extensions of the quantum harmonic oscilator and exceptional Hermite polynomials,'' \textit{J. Phys. A } \textbf{47} (2014), 015203,  

\noindent         27 pages arXiv:1306.5143v1\underbar{}

\noindent [77] N. Rosen and P. M. Morse, ``On the vibrations of polyatomic molecules,'' \textit{Phys. Rev. }\textbf{42} (1932), 210-217

\noindent [78] F. Cooper, A. Khare and U. P. Sukhatme, ``Supersymmetry and quantum mechanics,'' \textit{Phys. Rep.} \textbf{251} (1995), 267-385  arXiv:hep-th/9405029v2

\noindent [79] F. Cooper, A. Khare, and U. P. Sukhatme, \textit{Supersymmetry in Quantum Mechanics }(Denver: World Scientific, 2001)

\noindent [80] J. Ginocchio, ``A class of exactly solvable potentials: I. One-dimensional Schrödinger equation, '' \textit{Ann. Phys.} \textbf{152} (1984), 203-219

\noindent [81] M. Masjedjamei, ``Three finite classes of hypergeometric orthogonal polynomials and their application in functions approximation,'' \textit{Integral Transforms and Special Functions} \textbf{13}  (2002), 169-190

\noindent [82] B. V. Rudyak and B. N. Zakhariev, ``New exactly solvable models for Schrödinger Equation,'' \textit{Inverse Problems }\textbf{3} (1987), 125-133

\noindent 

\noindent

\end{document}